\documentclass[11pt]{article}

\usepackage{amsfonts,amssymb,amsmath,amsthm,epsfig,euscript,epstopdf}

\newcommand{\la}{\lambda}

\newcommand{\eumch}{eu\mathrm{mch}}
\newcommand{\umch}{u\mathrm{mch}}
\newcommand{\des}{\mathrm{des}}
\newcommand{\wdes}{\mathrm{wdes}}

\newcommand{\wris}{\mathrm{wrise}}
\newcommand{\ris}{\mathrm{rise}}

\newcommand{\sgn}{\mathrm{sgn}}

\newcommand{\red}{\mathrm{red}}

\newcommand{\sg}{\sigma}
\newcommand{\PP}{\mathbb{P}}
\newcommand{\lev}{\mathrm{lev}}

\newtheorem{theorem}{Theorem}
\newtheorem{lemma}[theorem]{Lemma}

\newcommand{\fig}[2]{\begin{figure}[ht]
\centerline{\scalebox{.66}{\epsfig{file=#1.eps}}}
\caption{#2}
\label{fig:#1}
\end{figure}}

\title{Generating functions for descents over words 
which avoid a consecutive pattern}

\author{Jeffrey Remmel\\
\small Department of Mathematics\\[-0.8ex]
\small University of California, San Diego\\[-0.8ex]
\small La Jolla, CA 92093-0112. USA\\[-0.8ex]
\small \texttt{jremmel@ucsd.edu}
\and
Luvreet Sangha \\
\small Department of Mathematics\\[-0.8ex]
\small University of California, San Diego\\[-0.8ex]
\small La Jolla, CA 92093-0112. USA\\[-0.8ex]
\small \texttt{lsangha@ucsd.edu}
}

\date{\small Submitted: Date 1;  Accepted: Date 2;
 Published: Date 3.\\
\small MR Subject Classifications: 05A15, 05E05 \\
keywords: compositions, words, pattern match, descent, symmetric functions, ordinary generating function}

\begin{document}

\maketitle

\begin{abstract} In this paper, we extend the reciprocity method 
introduced by Jones and Remmel to study the distributions of 
descents over words which have no $u$-matches for words 
$u$ that have at most one descent. 
\end{abstract}

\section{Introduction}

Let $\PP = \{1, 2, \ldots \}$ denote the set of positive integers and 
for any $k \in \PP$, 
let $[k]= \{1, \ldots, k\}$.  We let $\PP^*$ ($[k]^*$) denote the set of 
all words over the alphabet $\PP$ ($[k]$).  We let $\epsilon$ denote 
the empty word and we say $\epsilon$ has length 0. We let 
$\PP^+ = \PP^*-\{\epsilon\}$ and $[k]^+ = [k]^*-\{\epsilon\}$

If $u =u_1 \ldots 
u_j$ and $v=v_1 \ldots v_i$ are words in $\PP^*$, we let 
$uv = u_1 \ldots u_j v_1 \ldots v_i$ denote the concatenation 
of $u$ and $v$. Suppose that we fix $j \geq 1$. Then for 
any word $w =w_1 \ldots w_n$, we say that 
a word $u =u_1 \ldots u_j$ is a {\em prefix} of 
a word $w$ if 
there is a word $v$ such that $uv =w$, is a {\em suffix}  of $w$ if 
there is a word $v$ such that $vu =w$, and is a {\em factor} of $w$ if 
there are words $f$ and $v$ such that $fuv=w$.

Now suppose that $n \geq  1$ and 
$w = w_1 \ldots w_n \in \PP^n$. Then we let $|w| =n$ denote the length 
of $w$. We let 
$$\begin{array}{ll}
Des(w) = \{i:w_i > w_{i+1}\} & WDes(w) = \{i:w_i \geq w_{i+1}\} \\
Rise(w) = \{i:w_i < w_{i+1}\} & WRise(w) = \{i:w_i \leq w_{i+1}\} \\
\des(w) = |Des(w)| & \wdes(w) = |WDes(w)| \\
\ris(w) = |Rise(w)| & \wris(w) = |WRis(w)|\\
Lev(w) = \{i:w_i =w_{i+1}\} & \lev(w) = |Lev(w)|.
\end{array}$$
We shall refer to elements of $Des(w)$, $WDes(w)$, $Rise(w)$, 
$WRise(w)$, and $Lev(w)$ as descents, weak descents, 
rises, weak rises, and levels of $w$, respectively. 
We let $\overline{z}^w$ denote $z_{w_1} \ldots z_{w_n}$. 
We let $\red(w)$ denote the word that results from $w$ by replacing 
all occurrences of the $i$th smallest letter in $w$ by $i$. For example,
if $w = 44537792$, then $\red(w) =  33425561$.

Let $u=u_1 \ldots u_j \in \PP^j$ and $w=w_1 \ldots w_n \in \PP^n$. 
Then if $\red(u) =u$, 
a $u$-match in $w$ is a factor $v$ of $w$ such that $\red(v) =u$.  
An exact $u$-match in $w$ is a factor $v$ of $w$ such that 
$v =u$.  We let $\umch(w)$ denote the number of $u$-matches 
in $w$ if $\red(u) =u$ and $\eumch(w)$ denote the number of exact $u$-matches 
in $w$. For example, if $w = 31442521337792$ and $u=213$, then 
$w$ has three $u$-matches, namely 314, 425, and 213, but only one 
exact $u$-match.  Thus $\umch(w) =3$ and $\eumch(w) = 1$. 
For any word $w \in \PP^*$ and $i,j \in \PP$, we let $\mathbf{ij}(w)$ denote 
the number of exact matches of $ij$ in $w$.

For any word $u =u_1 \ldots u_j \in [k]^j$ such that $\red(u) =u$, let 
$St^{(\PP)}(u)$ ($St^{([k])}(u)$) 
equal the set of $1 < s \leq j$ such that there exists a word 
$w = w_1 \ldots w_{s+j-1}$ in  $\PP^*$ ($[k]^*$) such that 
$\red(w_1 \ldots w_j) = u$ and $\red(w_s \ldots w_{s+j-1}) =u$. 
That is, $St^{(\PP)}(u)$ ($St^{([k])}(u)$) 
is the set of positions $1<s \leq j$ such that 
there is a word $w$ in $\PP^*$ ($[k]^*$) in which there is a pair of 
overlapping $u$-matches such that the 
first $u$-match starts at position 1 and the second $u$-match starts at position $s$. We say that $u$ is {\em $\PP$-minimal overlapping} 
({\em $[k]$-minimal overlapping}) if $St^{(\PP)}(u) = \{j\}$ 
($St^{([k])}(u) = \{j\}$). Thus $u$ is $\PP$-minimal overlapping 
if any two consecutive $u$-matches can share 
at most one letter which must be the last letter of the first $u$-match 
and the first letter of the second $u$-match. 
We say that $u$ has the {\em $\PP$-weakly decreasing ($\PP$-weakly increasing, 
$\PP$-level)} overlapping property if $s \in St^{(\PP)}(u)$ implies 
that $u_1 \geq u_s$ ($u_1 \leq u_s$, $u_1 =u_s$). We say that $u$ has the {\em $[k]$-weakly decreasing ($[k]$-weakly increasing, 
$[k]$-level)} overlapping property if $s \in St^{([k])}(u)$ implies 
that $u_1 \geq u_s$ ($u_1 \leq u_s$, $u_1 =u_s$). If $k \geq 2$, then we 
say that $u$ is {\em $[k]$-non-overlapping} if $St^{([k])}(u) = \emptyset$. 
For example, suppose that $u = 123234$.  
Then 
\begin{enumerate}
\item $w^{(1)} =123234345$ witnesses that $4 \in St^{(\PP)}(u)$, 
\item $w^{(2)} =1232345456$ witnesses that $5 \in St^{(\PP)}(u)$, and 
\item $w^{(3)} =12323456567$ witnesses that $6 \in St^{(\PP)}(u)$.
\end{enumerate}
It is easy to see that in each case, $w^{(i)}$ uses the smallest 
alphabet possible. Clearly $2$ and $3$ are not in $St^{(\PP)}(u)$ or 
$St^{([k])}(u)$ for any $[k]$.  It thus follows that 
$St^{(\PP)}(u) =St^{([k])}(u) = \{4,5,6\}$ for any $k \geq 7$. 
However, $St^{([4])}(u) = \emptyset$ so that $u$ is 
$[4]$-non-overlapping, $St^{([5])}(u) = \{4\}$, and 
$St^{([6])}(u) = \{4,5\}$.  Note that 
$u$ has the $\PP$-weakly increasing overlapping property and 
the $[k]$-weakly increasing overlapping property for any $k \geq 5$. 
 Next suppose that $v = 345123$.  
Then 
\begin{enumerate}
\item $w^{(4)} =567345123$ witnesses that $4 \in St^{(\PP)}(v)$, 
\item $w^{(5)} =4561345123$ witnesses that $5 \in St^{(\PP)}(v)$, and 
\item $w^{(6)} =34512345123$ witnesses that $6 \in St^{(\PP)}(v)$.
\end{enumerate}
Again it is easy to see that in each case, $w^{(i)}$ uses the smallest 
alphabet possible and $2$ and $3$ are not in $St^{(\PP)}(u)$ or 
$St^{([k])}(v)$ for any $[k]$.  It follows that 
$St^{(\PP)}(v) =St^{([k])}(v) = \{4,5,6\}$ for any $k \geq 7$. 
However, $St^{([5])}(u)= \{6\}$ so that $u$ is $[5]$-minimal overlapping and $St^{([6])}(u) = \{5,6\}$. Note that 
$u$ has the $\PP$-weakly decreasing overlapping property and 
the $[k]$-weakly decreasing overlapping property for any $k \geq 5$ but 
that it also has the $[5]$-weakly increasing overlapping property and 
the $[5]$-level overlapping property.

We can also make similar definitions for exact matchings. That is, for 
$u =u_1 \ldots u_j \in [k]^j$, let 
$ESt^{(\PP)}(u)$ ($ESt^{([k])}(u)$) 
equal the set of $1 < s \leq j$ such that there exists a word 
$w = w_1 \ldots w_{s+j-1}$ in  $\PP^*$ ($[k]^*$) such that 
$w_1 \ldots w_j = u$ and $w_s \ldots w_{s+j-1} =u$. 
That is, $ESt^{(\PP)}(u)$ ($ESt^{([k])}(u)$) 
is the set of positions $1<s \leq j$ such that 
there is a word $w$ in $\PP^*$ ($[k]^*$) in which there is a pair of 
overlapping exact $u$-matches such that the 
first exact $u$-match starts at position 1 and the second exact $u$-match starts at position $s$. We say that $u$ is {\em exact $\PP$-minimal overlapping} 
({\em exact $[k]$-minimal overlapping} if $ESt^{(\PP)}(u) = \{j\}$ 
($ESt^{([k])}(u) = \{j\}$). Thus $u$ is exact $\PP$-minimal overlapping 
if any two consecutive exact $u$-matches can share 
at most one letter which must be the last letter of the first exact $u$-match 
and the first letter of the second exact $u$-match. 
For example $u=131$ is a word that has the exact $\PP$-minimal overlapping 
property.
 We say that $u$ is {\em exact $\PP$-non-overlapping} 
({\em exact $[k]$-non-overlapping} if $ESt^{(\PP)}(u) = \emptyset$ 
($ESt^{([k])}(u) = \emptyset$). 
For example $u=132$ is a word that has the exact $\PP$-non-overlapping 
property.

Let ${\bf z}_k = z_1, \ldots, z_k$ and 
${\bf z}_\infty = z_1,z_2, \ldots $.  Then for any $u \in [k]^j$, we let  
\begin{equation*}
EN^{(k)}_{n,u}(x,{\bf z}_k) =   
\sum_{w \in [k]^n, \eumch(w) =0} x^{\des(w)+1} \overline{z}^w \ \mbox{and} \ 
EN^{(\PP)}_{n,u}(x,{\bf z}_\infty) =  
\sum_{w \in \PP^n, \eumch(w) =0} x^{\des(w)+1} \overline{z}^w. 
\end{equation*}  
Similarly for $u \in [k]^j$ such that $\red(u) =u$, we let 
\begin{equation*}
N^{(k)}_{n,u}(x,{\bf z}_k) = \sum_{w \in [k]^n, \umch(w) =0} x^{\des(w)+1} \overline{z}^w \ \mbox{and} \ 
N^{(\PP)}_{n,u}(x,{\bf z}_\infty) = \sum_{w \in \PP^n, \umch(w) =0} x^{\des(w)+1} \overline{z}^w.
\end{equation*}

The main goal of this paper is to study the generating functions  
\begin{align*}
\mathcal{EN}^{(k)}_{u}(x,{\bf z}_k,t) &= 1 + \sum_{n \geq 1} 
EN^{(k)}_{n,u}(x,{\bf z}) t^n \ \mbox{and} \\  
\mathcal{EN}^{(\PP)}_{u}(x,{\bf z}_\infty,t) &= 1 + \sum_{n \geq 1} 
EN^{(\PP)}_{n,u}(x,{\bf z}) t^n. 
\end{align*}
in the case where $u$ is a word with $\des(u) \leq 1$ and the generating 
functions 
\begin{align*}
\mathcal{N}^{(k)}_{u}(x,{\bf z}_k,t) &= 1 + \sum_{n \geq 1} 
N^{(k)}_{n,u}(x,{\bf z}) t^n  \ \mbox{and}\\  
\mathcal{N}^{(\PP)}_{u}(x,{\bf z}_\infty,t) &= 1 + \sum_{n \geq 1} 
N^{(\PP)}_{n,u}(x,{\bf z}) t^n 
\end{align*}
in the case where $\red(u) =u$ and $\des(u) \leq 1$.

When $k$ and $|u|$ are small, there are well-known recursive methods 
to compute $N^{(k)}_{n,u}(x,{\bf z}_k)$ or $EN^{(k)}_{n,u}(x,{\bf z}_k)$.  
That is, suppose that $|u|=r$.  For any word $v \in [k]^*$, we let 
$\mathcal{B}^{(k)}_v = \{w \in [k]^*: v \ \mbox{is a prefix of} \ w\}$
and
$$\mathcal{N}^{(k)}_{u,v}(x,{\bf z}_k,t) = 
1 + \sum_{n \geq 1} t^n \sum_{w \in 
\mathcal{B}^{(k)}_v \cap [k]^n,\umch(w)=0} x^{\des(w)+1}\overline{z}^w.$$
For example, if $k=3$, $u=123$, and $v =12$, then 
the words in $\mathcal{B}^{(3)}_{12}$ are 
of the form $12$ or $1$ concatenated with 
either a word in $\mathcal{B}^{(3)}_{21}$, $\mathcal{B}^{(3)}_{22}$, 
or $\mathcal{B}^{(3)}_{23}$. Words of the form $1*\mathcal{B}^{(3)}_{23}$ 
cannot contribute to $\mathcal{N}^{(k)}_{u,v}(x,{\bf z}_k,t)$
since they all start with a  $123$-match. It follows that 
$$\mathcal{N}^{(3)}_{u,12}(x,{\bf z}_3,t) =
xz_1z_2t^2 + z_1t\mathcal{N}^{(3)}_{u,21}(x,{\bf z}_3,t) + 
z_1t\mathcal{N}^{(3)}_{u,22}(x,{\bf z}_3,t).$$
In this way, we can show that the functions 
$\mathcal{N}^{(3)}_{u,v}(x,{\bf z}_3,t)$ where $|v|=|u|-1$ 
satisfy simple 
recursions. Bringing the terms that do not involve 
the generating functions to one side, one can rewrite these 
equations in the form 
$$\vec{v} = M 
\begin{pmatrix}\mathcal{N}^{(3)}_{u,11}(x,{\bf z},t)\\
\mathcal{N}^{(3)}_{u,12}(x,{\bf z}_3,t)\\
\mathcal{N}^{(3)}_{u,13}(x,{\bf z}_3,t)\\
\mathcal{N}^{(3)}_{u,21}(x,{\bf z}_3,t)\\
\mathcal{N}^{(3)}_{u,22}(x,{\bf z}_3,t)\\
\mathcal{N}^{(3)}_{u,23}(x,{\bf z}_3,t)\\
\mathcal{N}^{(3)}_{u,31}(x,{\bf z}_3,t)\\
\mathcal{N}^{(3)}_{u,32}(x,{\bf z}_3,t)\\
\mathcal{N}^{(3)}_{u,33}(x,{\bf z}_3,t).
\end{pmatrix}.
$$
Then if one can invert the matrix $M$, one can 
solve for the generating functions 
$\mathcal{N}^{(3)}_{u,ij}(x,{\bf z}_3,t)$ from which 
one can easily recover the desired generating 
function $\mathcal{N}^{(3)}_{u}(x,{\bf z}_3,t)$. 
More details on the method can be found in \cite{HM}. 
The problem with this method   
is that it requires us to invert a 
$k^{|u|-1}$ matrix with multivariable entries 
which is impractical to compute as $k$ and $|u|$ get large.

The method that we will employ is what Jones and Remmel 
\cite{JR1,JR2} call the reciprocal method. The 
basic idea is the following. 
We assume 
that we can write the generating function 
$\mathcal{N}^{(\PP)}_{u}(x,{\bf z}_\infty,t)$ as 
\begin{equation}\label{eq:I1}
\mathcal{N}^{(\PP)}_{u}(x,{\bf z}_\infty,t)) =   \frac{1}{U^{(\PP)}_u(x,{\bf z}_\infty, t)}
\ \mbox{where} \ 
U^{(\PP)}_u(x,{\bf z}_\infty, t) = 1 + \sum_{n\geq 1} U^{(\PP)}_{u,n}(x,{\bf z}_\infty) t^n.
\end{equation}
Thus 
\begin{equation}\label{eq:I3}
U^{(\PP)}_u(x,{\bf z}_\infty, t) = \frac{1}{1+\sum_{n \geq 1} N^{(\PP)}_{u,n}(x,{\bf z}_\infty)t^n}.
\end{equation}
One can then use the homomorphism method to give a combinatorial   
interpretation to the right-hand side of (\ref{eq:I3}) which can 
be used to find a combinatorial 
interpretation for $U^{(\PP)}_u(x,{\bf z}_\infty, t)$. 
The homomorphism method derives generating functions for 
various statistics on permutations and words by 
applying a ring homomorphism defined on the 
ring of symmetric functions \begin{math}\Lambda\end{math}  
in infinitely many variables \begin{math}x_1,x_2, \ldots \end{math} 
to simple symmetric function identities such as 
\begin{equation}\label{conclusion2}
H(t) = 1/E(-t)
\end{equation}
where $H(t)$ and $E(t)$ are the generating functions for the homogeneous and elementary 
symmetric functions given by 
\begin{equation}\label{genfns}
H(t) = \sum_{n\geq 0} h_n t^n = \prod_{i\geq 1} \frac{1}{1-x_it} \ \mbox{and} \ E(t) = \sum_{n\geq 0} e_n t^n = \prod_{i\geq 1} 1+x_it.
\end{equation}
 See, for example, 
\cite{Bec4,Lan2,MenRem1,MenRem2,MenRem3,MRR} or the recent 
book by Mendes and Remmel \cite{MenRembook}.  In our 
case, we define a homomorphism \begin{math}\Theta_u\end{math} on 
\begin{math}\Lambda\end{math} by setting 
\begin{displaymath}\Theta_u(e_n) = (-1)^n N^{(\PP)}_{u,n}(x,{\bf z}_\infty)\end{displaymath}
Then 
\begin{displaymath} \Theta_u(E(-t)) = \sum_{n\geq 0} N^{(\PP)}_{u,n}(x,{\bf z}_\infty) t^n = \frac{1}{U^{(\PP)}_u(x,{\bf z}_\infty,t)}.\end{displaymath}
Hence 
$${U^{(\PP)}_u(x,{\bf z}_\infty,t)} = \frac{1}{\Theta_u(E(-t))} = \Theta_u(H(t))$$
which implies that 
\begin{equation}\label{eq:combhn}
\Theta_u(h_n) = U^{(\PP)}_{u,n}(x,{\bf z}_\infty).
\end{equation}
Thus if we can compute $\Theta_u(h_n)$ for all $n \geq 1$, then we can 
compute the polynomials $U^{(\PP)}_{u,n}(x,{\bf z}_\infty)$ and the generating function 
$U^{(\PP)}_u(x,{\bf z}_\infty,t)$ which in turn allows us to compute 
the generating function $\mathcal{N}^{(\PP)}_{u}(x,{\bf z}_\infty,t)$. 
The same method can be applied to find combinatorial 
interpretations for $U^{(k)}_u(x,{\bf z}_k,t)$, 
$EU^{(\PP)}_u(x,{\bf z}_\infty,t)$, and $EU^{(k)}_u(x,{\bf z}_k,t)$ 
where 
\begin{align*}
\mathcal{N}^{(k)}_{u}(x,{\bf z}_k,t)) &=   
\frac{1}{U^{(k)}_u(x,{\bf z}_k, t)}, \\
\mathcal{EN}^{(\PP)}_{u}(x,{\bf z}_\infty,t)) &=   \frac{1}{EU^{(\PP)}_u(x,{\bf z}_\infty, t)}, \ \mbox{and} \\
\mathcal{EN}^{(k)}_{u}(x,{\bf z}_k,t)) &=   
\frac{1}{EU^{(k)}_u(x,{\bf z}_k, t)}.
\end{align*}

The final steps of the reciprocity method 
 that we employ will be different from the 
ones used by Jones and Remmel \cite{JR1,JR2} for permutations. 
For the generating 
function for permutations that they studied, Jones and Remmel  used 
the combinatorial interpretation that arose 
from the analogue of $\Theta_u(h_n)$ to obtain simple recursions 
satisfied by their analogue of $U^{(\PP)}_{u,n}(x,{\bf z}_\infty)$.  
In our case, we shall use the combinatorial interpretation 
of $\Theta_u(h_n)$ that comes out of the homomorphism method 
plus a map which we call the ``collapse map'' to 
show that we can obtain a closed expression for the generating functions 
$U^{(\PP)}_u(x,{\bf z}_\infty,t)$ or $U^{([k])}_u(x,{\bf z}_k,t)$ 
by an appropriate substitution in certain other 
generating functions for words. 

The generating functions that we will substitute in 
will depend on the relative order of $u_1$ and $u_j$ where 
$u =u_1 \ldots u_j$.  In each case our generating 
function will be over the variables $x_{ij}$ where $i,j \in \PP$, 
the variables $z_i$ where $i \in \PP$, and $t$ which 
we denote as $(\mathbf{x}_\infty,\mathbf{z}_\infty,t)$. 
In the case where $u_1 > u_j$, our final expression 
for our desired generating functions $U^{(\PP)}_u(x,{\bf z}_\infty,t)$ or $U^{([k])}_u(x,{\bf z}_k,t)$ will be a substitution into 
the generating function 
$$\mathcal{D}^{\PP}(\mathbf{x}_\infty,\mathbf{z}_\infty,t) = 
\sum_{w \in \PP^*} t^{|w|}\overline{z}^w \prod_{i < j} 
x_{ji}^{\mathbf{ji}(w)}.$$

In the case where $u_1 < u_j$ and $u$ has 
the $\PP$-weakly increasing overlapping property 
($[k]$-weakly increasing overlapping property), our final expression 
for our desired generating function $U^{(\PP)}_u(x,{\bf z}_\infty,t)$ 
($U^{([k])}_u(x,{\bf z}_k,t)$) will be a substitution into 
the generating function 
$$\mathcal{R}^{\PP}(\mathbf{x}_\infty,\mathbf{z}_\infty,t) = 
\sum_{w \in w_1 \leq w_2 \leq \cdots \leq w_n \in \PP^*} t^{|w|}\overline{z}^w \prod_{i < j} 
x_{ij}^{\mathbf{ij}(w)}.$$

In the case were $u_1 =u_j$ and $u$ has 
the $\PP$-level overlapping property 
($[k]$-level overlapping property),  our final expression 
for our desired generating function $U^{(\PP)}_u(x,{\bf z}_\infty,t)$ 
($U^{([k])}_u(x,{\bf z}_k,t)$) will be a substitution into 
the generating function 
$$\mathcal{L}^{\PP}(\mathbf{x}_\infty,\mathbf{z}_\infty,t) = 
\sum_{w = w_1 \leq w_2 \leq \cdots w_n \in \PP^*} t^{|w|}\overline{z}^w 
\prod_{i} x_{ii}^{\mathbf{ii}(w)}.$$
If $u$ does not have the $\PP$-level overlapping property 
($[k]$-level overlapping property), it will still be the case 
that $u$ has the $\PP$-weakly decreasing  overlapping property 
($[k]$-weakly decreasing overlapping property). In such a  
case, our final expression 
for our desired generating function $U^{(\PP)}_u(x,{\bf z}_\infty,t)$ 
($U^{([k])}_u(x,{\bf z}_\infty,t)$) will be a substitution into 
the generating function 
$$\mathcal{WD}^{\PP}(\mathbf{x}_\infty,\mathbf{z}_\infty,t) = 
\sum_{w \in \PP^*} t^{|w|}\overline{z}^w \prod_{i \leq j} 
x_{ji}^{\mathbf{ji}(w)}.$$

We will prove the following theorems for the 
generating functions 
$\mathcal{D}^{\PP}(\mathbf{x}_\infty,\mathbf{z}_\infty,t)$, 
$\mathcal{L}^{\PP}(\mathbf{x}_\infty,\mathbf{z}_\infty,t)$, and 
$\mathcal{R}^{\PP}(\mathbf{x}_\infty,\mathbf{z}_\infty,t)$. 
Given a set 
$S \subseteq \PP$, we let 
\begin{equation}
DXZ(S) = \begin{cases}
z_j & \mbox{if} \ S = \{j\}, \ \mbox{and} \\
z_{j_1} \cdots z_{j_k} \prod_{i=1}^{k-1} (x_{j_{i+1}j_i}-1) & 
\mbox{if} \ S =\{j_1 < \cdots < j_k\} \ \mbox{where} \ k \geq 2.
\end{cases}
\end{equation}
and 
\begin{equation}
RXZ(S) = \begin{cases} 
\frac{z_j}{1-z_jt} & \mbox{if} \ S = \{j\}, \ \mbox{and} \\
\left(\prod_{i=1}^k \frac{z_{j_i}}{1-z_{j_i}t}\right) 
\prod_{i=1}^{k-1} x_{j_{i}j_{i+1}} & 
\mbox{if} \ S =\{j_1 < \cdots < j_k\} \ \mbox{where} \ k \geq 2.
\end{cases}
\end{equation}

Let $WD\PP^*$ ($WD[k]^*$) denote the set of all weakly decreasing words 
in $\PP^*$ ($[k]^*$). 
Given a nonempty word $v$ in 
$WD\PP^*$, we let 
\begin{equation}
WDXZ(v) = \begin{cases}
z_j & \mbox{if} \ v = j, \ \mbox{and} \\
z_{j_1} \cdots z_{j_k} \prod_{i=1}^{k-1} (x_{j_{i}j_{i+1}}-1) & 
\mbox{if} \ v =j_1 \geq \cdots \geq  j_k \ \mbox{where} \ k \geq 2.
\end{cases}
\end{equation}

\begin{theorem}\label{thm:labdes}
\begin{equation}\label{eq:labdes}
\mathcal{D}^\PP ({\bf x}_\infty,{\bf z}_\infty,t) = 
\frac{1}{1 - \sum_{n \geq 1} t^n \sum_{S \subseteq \PP, |S| =n} 
DXZ(S)}.
\end{equation}
\end{theorem}

\begin{theorem}\label{thm:labwdes}
\begin{equation}\label{eq:labwdes}
\mathcal{WD}^\PP ({\bf x}_\infty,{\bf z}_\infty,t) = 
\frac{1}{1 - \sum_{n \geq 1} t^n \sum_{v \in WD\PP^*, |v| =n} 
WDXZ(w)}.
\end{equation}
\end{theorem}

\begin{theorem}\label{thm:labris}
\begin{equation}\label{eq:labris}
\mathcal{R}^\PP ({\bf x}_\infty,{\bf z}_\infty,t) = 
1+ \sum_{n \geq 1} t^n \sum_{S \subseteq \PP, |S| =n} RXZ(S).
\end{equation}
\end{theorem}
\begin{theorem}\label{thm:lablev}
\begin{equation}\label{eq:lablev}
\mathcal{L}^\PP ({\bf x}_\infty,{\bf z}_\infty,t) = 
\prod_{i \geq 1} \left( 1+ \frac{z_it}{1-x_{ii}z_it}\right)
\end{equation}
\end{theorem}

The main advantage of our approach is that we obtain 
a uniform way to find expressions for the generating functions 
$U^{(\PP)}_u(x,{\bf z}_\infty,t)$, 
$U^{(k)}_u(x,{\bf z}_k,t)$, 
$EU^{(\PP)}_u(x,{\bf z}_\infty,t)$, and $EU^{(k)}_u(x,{\bf z}_k,t)$ 
which are independent of the length of $u$ as long as 
$\des(u) = 1$ and $u$ satisfies the appropriate overlapping 
conditions.  In fact our general methods can 
be applied even in cases where $\des(u) > 1$. However in 
such cases the combinatorial interpretation of $\Theta_u(h_n)$ 
that comes out the homomorphism method is significantly 
more complicated so that we will not pursue such results 
in this paper. 

The outline of this paper is as follows. 
In Section 2, we shall review the basic background on symmetric 
functions that one will need for the paper. 
In Section 3, we shall describe the use 
of the reciprocal method to obtain combinatorial 
interpretations for $U^{(\PP)}_{u,n}(x,{\bf z}_\infty,t)$, 
$U^{(k)}_{u,n}(x,{\bf z}_k,t)$, 
$EU^{(\PP)}_{u,n}(x,{\bf z}_\infty,t)$, and $EU^{(k)}_{u,n}(x,{\bf z}_k,t)$.
In section 4, we shall show how to use Theorem \ref{thm:labdes}  
 to find expressions for $U^{(\PP)}_u(x,{\bf z}_\infty,t)$, 
$U^{(k)}_u(x,{\bf z}_k,t)$, 
$EU^{(\PP)}_u(x,{\bf z}_\infty,t)$, and $EU^{(k)}_u(x,{\bf z}_k,t)$ 
in the case where $u =u_1 \ldots u_j$, $u_1 > u_j$,  and $\des(u)=1$.
In 
section 5, we shall show how to use Theorem \ref{thm:labris}  
 to find expressions for $U^{(\PP)}_u(x,{\bf z}_\infty,t)$, 
$U^{(k)}_u(x,{\bf z}_k,t)$, 
$EU^{(\PP)}_u(x,{\bf z}_\infty,t)$, and $EU^{(k)}_u(x,{\bf z}_k,t)$ 
in the case where $u =u_1 \ldots u_j$, $u_1 < u_j$, $\des(u) = 1$, 
and $u$ has the $\PP$-weakly increasing overlapping property or 
$[k]$-weakly increasing overlapping property.
In section 6, we shall show how to use Theorems  
\ref{thm:labwdes} and  \ref{thm:lablev} 
to find expressions for $U^{(\PP)}_u(x,{\bf z}_\infty,t)$, 
$U^{(k)}_u(x,{\bf z}_k,t)$, 
$EU^{(\PP)}_u(x,{\bf z}_\infty,t)$, and $EU^{(k)}_u(x,{\bf z}_k,t)$ 
in the case where $u =u_1 \ldots u_j$, $u_1 = u_j$, $\des(u)\leq 1$, and $u$ has the $\PP$-level overlapping property or 
$[k]$-level overlapping property.
In 
section 6, we shall prove Theorems \ref{thm:labdes},   \ref{thm:labwdes}, 
\ref{thm:labris}, and  \ref{thm:lablev}. 
Finally, in Section 7, we shall discuss some further extensions 
of our methods. For example, we will discuss how we 
can replace the statistic $\des(w)$ in our formulas 
by $\wdes(w)$ or $\lev(w)$ and we will discuss how we 
can extend our methods to handle cases where 
$u$ has more than one descent. 

\section{Symmetric Functions}
\label{section:sf}

In this section we give the necessary background on symmetric
functions needed for our proofs. 
We shall consider the ring of symmetric functions, $\Lambda$, over 
infinitely many variables $x_1, x_2, \ldots $.  The homogeneous 
symmetric functions, $h_n \in \Lambda$, and elementary symmetric functions, 
$e_n \in \Lambda$, are defined by the generating 
functions 
\begin{equation*}
H(t) = \sum_{n \geq 0} h_n t^n = \prod_{i=1}^\infty \frac{1}{1-x_it} 
\ \mbox{and} \ 
E(t) = \sum_{n \geq 0} e_n t^n = \prod_{i=1}^\infty (1+x_it).
\end{equation*} 
The $n$-th power symmetric function, $p_n \in \Lambda$, is defined as  
$\displaystyle p_n = \sum_{i=1}^\infty x_i^n$. 

Let $\la = (\la_1,\dots,\la_{\ell})$ be an integer partition; that is,
$\la$ is a finite sequence of weakly increasing non-negative
integers.  Let $\ell(\la)$ denote the number of nonzero integers in
$\la$. If the sum of these integers is $n$, we say that $\la$ is a
partition of $n$ and write $\la \vdash n$.  For any partition $\la =
(\la_1,\dots,\la_\ell)$, define $h_\la = h_{\la_1} \cdots
h_{\la_\ell}$, $e_\la = e_{\la_1} \cdots
e_{\la_\ell}$, and $p_\la = p_{\la_1} \cdots
p_{\la_\ell}$. The well-known fundamental theorem of symmetric
functions, see \cite{Mac}, says that 
$\{e_\la : \la \vdash n\}$ is a
basis for $\Lambda_n$, the space of symmetric functions 
which are homogeneous of degree $n$. Equivalently, the fundamental 
theorem of symmetric functions states that $\{e_0,e_1, \ldots \}$ 
is an algebraically independent set of generators for 
the ring $\Lambda$. It follows that one can completely specify a  
ring homomorphism $\Gamma: \Lambda \rightarrow R$ from 
$\Lambda$ into a ring $R$ by giving the values of $\Gamma(e_n)$ for 
$n \geq 0$.   

Next we give combinatorial interpretations 
to the expansion of $h_\mu$ in terms of the elementary symmetric 
functions.  Given partitions $\lambda =(\la_1, \ldots, \la_\ell)\vdash n$ and $\mu\vdash n$,
a {\em $\lambda$-brick tabloid of shape $\mu$} is a filling 
of the Ferrers diagram of shape $\mu$ with bricks of size $\lambda_1, \ldots, \lambda_\ell$ such that 
each brick lies in single row and no two bricks overlap. For 
example, Figure  \ref{fig:labrick} shows all the 
$\la$-brick tabloids of shape $\mu$ where $\la = (1,1,2,2)$ and $\mu = (2,4)$.

\fig{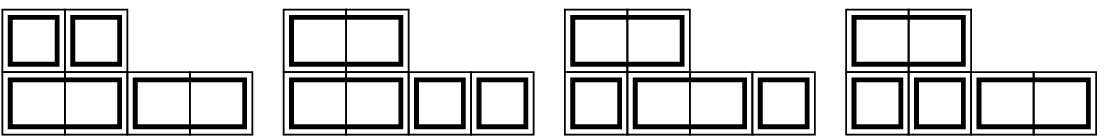}{The four $(1,1,2,2)$-brick tabloids of shape $(2,4)$.}

If $T$ is a brick tabloid of shape $(n)$ such that the 
lengths of the bricks, reading from left to right, are 
$b_1, \ldots, b_\ell$, then we shall write 
$T=(b_1, \ldots, b_\ell)$. For example, 
the brick tabloid $T =(2,3,1,4,2)$ is pictured in 
Figure \ref{fig:btab}. 

\fig{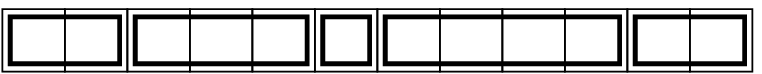}{The brick tabloid $T =(2,3,1,4,2)$.}

Let $\mathcal{B}_{\la,\mu}$ denote the set of all $\la$-brick tabloids of shape $\mu$ and let
$B_{\la,\mu} =|\mathcal{B}_{\la,\mu}|$.  E\u{g}ecio\u{g}lu and Remmel proved in
\cite{Eg1} that
\begin{equation}
\label{htoe}
h_\mu = \sum_{\la \vdash n} (-1)^{n - \ell(\la)} B_{\la,\mu} e_\la.
\end{equation}

\section{The reciprocal method}

In this section, we shall apply the reciprocal method to give 
combinatorial interpretations to $U^{(\PP)}_u(x,{\bf z}_\infty,t)$, 
$U^{(k)}_u(x,{\bf z}_k,t)$, 
$EU^{(\PP)}_u(x,{\bf z}_\infty,t)$, and $EU^{(k)}_u(x,{\bf z}_k,t)$. 

Fix a word $u$ such that $\des(u) \leq 1$. We will start out by considering 
$U^{(\PP)}_u(x,{\bf z}_\infty,t)$ as the other cases are similar.
Recall that 
\begin{equation}
U^{(\PP)}_u(x,{\bf z}_\infty,t) 
= \frac{1}{1+\sum_{n \geq 1} N^{(\PP)}_{u,n}(x,{\bf z}_\infty)t^n}.
\end{equation}
Thus if we let 
$\Theta_u(e_n) = (-1)^n N^{(\PP)}_{u,n}(x,{\bf z}_\infty)$ for 
$n \geq 1$ and $\Theta_u(e_0) =1$, we see that 
\begin{eqnarray*}
\Theta_u(H(t)) &=& 1 + \sum_{n \geq 1} \Theta_u(h_n) \\
&=& \Theta_u\left( \frac{1}{E(-t))}\right) 
=\frac{1}{1+ \sum_{n \geq 1} (-1)^n \Theta_u(e_n)} \\
&=& \frac{1}{1 + \sum_{n \geq 1} N^{(\PP)}_{u,n}(x,{\bf z}_\infty)t^n} 
= U^{(\PP)}_u(x,{\bf z}_\infty,t).
\end{eqnarray*}
Thus it follows that 
$\Theta_u(h_n) = U^{(\PP)}_{u,n}(x,{\bf z}_\infty)$.

By (\ref{htoe}), we have  that 
\begin{eqnarray}\label{eq:basic1}
\Theta_{u}(h_n) &=&  \sum_{\la \vdash n} (-1)^{n-\ell(\la)} 
B_{\la,n}~ \Theta_{u}(e_\la) \nonumber \\
&=&  \sum_{\la \vdash n} (-1)^{n-\ell(\la)} \sum_{(b_1, \ldots, 
b_{\ell(\la)}) \in \mathcal{B}_{\la,n}} \prod_{i=1}^{\ell(\la)}  
(-1)^{b_i}  N^{(\PP)}_{u,b_i}(x,{\bf z}_\infty) \nonumber \\
&=& \sum_{\la \vdash n} (-1)^{\ell(\la)} 
\sum_{(b_1, \ldots, b_{\ell(\la)}) \in \mathcal{B}_{\la,n}} 
\prod_{i=1}^{\ell(\la)} N^{(\PP)}_{u,b_i}(x,{\bf z}_\infty). 
\end{eqnarray}

Our next goal is to give a combinatorial interpretation to the 
right-hand side of (\ref{eq:basic1}). Fix a partition 
$\la$ of $n$ and a $\lambda$-brick tabloid 
$B=(b_1, \ldots, b_{\ell(\la)})$. We will 
interpret $\prod_{i=1}^{\ell(\la)}  N^{(\PP)}_{u,b_i}(x,{\bf z}_\infty) $ as 
the number of ways of picking words 
$(w^{(1)}, \ldots, w^{(\ell(\la))})$ such that for each 
$i$,  $w^{(i)} \in \PP^{b_i}$ is a word such 
that $\umch(w) = 0$ and assigning 
a weight to this $\ell(\la)$-tuple to be  
$\prod_{i=1}^{\ell(\la)}  x^{\des(w^{(i)})+1} \overline{z}^{w^{(i)}}$.

\begin{figure}[htbp]
  \begin{center}
    \includegraphics[width=0.6\textwidth]{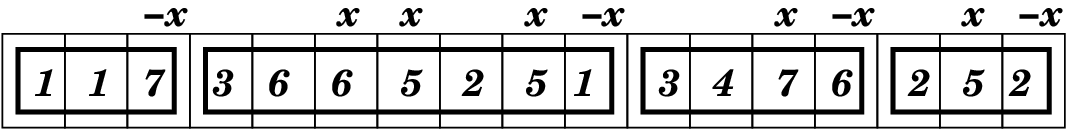}
    \caption{The construction of a filled-labeled-brick tabloid.}
    \label{fig:wFLBTab}
  \end{center}
\end{figure}

We can then use the pair $\langle B, (w^{(1)}, \ldots, w^{(\ell(\la))})\rangle$ to 
construct a filled-labeled-brick tabloid 
$O_{\langle B, (w^{(1)}, \ldots, w^{(\ell(\la))}\rangle}$ as 
follows. First for each brick $b_i$, we place 
the  word $w^{(i)}$ in the cells of the brick, reading from left to right.  
Then we label each cell 
of $b_i$ that starts a descent of $w^{(i)}$ with a $x$ and 
we also label the last cell of $b_i$ with $x$. This accounts 
for the factor $x^{\des(w^{(i)})+1}$. 
Finally, we use  
the factor $(-1)^{\ell(\la)}$ to change the label of 
the last cell of each brick from $x$ to $-x$. 
For example, suppose $n =17$, $u = 312$,  
$B = (3,7,4,3)$  
$w^{(1)} = 1~1~7$, 
$w^{(2)} = 3~6~6~5~2~5~1$, $w^{(3)} = 3~4~7~6$, 
and $w^{(4)} = 2~5~2$. Then 
we have pictured the filled-labeled brick tabloid $O_{\langle B, 
(w^{(1)}, \ldots, w^{(4)})\rangle}$ constructed 
from the pair  $\langle B, (w^{(1)}, \ldots, w^{(4)})\rangle$ in 
Figure \ref{fig:wFLBTab}.

Clearly, we can recover the 
pair $\langle B, (w^{(1)}, \ldots, w^{(\ell(\la))})\rangle$ and the labels 
on the cells from 
$B$ and the word  $w$ which is obtained 
by reading the elements in the cells of 
$O_{\langle B, (w^{(1)}, \ldots, w^{(\ell(\la))})\rangle}$ from left to right. 
Thus we shall specify the filled-labeled-brick tabloid 
$O_{\langle B, (w^{(1)}, \ldots, w^{(\ell(\la))})\rangle}$ by $(B,w)$. 
We let $\mathcal{O}^{(\PP)}_{u,n}$ denote the set of 
all filled-labeled-brick tabloids constructed in this way. 
That is, $\mathcal{O}^{(\PP)}_{u,n}$ consists of all pairs 
$O=(B,w)$ where 
\begin{enumerate}
\item $B=(b_1, \ldots, b_{\ell(\la)})$ is brick tabloid 
of shape $(n)$, 
\item $w =w_1 \ldots w_n \in \PP^n$ such that there is 
no $u$-match of $\sg$ which is entirely contained in a single 
brick of $B$, and 
\item if there is a cell $c$ such that a brick $b_i$ contains both cells $c$ and $c+1$
and $w_c > w_{c+1}$, then cell $c$ is labeled with a $x$ and the last cell of any brick is labeled with $-x$.
\end{enumerate}
The sign of $O$, $\sgn(O)$, is $(-1)^{\ell(\la)}$ and the weight 
of $O$, $wt(O)$,  is $x^{\ell(\la) +\mathrm{intdes}(\sg)}\overline{z}^w$ where 
$\mathrm{intdes}(w)$ denotes the number of $i$ such that 
$w_i > w_{i+1}$ and $w_i$ and $w_{i+1}$ lie in the same brick. 
We shall refer to such $i$ as an {\em internal descent} of $O$. 
Note that the labels on $O$ are completely determined 
by the underlying brick tabloid $B = (b_1, \ldots ,b_{\ell(\la)})$ and 
the underlying word $w$.  Thus the filled-labeled-brick 
tabloid $O$ pictured in Figure \ref{fig:wFLBTab} equals
 $((3,7,4,3),1~1~7~3~6~6~5~2~5~1~3~4~7~6~2~5~2)$. 

It follows that 
\begin{equation}\label{eq:basic2}
\Theta_{u}(h_n) = \sum_{O \in \mathcal{O}^{(\PP)}_{u,n}} sgn(O) wt(O).
\end{equation}

\begin{figure}[htbp]
  \begin{center}
    \includegraphics[width=0.6\textwidth]{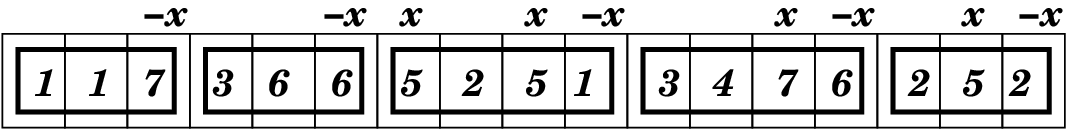}
    \caption{$I_{u}(O)$ for $O$ in Figure \ref{fig:wFLBTab}.}
    \label{fig:wFLBTab2}
  \end{center}
\end{figure}

Next we define a weight-preserving, sign-reversing involution 
$I_{u}$ on \begin{math}\mathcal{O}^{(\PP)}_{u,n}\end{math}. Given an element 
\begin{math}O =(B,w) \in \mathcal{O}^{(\PP)}_{u,n}\end{math} 
where $B=(b_1, \ldots, b_k)$ and $w = w_1 \ldots w_n$, 
scan the cells of \begin{math}O\end{math} from left to right 
looking for the first 
cell \begin{math}c\end{math} such that either 
\begin{enumerate}
\item[(i)] $c$ is labeled with a $x$ or 

\item[(ii)] $c$ is a cell at the end of a brick \begin{math}b_i\end{math}, 
$w_c > w_{c+1}$,  
and there is no $u$-match of $w$ that lies entirely 
in the cells of bricks \begin{math}b_i\end{math} and \begin{math}b_{i+1}\end{math}. 
\end{enumerate} 
In case (i), if \begin{math}c\end{math} is a cell in brick \begin{math}b_j\end{math}, then we split 
\begin{math}b_j\end{math} into two bricks \begin{math}b_j^\prime\end{math} and \begin{math}b_j^{\prime \prime}\end{math} where 
\begin{math}b_j^\prime\end{math} contains all the cells of \begin{math}b_j\end{math} up to and including 
cell \begin{math}c\end{math} and \begin{math}b_j^{\prime \prime}\end{math} consists of the remaining cells 
of \begin{math}b_j\end{math} and we change the label on cell \begin{math}c\end{math} from \begin{math}x\end{math} to \begin{math}-x\end{math}. 
In case (ii), we combine the two bricks \begin{math}b_i\end{math} and \begin{math}b_{i+1}\end{math} into 
a single brick \begin{math}b\end{math} and change the label on cell \begin{math}c\end{math} from \begin{math}-x\end{math} to \begin{math}x\end{math}.  If neither case (i) nor case (ii) applies, then we define $I_u(O) = O$. 
For example, consider the element 
\begin{math}O \in \mathcal{O}^{(\PP)}_{312,17}\end{math} pictured in Figure 
\ref{fig:wFLBTab}.  Note that even 
though the number in the last cell of brick 1 is greater than the number in the first cell of brick 2, we can not combine these two bricks because 
\begin{math}7~3~6\end{math} would be a 312-match. Thus the first 
place that we can apply the involution is on cell 6 which is labeled 
with an \begin{math}x\end{math} so  that \begin{math}I_{u}(O)\end{math} is the object pictured in Figure \ref{fig:wFLBTab2}. 

We claim that whenever $u$ is a word such that
 $\red(u) =u$ and $\des(u) \leq 1$, 
\begin{math}I_{u}\end{math} is an involution, i.e. \begin{math}I_{u}^2\end{math} is the identity. 
First we consider the case where $\des(u) =1$. 
Now suppose that we are in case (i) where we split a brick \begin{math}b_j\end{math} at 
cell \begin{math}c\end{math} which is labeled with a \begin{math}x\end{math}.  
In that case, we let 
\begin{math}a\end{math} be the number in cell \begin{math}c\end{math} and \begin{math}a^\prime\end{math} be the number
in cell \begin{math}c+1\end{math} which must also be in brick \begin{math}b_{j}\end{math}. It must be the case 
that there is no cell labeled \begin{math}x\end{math} before cell \begin{math}c\end{math} since otherwise 
we would not use cell \begin{math}c\end{math} to define the involution. 
However, we  have to consider the possibility 
that when we split \begin{math}b_j\end{math} into \begin{math}b_j^\prime\end{math} and \begin{math}b_j^{\prime \prime}\end{math}, we might then be able to combine the brick \begin{math}b_{j-1}\end{math} with  
\begin{math}b_j^\prime\end{math} because the number in that last cell of 
\begin{math}b_{j-1}\end{math} is greater than the number in the first cell of 
\begin{math}b_j^\prime\end{math} and there is no $u$-match in the cells of 
\begin{math}b_{j-1}\end{math} and  \begin{math}b_j^\prime\end{math}.  
Since we always take an action on the left most cell possible 
when defining $I_u(O)$, we know that we cannot combine  \begin{math}b_{j-1}\end{math} and \begin{math}b_j\end{math} so that 
there must be a $u$-match in the cells of \begin{math}b_{j-1}\end{math} and \begin{math}b_j\end{math}. Clearly, 
that $u$-match must have involved the number \begin{math}a^\prime\end{math} and 
the number in cell $d$ which is the last cell in brick $b_{j-1}$. But that is impossible because then there would be two descents among the numbers
between cell $d$ and cell $c+1$ which would violate our assumption 
that $u$ has only one descent. Thus whenever we apply case (i) to 
define $I_u(O)$, 
the first action that we can take is to combine bricks 
\begin{math}b_j^\prime\end{math} and \begin{math}b_j^{\prime \prime}\end{math} so that \begin{math}I_{u}^2(O) = O\end{math}. 

If we are in case (ii), then again we can assume that there 
are no cells labeled \begin{math}x\end{math} that occur before cell \begin{math}c\end{math}.  When we combine 
bricks \begin{math}b_i\end{math} and \begin{math}b_{i+1}\end{math}, then we will label cell \begin{math}c\end{math} with a \begin{math}x\end{math}. 
It is clear that combining the cells of \begin{math}b_i\end{math} and \begin{math}b_{i+1}\end{math} cannot 
help us combine the resulting brick \begin{math}b\end{math} with an earlier brick since 
it will be harder to have no $u$-matches with the larger brick \begin{math}b\end{math}. 
Thus the first place cell \begin{math}c\end{math} where we can apply the involution to $I_u(O)$ 
will again be cell \begin{math}c\end{math} which is now labeled with a \begin{math}x\end{math} so that \begin{math}I_{u}^2(O) =O\end{math} if we are in case (ii).

The case where $\des(u) =0$ is even easier.  Suppose that $a$ is number in the 
the last cell of $b_j$ and $a'$ is the number in the first cell of $b_{j+1}$ and $a > a'$. 
Then there can be no $u$-match of $w$ that is contained in the cells of $b_j$ and $b_{j+1}$ because by our 
definitions there is no $u$-match in the cells of $b_j$ and there is no $u$-match in the cells of $b_{j+1}$ so that the only possible $u$-match 
in the cells of $b_j$ and $b_{j+1}$ would have to involve $a$ and $a'$ 
which is impossible if $\des(u) =0$. It easily follows that  
we will apply the involution to the first possible cell $c$ which is labeled with either $x$ or $-x$ and what ever action we take at cell $c$ to 
create $I_u(O)$, we will come back to cell $c$ to undo that 
action to define $I^2(O)$.

Our definitions ensure that if 
\begin{math}I_{u}(O) \neq O\end{math}, 
then 
\begin{math}sgn(O)wt(O) = - sgn(I_{u}(O))wt(I_{u}(O))\end{math}. Hence, 
if we let $\mathcal{IO}^{(\PP)}_{u,n}$ denote 
set set all $O =(B,w) \in \mathcal{O}^{(\PP)}_{u,n}$ such that 
$I_u(O) =O$, then 
\begin{equation}\label{eq:basic3} 
\Theta_u(h_n) =  \sum_{O \in \mathcal{O}^{(\PP)}_{u,n}} sgn(O) wt(O) = 
\sum_{O \in \mathcal{IO}^{(\PP)}_{u,n}} sgn(O) wt(O). 
\end{equation}
Thus we must examine the fixed points of 
\begin{math}I_{u}\end{math}. So assume 
that \begin{math}(B,w)\end{math} is a fixed point of \begin{math}I_{u}\end{math}.

There are two cases to consider. \\
\ \\
{\bf Case 1.}  $\des(u) =0$.  \\

  Suppose that $(B,w) \in \mathcal{IO}_{u,n}$
where $B=(b_1, \ldots, b_k)$ and $w= w_1 \ldots w_n$. 
There can be no cell $c$ which is labeled with $x$ in $(B,w)$ since 
we could use such a cell to define $I_u$ which would violate our 
assumption that $(B,w)$ is a fixed point of $I_u$. Similarly there 
can be no cell $c$ which is at the end of a brick $b_j$ 
such that $w_c > w_{c+1}$ since 
again we could use such a cell to  define $I_u(O)$.
This means that $w$ must be weakly increasing within any brick and if 
$c$ is a cell at the  end of brick $b_j$ which is followed by another brick 
$b_{j+1}$, then $w_c \leq w_{c+1}$. Thus 
$(B,w)$ is a fixed point if and only if $w$ is a weakly increasing word such that $w$ has 
no $u$-match that lies entirely within one of the brick of $B$. 
If $B$ has $k$ bricks,
then then weight of $(B,w)$ is just $(-x)^k \overline{z}^w$. 
We let $\mathcal{WIO}_{u,n} = \{(B,w) \in \mathcal{IO}^{(\PP)}_{u,n}:w_1 \leq w_2 \leq \cdots \leq w_n\}$ 
denote the set of elements of $\mathcal{IO}^{(\PP)}_{u,n}$ where $w$ is weakly increasing. 
Then we have the following lemma. Let $\mathbb{Q}(x,{\bf z}_\infty)$ be the set of rational functions in 
the variables $x$ and ${\bf z}_\infty$ over the rationals $\mathbb{Q}$. 

\begin{lemma}\label{lem:uinc} 
Suppose that $u$ is a word in $\PP^+$ such that $\red(u) =u$ and $\des(u) =0$. 
Let $\Theta_u:\Lambda \rightarrow 
\mathbb{Q}(x,{\bf z}_\infty)$ be the  ring homomorphism defined by setting   
$\Theta_{u}(e_0) =1$ and  
$\Theta_{u}(e_n) = (-1)^n N^{(\PP)}_{u,n}(x,{\bf z}_\infty)$ for $n \geq 1$.
Then 
\begin{equation}
U^{(\PP)}_{u,n}(x,{\bf z}_\infty) = \Theta_u(h_n) = \sum_{((b_1, \ldots, b_k),w) \in \mathcal{WIO}_{u,n}} (-x)^k \overline{z}^w.
\end{equation}
\end{lemma}
\ \\
{\bf Case 2.} $\des(u) = 1$. \\

First it is easy to see that there can be no cells which are labeled with 
\begin{math}x\end{math} so that numbers in each brick of \begin{math}O\end{math} must be weakly increasing. 
Second we cannot combine two consecutive bricks \begin{math}b_i\end{math} and \begin{math}b_{i+1}\end{math} in 
\begin{math}O\end{math} which means that either  there is an weak increase between the bricks 
\begin{math}b_i\end{math} and \begin{math}b_{i+1}\end{math} or there is a decrease between the bricks 
\begin{math}b_i\end{math} and \begin{math}b_{i+1}\end{math}, but there is a 
$u$-match in the cells of the bricks 
\begin{math}b_i\end{math} and \begin{math}b_{i+1}\end{math}. 
Thus we have proved the following. 

\begin{lemma} \label{lem:udes1}
Suppose that $u \in \PP^+$, $\red(u) =u$, and $\des(u)=1$.  Let $\Theta_u:\Lambda \rightarrow 
\mathbb{Q}(x,{\bf z}_\infty)$ be the  ring homomorphism defined by setting   
$\Theta_{u}(e_0) =1$ and  
$\Theta_{u}(e_n) = (-1)^n N^{(\PP)}_{u,n}(x,{\bf z}_\infty)$ for $n \geq 1$.
 Then 
\begin{equation}\label{eq:udes1}
U^{(\PP)}_{u,n}(x,{\bf z}_\infty) = \Theta_u(h_n) = \sum_{O \in \mathcal{O}^{(\PP)}_{u,n},I_u(O) = O}sgn(O)wt(O)
\end{equation}
where $\mathcal{O}^{(\PP)}_{u,n}$ is the set of objects and 
$I_u$ is the involution defined above. Moreover  $O =(B,w)$, where 
$B=(b_1, \ldots, b_k)$ and $w =w_1 \ldots w_n$,  is a 
fixed point of $I_u$ if and only if it has the following two properties: 
\begin{enumerate}
\item there are no cells labeled with $x$ in $O$, i.e., the 
elements of $w$ in each brick of $O$ are weakly increasing  and 
\item if \begin{math}b_i\end{math} and \begin{math}b_{i+1}\end{math} are two consecutive bricks in \begin{math}O\end{math}, then 
either (a) there is a weak increase between \begin{math}b_i\end{math} and \begin{math}b_{i+1}\end{math}, i.e., $w_{\sum_{j=1}^i |b_j|} \leq w_{1+\sum_{j=1}^i |b_j|}$, or (b) 
there is a decrease  between \begin{math}b_i\end{math} and \begin{math}b_{i+1}\end{math}, i.e., $w_{\sum_{j=1}^i |b_j|} > w_{1+\sum_{j=1}^i |b_j|}$,
 but there is a $u$-match contained in the elements of the cells of \begin{math}b_i\end{math} and \begin{math}b_{i+1}\end{math} which must necessarily involve 
$w_{\sum_{j=1}^i |b_j|}$ and  $w_{1+\sum_{j=1}^i |b_j|}$. 
\end{enumerate}
\end{lemma}

Clearly, if we restrict to the alphabet $[k]$ instead of $\PP$, we will get the same two lemmas except that the words all have to be in $[k]^*$ rather 
than in $\PP^*$.

Next we want to consider what happens when we replace $u$-matches by exact $u$-matches. 
We can follow the same steps to interpret  
$EU^{(\PP)}_u(x,{\bf z}_\infty,t)$.
That is,  
\begin{equation}
EU^{(\PP)}_u(x,{\bf z}_\infty,t) 
= \frac{1}{1+\sum_{n \geq 1} EN^{(\PP)}_{u,n}(x,{\bf z}_\infty)t^n}.
\end{equation}
Thus if we let 
$\Gamma_u(e_n) = (-1)^n EN^{(\PP)}_{u,n}(x,{\bf z}_\infty)$ for 
$n \geq 1$ and $\Gamma_u(e_0) =1$, we see that 
\begin{eqnarray*}
\Gamma_u(H(t)) &=& 1 + \sum_{n \geq 1} \Gamma_u(h_n) \\
&=& \Gamma_u\left( \frac{1}{E(-t))}\right) = \frac{1}{1+ \sum_{n \geq 1} (-1)^n \Gamma_u(e_n)} \\
&=& \frac{1}{1 + \sum_{n \geq 1} EN^{(\PP)}_{u,n}(x,{\bf z}_\infty)t^n} = EU^{(\PP)}_u(x,{\bf z}_\infty,t).
\end{eqnarray*}
Thus it follows that 
$\Gamma_u(h_n) = EU^{(\PP)}_{u,n}(x,{\bf z}_\infty)$.

By (\ref{htoe}), we have  that 
\begin{eqnarray}\label{eq:basic1e}
\Gamma_{u}(h_n) &=&  \sum_{\la \vdash n} (-1)^{n-\ell(\la)} 
B_{\la,n}~ \Gamma_{u}(e_\la) \nonumber \\
&=&  \sum_{\la \vdash n} (-1)^{n-\ell(\la)} \sum_{(b_1, \ldots, 
b_{\ell(\la)}) \in \mathcal{B}_{\la,n}} \prod_{i=1}^{\ell(\la)}  
(-1)^{b_i}  EN^{(\PP)}_{u,b_i}(x,{\bf z}_\infty) \nonumber \\
&=& \sum_{\la \vdash n} (-1)^{\ell(\la)} 
\sum_{(b_1, \ldots, b_{\ell(\la)}) \in \mathcal{B}_{\la,n}} 
\prod_{i=1}^{\ell(\la)} EN^{(\PP)}_{u,b_i}(x,{\bf z}_\infty)
\end{eqnarray}

Again we can give a combinatorial interpretation to the 
right-hand side of (\ref{eq:basic1e}). Fix a partition 
$\la$ of $n$ and a $\lambda$-brick tabloid 
$B=(b_1, \ldots, b_{\ell(\la)})$. We will 
interpret $\prod_{i=1}^{\ell(\la)}  EN^{(\PP)}_{u,b_i}(x,{\bf z}_\infty) $ as 
the number of ways of picking words 
$(w^{(1)}, \ldots, w^{(\ell(\la))})$ such that for each 
$i$,  $w^{(i)} \in \PP^{b_i}$ is a word such 
that $\eumch(w) = 0$ and assigning 
a weight to this $\ell(\la)$-tuple to be  
$\prod_{i=1}^{\ell(\la)} x^{\des(w^{(i)})+1} \overline{z}^{w^{(i)}}$.

Following the same steps that we did to interpret $\Theta_u(h_n)$, 
we let $\mathcal{EO}^{(\PP)}_{u,n}$ denote the set of 
all filled-labeled-brick tabloids constructed in this way. 
That is, $\mathcal{EO}^{(\PP)}_{u,n}$ consists of all pairs 
$O=(B,w)$ where 
\begin{enumerate}
\item $B=(b_1, \ldots, b_{\ell(\la)})$ is brick tabloid 
of shape $(n)$, 
\item $w =w_1 \ldots w_n \in \PP^n$ such that there is 
no exact $u$-match of $\sg$ which is entirely contained in a single 
brick of $B$, and 
\item if there is a cell $c$ such that a brick $b_i$ contains both cells $c$ and $c+1$
and $w_c > w_{c+1}$, then cell $c$ is labeled with a $x$ and the last cell of any brick is labeled with $-x$.
\end{enumerate}
The sign of $O$, $\sgn(O)$, is $(-1)^{\ell(\la)}$ and the weight 
of $O$, $wt(O)$,  is $x^{\ell(\la) +\mathrm{intdes}(\sg)}\overline{z}^w$.  
Then as before we can conclude  
\begin{equation}\label{eq:ebasic2}
\Gamma_{u}(h_n) = \sum_{O \in \mathcal{EO}^{(\PP)}_{u,n}} sgn(O) wt(O).
\end{equation}

At this point, we can define an involution $J_u$ exactly as we did for 
$I_u$ except replace $u$-match by exact $u$-matches in the definitions. 
This will allow us to prove the following two lemmas.

\begin{lemma}\label{lem:euinc} 
Suppose that $u$ is a word in $\PP^+$ such that $\des(u) =0$. 
Let $\Gamma_u:\Lambda \rightarrow 
\mathbb{Q}(x)$ be the  ring homomorphism defined by setting  
$\Gamma_{u}(e_0) =1$ and  
$\Gamma_{u}(e_n) = (-1)^n EN^{(\PP)}_{u,n}(x,{\bf z}_\infty)$ for $n \geq 1$.
Then 
\begin{equation}\label{eq:euinc}
EU^{(\PP)}_{u,n}(x,{\bf z}_\infty) = \theta_u(h_n) = \sum_{((b_1, \ldots, b_k),w) \in \mathcal{WIEO}_{u,n}} 
(-x)^k \overline{z}^w
\end{equation}
where $\mathcal{WIEO}_{u,n}$ is the set all $(B,w) \in \mathcal{EO}_{u,n}$ 
such that $J_u(B,w) = (B,w)$ and $w$ is weakly increasing.
\end{lemma}

\begin{lemma} \label{lem:eudes1}
Suppose that $u \in \PP^+$ and $\des(u)=1$.  Let $\Gamma_u:\Lambda \rightarrow 
\mathbb{Q}(y)$ be the  ring homomorphism defined by setting  
$\Gamma_{u}(e_0) =1$ and  
$\Gamma_{u}(e_n) = (-1)^nEN^{(\PP)}_{u,n}(x,{\bf z}_\infty)$ for $n \geq 1$.
 Then 
\begin{equation} \label{eq:eudes1}
EU^{(\PP)}_{u,n}(x,{\bf z}_\infty) = \Gamma_u(h_n) = \sum_{O \in \mathcal{EO}^{(\PP)}_{u,n},J_u(O) = O}sgn(O)wt(O)
\end{equation}
where $\mathcal{EO}^{(\PP)}_{u,n}$ is the set of objects and 
$J_u$ is the involution defined above. Moreover $O =(B,w)$, where 
$B=(b_1, \ldots, b_k)$ and $w =w_1 \ldots w_n$,  is a 
fixed point of $J_u$ if and only if it has the following two properties: 
\begin{enumerate}
\item there are no cells labeled with $x$ in $O$, i.e., the 
elements of $w$ in each brick of $O$ are weakly increasing  and 
\item if \begin{math}b_i\end{math} and \begin{math}b_{i+1}\end{math} are two consecutive bricks in \begin{math}O\end{math}, then 
either (a) there is a weak increase between \begin{math}b_i\end{math} and \begin{math}b_{i+1}\end{math}, i.e., $w_{\sum_{j=1}^i |b_j|} \leq w_{1+\sum_{j=1}^i |b_j|}$, or (b) 
there is a decrease  between \begin{math}b_i\end{math} and \begin{math}b_{i+1}\end{math}, i.e., $w_{\sum_{j=1}^i |b_j|} > w_{1+\sum_{j=1}^i |b_j|}$,
 but there is an exact $u$-match contained in the elements of the cells of \begin{math}b_i\end{math} and \begin{math}b_{i+1}\end{math} which must necessarily involve 
$w_{\sum_{j=1}^i |b_j|}$ and  $w_{1+\sum_{j=1}^i |b_j|}$. 
\end{enumerate}
\end{lemma}

\section{The case where $u = u_1 \ldots u_j$, $\des(u)=1$, and $u_1 > u_j$}

In this section, we shall consider the problem 
of computing the generating functions \\
$\mathcal{N}^{(\PP)}_u(x,{\bf z}_\infty,t)$, 
$\mathcal{N}^{(k)}_u(x,{\bf z}_k,t)$, $\mathcal{EN}^{(\PP)}_u (x,{\bf z}_\infty,t)$, and 
$\mathcal{EN}^{(k)}_u(x,{\bf z}_k,t)$ for  
$u = u_1 \ldots u_j$ such that $\des(u)=1$, and $u_1 > u_j$. 

Now suppose that $u = u_1 \ldots u_j$, $\red(u) =u$, $u_1 > u_j$, and 
$\des(u)=1$.  Let $1 \leq s <j$ be the position 
such that $u_s>u_{s+1}$ so that $u_1 \leq \cdots \leq u_s>u_{s+1} \leq \cdots 
\leq u_j$.  Then  
$St^{(\PP)}(u) \subseteq \{s+1, \ldots, j\}$ since  if we try to start 
a match at one of the positions $2, \ldots, s$, the 
descent in the second match 
would not be in the right place. It follows that $u$ automatically 
has the $\PP$-weakly decreasing overlapping property and 
the $[k]$-weakly decreasing overlapping property for any $k \geq 2$.  
 
We start by considering a special class of words $u = u_1 \ldots u_j$
which have the $\PP$-minimal overlapping property ($[k]$-minimal 
overlapping property). 
This means that any two consecutive $u$-matches can share at most one letter. 
For example $u=2341$ has the $\PP$-minimal overlapping property while 
$u=3412$ does not have the $\PP$-minimal overlapping property since in 
the word $w= 563412$, the $u$-matches $5634$ and $3412$ share two letters. 

\begin{figure}[htbp]
  \begin{center}
    \includegraphics[width=0.6\textwidth]{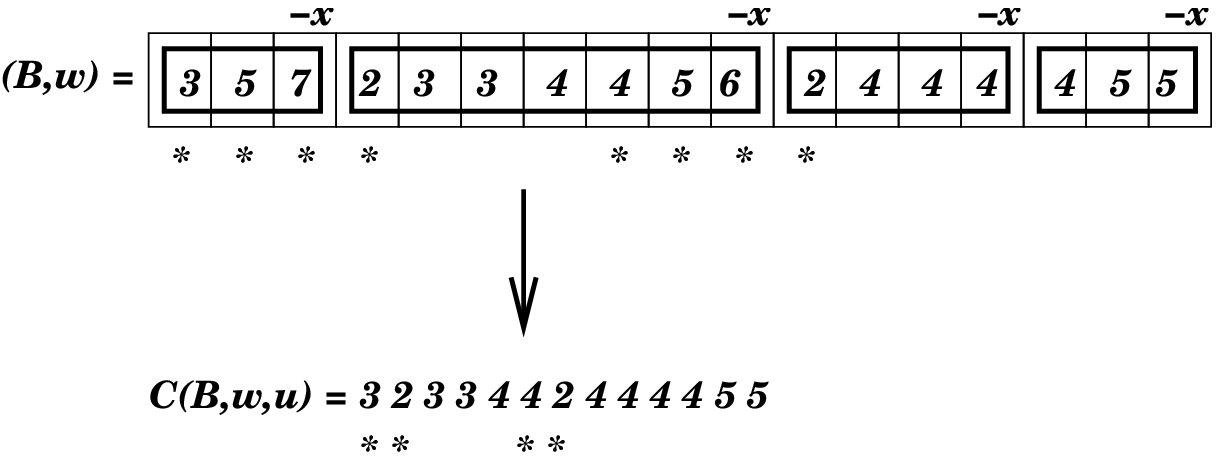}
    \caption{A fixed point of $I_{2341}$.}
    \label{fig:Collapse2341}
  \end{center}
\end{figure}

Thus assume that  $u = u_1 \ldots u_j$, $\red(u) =u$, $\des(u)=1$,  $u_1 > u_j$, and $u$ has the 
$\PP$-minimal overlapping property. First we introduce what we shall call the 
{\em collapse map} which 
maps fixed points of $I_u$ or $J_u$ to a certain subset of words 
in $\PP^*$.  This is best explained through an example. Suppose 
that $u = 2341$ and we want to compute $U^{(7)}_{2341}(x,{\bf z}_7,t)$.
By (\ref{eq:udes1}), we know 
that 
\begin{equation}\label{eq:kudes1}
U^{(7)}_{u,n}(x,{\bf z}_7) = \sum_{O \in \mathcal{O}^{(7)}_{u,n},I_u(O) = O}sgn(O)wt(O).
\end{equation}
Now suppose that we are given a fixed point $(B,w)$ of $I_u$ where 
$B=(b_1, \ldots, b_k)$ and $w=w_1 \ldots w_n$ such as the one pictured 
in Figure \ref{fig:Collapse2341}.  We know that to be a fixed point of $I_u$, 
$w$ must be  weakly increasing within bricks of $B$ and that for any $i <k$, 
if $c$ is last cell in brick $b_i$ and $w_c > w_{c+1}$, then there must 
be a $u$-match in $w$ which is contained in the cells of $b_i$ and 
$b_{i+1}$.  In our particular example, since $u=2341$ has a single descent, 
this match must involve the last three cells of $b_i$ and the first cell 
of $b_{i+1}$.  In Figure \ref{fig:Collapse2341}, we have indicated 
the two such $2341$-matches in our example by placing stars below the cells in the 
$2341$-matches. In this case the collapse map just maps $(B,w)$ to the 
word $v =C(B,w,u)$ which is the result of starting with 
$w$ and removing the letters in all such matches that do not 
correspond to the end points of the match.  This process is pictured 
in Figure \ref{fig:Collapse2341} where again we have starred the elements 
in $C(B,w,u)$ that remain from the original $2341$-matches in $w$. What 
makes the case where $u$ has the minimal overlapping property easier 
is that, since any two consecutive $u$-matches can share at most letter, 
there is no possibility that an end point of $u$-match in $w$ occurs in 
the middle of another $u$-match in $w$ so that the letters that 
we remove from $w$ for any pair of $u$-matches are disjoint from each other.

The next question that we want to consider is how can 
we construct all the fixed points of $(B,w)$ of $I_u$ such that 
$C(B,w,v)$ is equal to a given word $v=v_1 \ldots v_n$. 
First it easy to see that the only descents that 
appear in a word $C(B,w,u)$ must come from  $2341$-matches that straddled two 
bricks in $B$.  Thus if $v_s > v_{s+1}$, then $v_s$ must have played the role 
of 2 in the original $2341$-match and $v_{s+1}$ must have played the role 
of 1 in the original $2341$-match. Such a requirement rules out certain 
words from being in the range of the collapse map $C$. For example, 
suppose that the underlying alphabet is $[7]$. Then    
if  $v_s =6$ and $v_{s+1} =i$ where $i < 6$, then $v$ 
could not have come from 
the collapse of $2341$-match because we can not add two letters which 
could play the role of 3 and 4  in the 2341-match.    If 
we consider the first descent $32$ in the $C(B,w,u)$ of 
Figure \ref{fig:Collapse2341}, then we see there are many ways that we 
could add the two middle letters. That is, the 
original 2341-match could have been any $3cd2$ where 
$c < d$ and $c,d \in \{4,5,6,7\}$. It follows 
that the extra weight from these possibilities that 
is not included in $\overline{z}^{C(B,w,u)}t^{|C(B,w,u)|}$ in this case 
would be $-x t^2 \sum_{4 \leq c < d \leq 7} z_c z_d$. Here the $-x$ comes from 
the fact that we know that the original match straddled two bricks 
and there is a weight of $-x$ associated with the end point of 
the first of those two bricks. On the other hand, if $v_s \leq v_{s+1}$, 
then we have only two choices. That is, either cell $s$ was the end of a brick 
or cell $s$ was an internal cell of a brick.  This implies that 
each weak rise in $v$ contributes a factor of $(1-x)$ since if 
$s$ is at the end of a brick, there is a weight of $-x$ associated with 
the last cell of a brick. In this way, we can associate a 
weight with each weak rise or descent of $v$ which 
will allow us to compute
$$\sum_{\overset{(B,w) \ \mathrm{is\ a\ fixed\ point\ of} \ I_u}{C(B,w,u) =v}} sgn(B,w) wt(B,w).$$

In our case where $u =2341$ and $k=7$, 
the weights associated with the descents are 
given in the table below.

\begin{center}
\begin{tabular}{|l|l|}
\hline
Descents & $wt_{2341,7}(ji)$\\
\hline
$7i\  (i < 7)$ & 0 \\
\hline
$6i\  (i < 6)$ & $0$\\
\hline 
$5i \ (i < 5)$ & $-x z_6z_7t^2$ \\
\hline
$4i \ (i < 4)$ & $-x (z_5z_6+z_5z_7+z_6z_7)t^2$ \\
\hline
$3i \ (i < 3)$ & $-x (z_4z_5+z_4z_6+z_4z_7+z_5z_6+z_5z_7+z_6z_7)t^2$ \\
\hline
$21$ & $-x (z_3z_4+ z_3z_5+z_3z_6+z_3z_7+z_4z_5+z_4z_6+z_4z_7+z_5z_6+z_5z_7 + z_6z_7)t^2$ \\
\hline
\end{tabular}\\
The weights $wt_{2341,7}(ji)$.
\end{center}

However, if $u=2341$ and we want to compute 
$U^{(\PP)}_{u,n}(x,{\bf z}_{\infty})$, the 
weights for any descent $ji$ would be $-xt^2\sum_{j < c <d} z_cz_d$ which 
is an infinite sum.

Going back to our example where $u =2341$ and $k =7$, 
it follows 
that for any $v \in [7]^+$, 
\begin{equation}\label{eq:v72341}
\sum_{\overset{(B,w) \ \mathrm{is\ a\ fixed\ point\ of} \ I_u}{C(B,w,u) =v}}
sgn(B,w) wt(B,w) =
-x \overline{z}^v t^{|v|}(1-x)^{\wris(v)} \prod_{s \in Des(v)} wt_{2341,7}(v_s v_{s+1}).
\end{equation}
Here the initial $-x$ comes from the fact that the last cell of $(B,w)$ always 
contributes a $-x$ since the last cell is at the end of a brick. 
It follows that 
\begin{eqnarray}\label{eq:U72341}
U^{(7)}_{2341}(x,{\bf z}_7,t) &=& 1 + \sum_{n \geq 1} U^{(7)}_{2341,n}(x,{\bf z}_7) t^n \nonumber \\
&=& 1 + \sum_{v \in [7]^+} -x(1-x)^{\wris(v)} \overline{z}^v t^{|v|}\prod_{s \in Des(v)} wt_{2341,7}(v_s v_{s+1}).
\end{eqnarray}
Hence we could compute 
$\displaystyle \mathcal{N}^{(7)}_{2341,n}(x,{\bf z}_7,t) = \frac{1}{U^{(7)}_{2341,n}(x,{\bf z}_7,t)}$
if we can compute the right-hand side of (\ref{eq:U72341}).

The case of exact matches is even simpler. In that case, we want to compute 
$$\sum_{\overset{(B,w) \ \mathrm{is \ a \ fixed \ point \ of} \ J_u}{C(B,w,u) =v}} sgn(B,w) wt(B,w).$$
Going back to our example of $u=2341$ over the alphabet 
$[7]$, we see the only descents that 
appear in a word $v= C(B,w,u)$ must come from  exact $2341$-matches that straddled two 
bricks in $B$.  Thus if $v_s > v_{s+1}$, then it 
must be the case that $v_s=2$, $v_{s+1}=1$ and we 
must have eliminated a 3 and 4 from $w$.  Thus we want to compute 
$EU^{(\PP)}_{2341,n}(x,{\bf z}_{\infty})$ or 
$EU^{(k)}_{2341,n}(x,{\bf z}_k)$ for $k \geq 4$, the 
weights would be the following. 
\begin{center}
\begin{tabular}{|l|l|}
\hline
Descents & weight $ewt_{2341,\PP}(ji)$\\
\hline
$ji\  \mbox{where either} \ j \neq 2 \ \mbox{or} \ i \neq 1$  & 0 \\
\hline 
21 & $-x z_3z_4t^2$ \\
\hline
\end{tabular}\\
The weights $ewt_{2341}(ji)$.
\end{center}
It follows 
that for any $v \in \PP^+$, 
\begin{equation}\label{eq:ev72341}
\sum_{\overset{(B,w) \ \mathrm{is \ a \ fixed \ point \ of} \ J_{2341}}{C(B,w,2341) =v}} sgn(B,w) wt(B,w) =
-x\overline{z}^vt^{|v|}(1-x)^{\wris(v)} \prod_{s \in Des(v)} ewt_{2341,\PP}(v_s v_{s+1}).
\end{equation}
and
\begin{eqnarray}\label{eq:EU72341}
EU^{(\PP)}_{2341,n}(x,{\bf z}_\infty,t) &=& 1 + \sum_{n \geq 1} EU^{(\PP)}_{2341,n}(x,{\bf z}_\infty) t^n \nonumber \\
&=& 1 + \sum_{v \in \PP^+} -x\overline{z}^vt^{|v|}(1-x)^{\wris(v)} \prod_{s \in Des(v)} ewt_{2341,\PP}(v_s v_{s+1}).
\end{eqnarray}

In our case, $\prod_{s \in Des(v)} ewt_{2341,\PP}(v_s v_{s+1}) =0$ 
unless the only descents in $v$ are of the from $21$. It follows 
that the only nonempty words $v$ that can contribute to 
(\ref{eq:EU72341}) are words $v$ of the form 
$w$ or of the form 
$1^{a_1}2^{b_1}21 1^{a_2}2^{b_2}21 \ldots 1^{a_r}2^{b_r}21w$ for some 
$r \geq 1$ where $w$ is a weakly increasing word.
Let 
$$
W(x,\mathbf{z}_\infty,t) := \prod_{i=1}^\infty \frac{1}{(1-(1-x)z_it)}.$$
The generating function of $-x\overline{z}^v t^{|v|}(1-x)^{\wris(v)} \prod_{s \in Des(v)} ewt_{2341,7}(v_s v_{s+1})$ over all nonempty weakly increasing words 
is just 
\begin{equation}\label{2341wi}
\frac{-x}{(1-x)}\left( -1 +W(x,\mathbf{z}_\infty,t)\right).
\end{equation} 
The generating function of $-x\overline{z}^v t^{|v|}(1-x)^{\wris(v)} \prod_{s \in Des(v)} ewt_{2341,7}(v_s v_{s+1})$ over all words $v$ of the form 
$1^{a_1}2^{b_1}2$ is 
$$\frac{z_2t}{(1-(1-x)z_1t)(1-(1-x)z_2t)}.$$
The generating function of $-x\overline{z}^v t^{|v|}(1-x)^{\wris(v)} \prod_{s \in Des(v)} ewt_{2341,7}(v_s v_{s+1})$ over all words $v$ of the form 
$1^{a_1}2^{b_1}211^{a_2}2^{b_2}21 \ldots 1^{a_r}2^{b_r}21w$ where 
$w$ is weakly increasing is 
$$-xW(x,\mathbf{z}_\infty,t)\left(\frac{z_2t}{(1-(1-x)z_1t)(1-(1-x)z_2t)}\right)^r 
(-xz_1z_3z_4t^3)^r(1-x)^{r-1}.$$
Here the term $(-xz_1z_3z_4t^3)^r$ comes from the weights 
$ewt_{2341,\PP}(21)$ that arise from the descents 21 and 
$(1-x)^{r-1}$ comes from the weights of the rises coming from the first 
$r-1$ 1s which are the second elements of the descents 21. 
It follows that the generating function of $-x\overline{z}^v t^{|v|}(1-x)^{\wris(v)} \prod_{s \in Des(v)} ewt_{2341,7}(v_s v_{s+1})$ over all $v$ such 
that $v$ is of the form 
$1^{a_1}2^{b_1}21 1^{a_2}2^{b_2}21 \ldots 1^{a_r}2^{b_r}21w$ for some 
$r \geq 1$ where $w$ is a weakly increasing word is equal to 
\begin{align}\label{2341nwi}
& -xW(x,\mathbf{z}_\infty,t)\sum_{r \geq 1}
\left(\frac{-xz_1z_2z_3z_4t^4}{(1-(1-x)z_1t)(1-(1-x)z_2t)}\right)^r 
(1-x)^{r-1} = \nonumber \\
&\frac{-xW(x,\mathbf{z}_\infty,t)}{(1-x)}\left( -1 + 
\frac{1}{1-\frac{-xz_1z_2z_3z_4t^4(1-x)}{(1-(1-x)z_1t)(1-(1-x)z_2t)}}\right)
\end{align}

Putting (\ref{2341wi}) and (\ref{2341nwi}) together we see that 
\begin{eqnarray}
EU^{(\PP)}_{2341,n}(x,{\bf z}_\infty,t) &=& 
1+ \frac{-x}{(1-x)}\left( -1 +W(x,\mathbf{z}_\infty,t)\right) + 
\nonumber \\
&&\frac{-xW(x,\mathbf{z}_\infty,t)}{(1-x)}\left( -1 + 
\frac{1}{1-\frac{-xz_1z_2z_3z_4t^4(1-x)}{(1-(1-x)z_1t)(1-(1-x)z_2t)}}\right)
\nonumber \\
&=& 1+\frac{x}{(1-x)}-\frac{xW(x,\mathbf{z}_\infty,t)}{(1-x)}
\frac{1}{1+\frac{xz_1z_2z_3z_4t^4(1-x)}{(1-(1-x)z_1t)(1-(1-x)z_2t)}}.
\end{eqnarray}

Thus 
\begin{equation}\label{finale2341}
\mathcal{EN}_{2341}^{(\PP)}(x,\mathbf{z}_\infty,t) = 
\frac{1}{1+\frac{x}{(1-x)}-\frac{xW(x,\mathbf{z}_\infty,t)}{(1-x)}
\frac{1}{1+\frac{xz_1z_2z_3z_4t^4(1-x)}{(1-(1-x)z_1t)(1-(1-x)z_2t)}}}.
\end{equation}

It should be clear from our arguments that the only role that 
the 2 and 3 played in the final form 
$\mathcal{EN}_{2341}^{(\PP)}(x,\mathbf{z}_\infty,t)$ was to contribute 
a factor of $z_3z_4t^2$ to the expression 
$\frac{xz_1z_2z_3z_4t^4(1-x)}{(1-(1-x)z_1t)(1-(1-x)z_2t)}$ on 
the right hand side of (\ref{finale2341}). Thus our arguments show 
that if $u =2 \alpha 1$ where $\alpha$ is non-empty weakly increasing word in 
$\{2,3, \ldots \}^*$, then we have the following theorem. 

\begin{theorem} Let $u =2 \alpha 1$ where 
$\alpha$ is non-empty weakly increasing word in 
$\{2,3, \ldots \}^*$. Then 
\begin{equation}\label{finale2alpha1}
\mathcal{EN}_{2\alpha 1}^{(\PP)}(x,\mathbf{z}_\infty,t) = 
\frac{1}{1+\frac{x}{(1-x)}-\frac{xW(x,\mathbf{z}_\infty,t)}{(1-x)}
\frac{1}{1+\frac{xz_1z_2\overline{z}^\alpha t^{2+|\alpha|}(1-x)}{(1-(1-x)z_1t)(1-(1-x)z_2t)}}}.
\end{equation}
\end{theorem}

\begin{figure}[htbp]
  \begin{center}
    \includegraphics[width=0.6\textwidth]{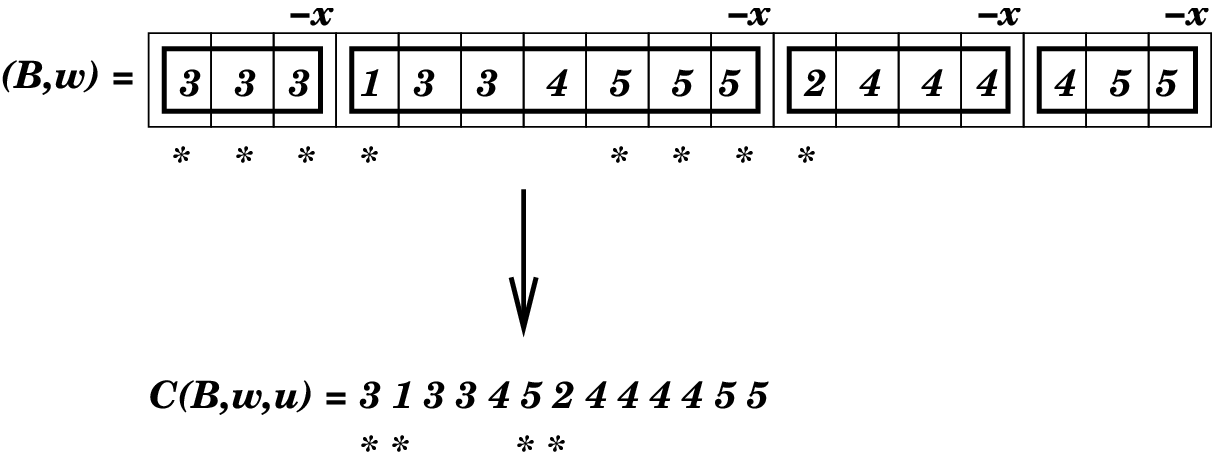}
    \caption{A fixed point of $I_{2221}$.}
    \label{fig:Collapse2221}
  \end{center}
\end{figure}

Other examples where the weights $wt_{u,\PP}(ij)$ 
are easy to compute are words 
of the form $u=2^r1$ or $u=21^r$ where $r \geq 2$.  
It is easy to see that both $2^r1$ and $21^r$ have the minimal overlapping property. 
In this case, the only $u$-matches are of the form 
$b^ra$ where $b > a \geq 1$ if $u = 2^r1$ or $ba^r$ where $b > a \geq 1$ if $u = 21^r$.  
For example, suppose that $u =2^3 1$. Then 
in Figure \ref{fig:Collapse2221}, we have pictured a fixed point of $I_u$ 
where we have indicated 
the two $2^31$-matches in our example by placing stars below the cells in the 
$2^31$-matches. In the case the collapse map just maps $(B,w)$ to the 
word $v =C(B,w,u)$ which is the result of starting with 
$w$ and removing the letters in all such matches that do not 
correspond to the end points of the match.  This process is pictured 
in Figure \ref{fig:Collapse2221} where again we have starred the elements 
in $C(B,w,u)$ that remain from the original $2^31$-matches in $w$.

In this case, if we have a descent $ji$, 
then $wt_{2^31,k}(ji) = wt_{2221,\PP}(ji) = -xz_j^2t^2$ since 
we will always add back two $j$s for each descent of 
of the form $ji$. Thus if $u =2^31$  
it follows 
that for any $v \in \PP^+$, 
\begin{equation}\label{eq:v72221}
\sum_{\overset{(B,w) \ \mbox{is a fixed point of} \ I_{2221}}{C(B,w,2221) =v}} sgn(B,w) wt(B,w) =
-x(1-x)^{\wris(v)} \prod_{s \in Des(v)} -x z_{v_s}^2t^2.
\end{equation}
and 
\begin{eqnarray}\label{eq:U72221}
U^{(\PP)}_{2221,n}(x,{\bf z}_\infty,t) &=& 1 + \sum_{n \geq 1} U^{(\PP)}_{2221,n}(x,{\bf z}_\infty) t^n \nonumber \\
&=& 1 + \sum_{v \in \PP^+} -x(1-x)^{\wris(v)} \prod_{s \in Des(v)} 
-x z_{v_s}^2t^2.
\end{eqnarray}
Hence we could compute 
$\displaystyle \mathcal{N}^{(\PP)}_{2221,n}(x,{\bf z}_\infty,t) = \frac{1}{U^{(\PP)}_{2221,n}(x,{\bf z}_\infty,t)}$
if we can compute the right-hand side of (\ref{eq:U72221}).

When $u$ does not have the minimal overlapping property, 
we can obtain similar results but the collapse maps 
and the weights $wt_u(ji)$ are more complicated. Again 
this is best explained through an example.
Suppose that $u =3412$ and $k=8$.

\begin{figure}[htbp]
  \begin{center}
    \includegraphics[width=0.6\textwidth]{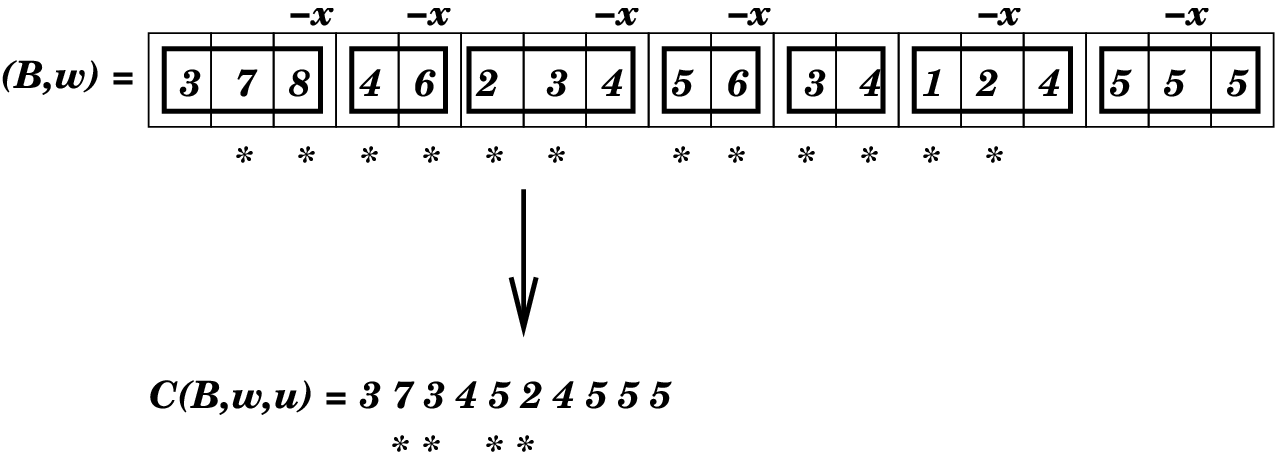}
    \caption{A fixed point of $I_{3412}$.}
    \label{fig:Collapse3412}
  \end{center}
\end{figure}

When $u$ does not have the minimal overlapping property, then we 
can have a situation such as the one pictured in 
Figure \ref{fig:Collapse3412}.  If we look 
at the descents between bricks 1 and 2 which correspond 
to the $u$-match 7846, we see that we would like to 
eliminate the 8 and 4.  However, this $u$-match overlaps 
the $u$-match associated with the descent between bricks 2 and 3 
which is 4623.  Thus we would also like to eliminate the 6 and 2. 
We will say that two such matches are {\em linked} if one of 
the end points of first match is one of middle elements of 
the second match. Depending on the pattern we could 
have a series of $u$-matches in a fixed point of $(B,w)$ which 
are linked. The fact that we are assuming that 
$u= u_1 \ldots u_j$ where $u_1 > u_j$ and $\des(u) =1$ implies 
that $u$ has the $\PP$-weakly decreasing overlapping property or $[k]$-weakly decreasing overlapping property. 
It is then easy to see that if $w=w_1 \ldots w_n$ is a word such that there is a 
$u$-match starting position 1 and a $u$-match ending at position $n$ 
and any two consecutive $u$-matches in $w$ are linked, then 
$w_1 > w_n$.  In such a situation,  the collapse map will 
eliminate all the symbols except for the first element 
of the first match and last element of the last match 
in a maximal sequence of linked $u$-matches which will result 
in a descent. This is illustrated 
in Figure \ref{fig:Collapse3412} where we have two maximal blocks of linked $3412$-matches. 
Thus in the linked $3412$-matches in cells 2 through 
7, we keep only the 7 and the 3 and in the linked matches 
in cells 9 through 14, we keep only the 5 and the 2.  
Because we are assuming that $u_1 > u_j$, we know 
that maximal blocks of linked $u$-matches must be finite 
since the end point of such matches must strictly decrease. 
When we see a descent $ji$ in a word $C(B,w,u)$, the 
weight associated with such a decent is now more complicated. 
For example, in our case where $u =3412$ and $k=8$, 
a decent of the form $72$ can correspond to 
a single $3412$-match which would have to be of the 
form $7812$, it could correspond to a maximum block 
with 2 linked $3412$-matches in which case it must 
be of the form $78cd12$ where $3 \leq c < d \leq 6$, or 
it could correspond to a maximum block with 3 linked 
$3412$-matches in which case it must 
be $78563412$.  Thus 
$$wt_{3412}(72) = -xz_1z_8t^2 + x^2 t^4 z_1z_8 
\left( \sum_{3 \leq c < d \leq 6} 
z_c z_d \right) -x^3t^6 z_1z_3z_4z_5 z_6 z_8$$
On the other hand a descent of 
the form $ji$ where $j-i \leq 2$ can only correspond to a 
single $3412$-match so that $wt_{3412}(ji) = 
-xt^2(\sum_{j < s \leq 8}z_s)(\sum_{1 \leq t < i} z_t)$. 

We give the weights associated with the descents 
for $u=3412$ and $k=8$ in the following table. 

\begin{center}
\begin{tabular}{|l|l|}
\hline
Descents & $wt_{3412,8}(ji)$\\
\hline
$8i\  (i < 8)$ & 0 \\
\hline
$j1\  (j > 2)$ & 0\\
\hline 
$ji \ (j > i \ \& \ j-i \leq 2$ & 
$-xt^2(\sum_{j < s \leq 8}z_s)(\sum_{1 \leq t < i} z_t)$ \\
\hline
72 & $-xz_1z_8t^2 + x^2 t^4 z_1z_8 (\sum_{3 \leq c < d \leq 6} 
z_c z_d) -x^3t^6 z_1z_3z_4z_5 z_5 z_6 z_8$\\
\hline
62 & 
$-xt^2(z_7+z_8)z_1 +x^2t^4(z_7+z_8)z_1(\sum_{3 \leq c < d \leq 5} z_cz_d)$ \\
\hline
52 & 
$-xt^2(z_6+z_7+z_8)z_1 +x^2t^4(z_6+z_7+z_8)z_1z_3z_4$ \\
\hline
73& $-xt^2z_8(z_1+z_2) +x^2t^4z_8(z_1+z_2) +x^2t^4z_8(z_1+z_2) 
(\sum_{4 \leq c < d \leq 6} z_cz_d)$ \\
\hline
63 & $-xt^2(z_7+z_8)(z_1+z_2)+x^2t^4(z_7+z_8)(z_1+z_2)z_4z_5$\\
\hline
74 & $-x t^2z_8(z_1+z_2+z_3) + x^2t^4z_8(z_1+z_2+z_3)z_5z_6 $ \\
\hline
\end{tabular}\\
The weights $wt_{3412,8}(ji)$.
\end{center}

It follows 
that for any $v \in [8]^+$, 
\begin{equation}\label{eq:v83412}
\sum_{\overset{(B,w) \ \mathrm{is \ a \ fixed \ point \ of} \ I_{3412}}{C(B,w,3412) =v}} sgn(B,w) wt_{3412}(B,w) =
-x\overline{z}^v t^{|v|}(1-x)^{\wris(v)} \prod_{s \in Des(v)} wt_{3412,8}(v_s v_{s+1}).
\end{equation}
and
\begin{eqnarray}\label{eq:U83412}
U^{(8)}_{3412,n}(x,{\bf z}_8,t) &=& 1 + \sum_{n \geq 1} U^{(8)}_{3412,n}(x,{\bf z}_8) t^n \nonumber \\
&=& 1 + \sum_{v \in [8]^+} -x\overline{z}^v t^{|v|}(1-x)^{\wris(v)} \prod_{s \in Des(v)} 
wt_{3412,8}(v_s v_{s+1}).
\end{eqnarray}

What we need to be able to compute the right-hand sides of either 
(\ref{eq:U72341}), (\ref{eq:EU72341}), (\ref{eq:U72221}), or (\ref{eq:U83412}) 
is the generating functions over all words $v \in \PP^*$ where we not only keep track of the descents of $P$ but also 
of type of descents of $P$. 

By Theorem \ref{thm:labdes}, 
we know that

\begin{equation}\label{PPDES1}
\frac{1}{1 - \sum_{n \geq 1} t^n \sum_{S \subseteq \PP, |S| =n} 
DXZ(S)} = 1+\sum_{w = w_1 \ldots w_n \in \PP^+} t^{|w|}\overline{z}^w \prod_{i \in Des(w)}
x_{w_i w_{i+1}}.
\end{equation}
Hence 
\begin{align}\label{PPDES2}
\sum_{w = w_1 \ldots w_n \in \PP^+} t^{|w|}\overline{z}^w \prod_{i \in Des(w)}
x_{w_i w_{i+1}} &= 
\left(\frac{1}{1 - \sum_{n \geq 1} t^n \sum_{S \subseteq \PP, |S| =n} 
DXZ(S)}\right) -1 \nonumber \\
&=\frac{\sum_{n \geq 1} t^n \sum_{S \subseteq \PP, |S| =n} 
DXZ(S)}{1 - \sum_{n \geq 1} t^n \sum_{S \subseteq \PP, |S| =n} 
DXZ(S)}.
\end{align}
Next suppose that  we replace $t$ by $yt$ and $x_{ij}$ by 
$\frac{x_{ij}}{y}$.  Under this substitution 
the left-hand side in (\ref{PPDES2}) becomes 
 
$$\sum_{w = w_1 \ldots w_n \in \PP^+} t^{|w|}y^{\wris(w)+1}\overline{z}^w \prod_{i \in Des(w)}x_{w_i w_{i+1}}.$$
Note that for $S=\{j_1 < \cdots < j_k\}$ where $k \geq 2$, our 
substitution replaces 
$t^k DXZ(S)$ by 
$$y^kt^kz_{j_1} \cdots z_{j_k} \prod_{i=1}^{k-1} 
\left(\frac{x_{j_{i+1}j_i}}{y} -1\right) = 
yt^k z_{j_1} \cdots z_{j_k} \prod_{i=1}^{k-1} 
(x_{j_{i+1}j_i} -y).$$
Thus if we let 
\begin{equation}
DXYZ(S) = \begin{cases} 
z_j & \mbox{if} \ S = \{j\}, \ \mbox{and} \\
z_{j_1} \cdots z_{j_k} \prod_{i=1}^{k-1} (x_{j_{i+1}j_i}-y) & 
\mbox{if} \ S =\{j_1 < \cdots < j_k\} \ \mbox{where} \ k \geq 2,
\end{cases}
\end{equation}
then we see that right-hand side of (\ref{PPDES2}) becomes
$$
\frac{y \sum_{n \geq 1} t^n \sum_{S \subseteq \PP, |S| =n} 
DXYZ(S)}{1 - y \sum_{n \geq 1} t^n \sum_{S \subseteq \PP, |S| =n} 
DXYZ(S)}.
$$
It follows that
\begin{equation*}
-x\sum_{w = w_1 \ldots w_n \in \PP^+} t^{|w|}y^{\wris(w)}\overline{z}^w \prod_{i \in Des(w)}x_{w_i w_{i+1}} =
\frac{-x\sum_{n \geq 1} t^n \sum_{S \subseteq \PP, |S| =n} 
DXYZ(S)}{1 - y \sum_{n \geq 1} t^n \sum_{S \subseteq \PP, |S| =n} 
DXYZ(S)}.
\end{equation*}
Thus 
\begin{equation}\label{PPDES3}
1-x\sum_{w = w_1 \ldots w_n \in \PP^+} t^{|w|}y^{\wris(w)}\overline{z}^w \prod_{i \in Des(w)}x_{w_i w_{i+1}} =
\frac{1-(x+y)\sum_{n \geq 1} t^n \sum_{S \subseteq \PP, |S| =n} 
DXYZ(S)}{1 - y \sum_{n \geq 1} t^n \sum_{S \subseteq \PP, |S| =n} 
DXYZ(S)}.
\end{equation}
By setting $z_i=0$ for $i > k$, we also 
obtain that 
\begin{equation}\label{kDES3}
1-x\sum_{w = w_1 \ldots w_n \in [k]^+} t^{|w|}y^{\wris(w)}\overline{z}^w \prod_{i \in Des(w)}x_{w_i w_{i+1}} =
\frac{1-(x+y)\sum_{n = 1}^{k} t^n \sum_{S \subseteq [k], |S| =n} 
DXYZ(S)}{1 - y \sum_{n = 1}^{k} t^n \sum_{S \subseteq [k], |S| =n} 
DXYZ(S)}.
\end{equation}

Note that if we replace $y$ by $(1-x)$ 
and $x_{ji}$ by $wt_u(ji)$, the left-hand side 
of (\ref{PPDES3}) becomes 
$U^{(\PP)}_u(x,{\bf z}_\infty,t)$ and the left-hand side 
of (\ref{kDES3}) becomes 
$U^{(k)}_u(x,{\bf z}_k,t)$. Similarly, if we replace $y$ by $(1-x)$ and 
$x_{ji}$ by $ewt_u(ji)$, the left-hand side 
of (\ref{PPDES3}) becomes 
$EU^{(\PP)}_u(x,{\bf z}_\infty,t)$ and the left-hand side 
of (\ref{kDES3}) becomes 
$EU^{(k)}_u(x,{\bf z}_k,t)$.
Then using the fact that 
$\mathcal{N}^ {(\PP)}_u(x,{\bf z}_\infty,t) = 
1/U^{(\PP)}_u(x,{\bf z}_\infty,t)$ and that 
$\mathcal{EN}^ {(\PP)}_u(x,{\bf z}_\infty,t) = 
1/EU^{(\PP)}_u(x,{\bf z}_\infty,t)$,  
we have the following 
theorem.

\begin{theorem}\label{thm:u1>ujPP} 
Suppose that $u = u_1 \ldots u_j \in \PP^*$, $\red(u) =u$, 
$\des(u) =1$, $u_1 > u_j$. Then 
\begin{equation} 
\mathcal{N}^ {(\PP)}_u(x,{\bf z}_\infty,t) = 
\frac{1 - (1-x) \sum_{n \geq 1} t^n \sum_{S \subseteq \PP, |S| =n} 
DXTZ_u(S)}{1-\sum_{n \geq 1} t^n \sum_{S \subseteq \PP, |S| =n} 
DXTZ_u(S)}
\end{equation}
and
\begin{equation}
\mathcal{EN}^ {(\PP)}_u(x,{\bf z}_\infty,t) = 
\frac{1 - (1-x) \sum_{n \geq 1} t^n \sum_{S \subseteq \PP, |S| =n} 
EDXTZ_u(S)}{1-\sum_{n \geq 1} t^n \sum_{S \subseteq \PP, |S| =n} 
EDXTZ_u(S)}
\end{equation}
where 
\begin{equation}
DXTZ_u(S) = \begin{cases} 
z_j & \mbox{if} \ S = \{j\}, \ \mbox{and} \\
z_{j_1} \cdots z_{j_k} \prod_{i=1}^{k-1} (wt_u(j_{i+1}j_i)+x-1) & 
\mbox{if} \ S =\{j_1 < \cdots < j_k\} \ \mbox{where} \ k \geq 2
\end{cases}
\end{equation}
and
\begin{equation}
EDXTZ_u(S) = \begin{cases} 
z_j & \mbox{if} \ S = \{j\}, \ \mbox{and} \\
z_{j_1} \cdots z_{j_k} \prod_{i=1}^{k-1} (ewt_u(j_{i+1}j_i)+x-1) & 
\mbox{if} \ S =\{j_1 < \cdots < j_k\} \ \mbox{where} \ k \geq 2.
\end{cases}
\end{equation}
\end{theorem}

If we set $z_i =0$ for 
all $i > k$, then we obtain the following theorem.

\begin{theorem}\label{thm:u1>ujk} 
Now suppose that $u = u_1 \ldots u_j \in [k]^*$, $\red(u) =u$, 
$\des(u) =1$, $u_1 > u_j$. Then 
\begin{equation} 
\mathcal{N}^ {(k)}_u(x,{\bf z}_k,t) = 
\frac{1 - (1-x) \sum_{n = 1}^{k} t^n \sum_{S \subseteq [k], |S| =n} 
DXTZ_u(S)}{1-\sum_{n = 1}^{k} t^n \sum_{S \subseteq [k], |S| =n} 
DXTZ_u(S)}
\end{equation}
and
\begin{equation}
\mathcal{EN}^ {(k)}_u(x,{\bf z}_k,t) = 
\frac{1 - (1-x) \sum_{n \geq 1} t^n \sum_{S \subseteq [k], |S| =n} 
EDXTZ_u(S)}{1-\sum_{n \geq 1} t^n \sum_{S \subseteq [k], |S| =n} 
EDXTZ_u(S)}.
\end{equation}
\end{theorem}

It follows from Theorem \ref{thm:u1>ujk} that 
to compute the generating function $\mathcal{N}^{(k)}_u(x,{\bf z}_k,t)$, we need 
only compute sums of the form 
$$P_{n,u}(x,t) = \sum_{S \subseteq [k], |S| =n} DXTZ_u(S)$$
for $1 \leq n \leq k$.
As an example, suppose that $k=5$ and we want to compute 
$\mathcal{N}^ {(5)}_{2341}(x,{\bf z}_5,t)$ where 
we set $z_i=1$ for all $i$.
Then with this specialization, it is easy to see that 
\begin{enumerate}
\item $wt_{2341}(21) = -3xt^2$,
\item $wt_{2341}(3i) = -xt^2$ for all $i <3$, 
\item $wt_{2341}(4i) = 0$ for all $i <4$, and 
\item $wt_{2341}(5i) = 0$ for all $i <5$.
\end{enumerate}
It follows that that 
\begin{enumerate}
\item $DXTZ_{2341}(\{1,2\})|_{z_i =1} = -3xt^2+x-1$,
\item $DXTZ_{2341}(\{i,3\})|_{z_i =1} = -xt^2+x-1$ for all $i <3$, 
\item $DXTZ_{2341}(\{i,4\})|_{z_i =1} = x-1$ for all $i <4$, and 
\item $DXTZ_{2341}(\{i,5\})|_{z_i =1} = x-1$ for all $i <5$.
\end{enumerate}
One can then compute that 
\begin{enumerate}
\item $P_{1,2341}(x,t)=5$,
\item $P_{2,2341}(x,t)=-10+10x-5xt^2$,
\item $P_{3,2341}(x,t) = 
10-20x+14t^2x+10x^2-14t^2x^2+3t^4x^2$, 
\item $P_{4,2341}(x,t)=-5+15x-13t^2x-15x^2+26t^2x^2-6t^4x^2+5x^3-13t^2x^3+6t^4x^3$, and 
\item $P_{5,2341}(x,t) = (-3xt^2+x-1)(-xt^2+x-1)(x-1)^2$.
\end{enumerate}

Thus 
\begin{equation}
\sum_{w \in [5]^*,2341\mbox{-mch}(w)=0} x^{\des(w)+1} = 
\frac{1-(1-x)(\sum_{n=1}^5 P_{n,2341}(x,t)t^n)}{1-\sum_{n=1}^5 P_{n,2341}(x,t)t^n}.
\end{equation}
We have computed that the initial terms of this series are 
\begin{align*}
& 1+5xt+5(3x+2x^2)t^2+5(7x+16x^2+2x^3)t^3 + 5(14x+72x^2+37x^3+x^4)t^4+\\
& (126x+1210x^2+1492x^3+246x^4+x^5)t^5+ 
(210x+3387x^2+7921x^3+3522x^4+210x^5)t^6 + \\
&(330x+8344x^2+32461x^3+28902x^4+5471x^5+120x^6)t^7+\cdots.
\end{align*}

One can obtain several interesting generating functions from 
$\mathcal{N}^ {(5)}_{2341}(x,{\bf z}_5,t)$. For example 
setting $x=0$ in $\frac{1}{2} \frac{\partial^2}{\partial x^2} \mathcal{N}^ {(5)}_{2341}(x,{\bf z}_5,t)$, one finds that the generating function 
for the number of words $w$ in $[5]^*$ such that 
$\des(w) =1$ and 2341-$\mathrm{mch}(w) =0$ is 
\begin{equation*}
\frac{t^2(10-20t+10t^2+10t^3-13t^4+4t^5)}{(1-t)^{10}}.
\end{equation*}
Similarly setting $x=0$ in $\frac{1}{6} \frac{\partial^3}{\partial x^3} \mathcal{N}^ {(5)}_{2341}(x,{\bf z}_5,t)$, one finds that the generating function 
for the number of words $w$ in $[5]^*$ such that 
$\des(w) =2$ and 2341-$\mathrm{mch}(w) =0$ is 
\begin{equation*}
\frac{t^3Q(t)}{(1-t)^{15}}
\end{equation*}
where 
$$Q(t) = 
10+35t-233t^2+ 416t^3-219t^4-266t^5+458t^6-167t^7-161t^8+198t^9-83t^{10}+13t^{11}.$$

If $u=2^s1$ where $s \geq 2$ and we set 
$z_i =1$ for all $i$, then it is easy to see 
that $wt_{2^s1}(ji) = -xt^{s-1}$ for all $j > i$ and that 
$DXTZ(S)|_{z_i =1} = (-xt^{s-1}+1-x)^{|S|-1}$ for all $S \subseteq [k]$ 
where $|S| \geq 1$. It then easily follows from 
Theorem \ref{thm:u1>ujk} 

\begin{equation}
\sum_{w \in [k]^*,2^s1\mbox{-mch}(w)=0} x^{\des(w)+1} = 
\frac{1-(1-x)(\sum_{n=1}^k \binom{k}{n} (-xt^{s-1}+1-x)^{n-1}t^n)}
{1-\sum_{n=1}^k \binom{k}{n} (-xt^{s-1}+1-x)^{n-1} t^n}.
\end{equation}

As an example, 
\begin{equation}
\sum_{w \in [5]^*,2^31\mbox{-mch}(w)=0} x^{\des(w)+1} = 
\frac{1-(1-x)(\sum_{n=1}^5 \binom{5}{n} (-xt^2+x-1)^{n-1}t^n)}{1-\sum_{n=1}^k \binom{5}{n} (-xt^2+x-1)^{n-1} t^n}.
\end{equation}
We have computed that the initial terms of this series are 
\begin{align*}
& 1+5xt+5(3x+2x^2)t^2+5(7x+16x^2+2x^3)t^3 + 5(14x+71x^2+37x^3+x^4)t^4+\\
& (126x+1166x^2+1486x^3+246x^4+x^5)t^5+ 
5(42x+634x^2+1553x^3+704x^4+42x^5)t^6 + \\
&(330x+7554x^2+30998x^3+28662x^4+5471x^5+120x^6)t^7+O[t]^8.
\end{align*}

\section{The case $u = u_1 \ldots u_j$, $\des(u)=1$, and $u_1 < u_j$}

In this section, we shall consider the problem 
of computing the generating functions \\
$\mathcal{N}^{(\PP)}_u(x,{\bf z}_\infty,t)$, 
$\mathcal{N}^{(k)}_u(x,{\bf z}_k,t)$, $\mathcal{EN}^{(\PP)}_u (x,{\bf z}_\infty,t)$, and 
$\mathcal{EN}^{(k)}_u(x,{\bf z}_k,t)$ for  
$u = u_1 \ldots u_j$ such that $\des(u)=1$, $u_1 < u_j$, 
and $u$ has the $\PP$-weakly increasing overlapping property 
($[k]$-weakly increasing overlapping property).

\begin{figure}[htbp]
  \begin{center}
    \includegraphics[width=0.6\textwidth]{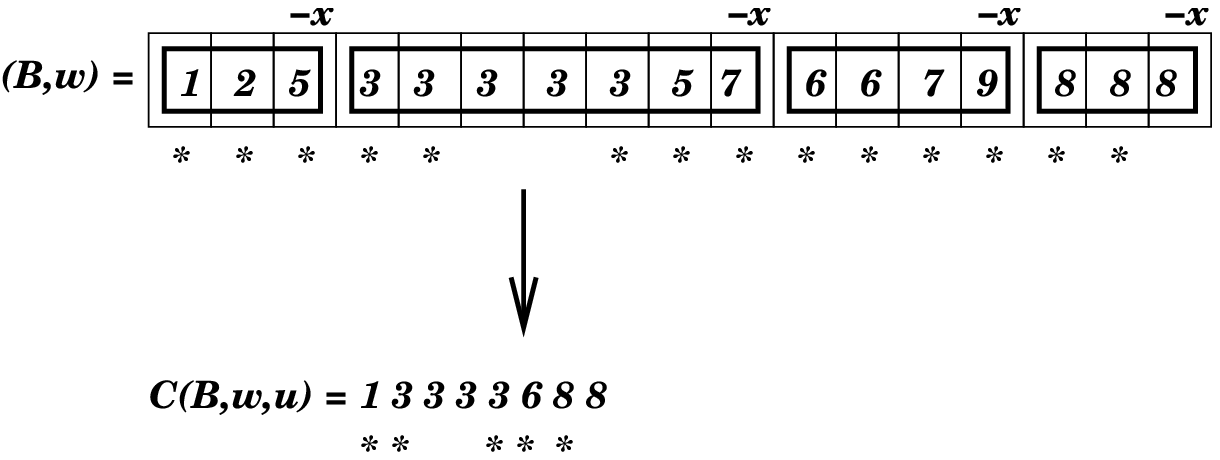}
    \caption{A fixed point of $I_{12433}$.}
    \label{fig:Collapse12433}
  \end{center}
\end{figure}

Again the simplest case is when $u$ has the $\PP$-minimal 
overlapping property in which case $u$ automatically has 
the $\PP$-weakly increasing overlapping property. 
For example, suppose that $u = 12433$. 
Now suppose that we are given a fixed point $(B,w)$ of $I_u$, where 
$B=(b_1, \ldots, b_k)$ and $w=w_1 \ldots w_n$, such as the one pictured 
in Figure \ref{fig:Collapse12433}.  We know that to be a fixed point of $I_u$, 
$w$ must be  weakly increasing within bricks of $B$ and that for any $i <k$, 
if $c$ is last cell in brick $b_i$ and $w_c > w_{c+1}$, then there must 
be a $u$-match in $w$ which is contained in the cells of $b_i$ and 
$b_{i+1}$.  In our particular example, since $u=12433$ has a single descent, 
this match must involve the last three cells of $b_i$ and the first two cells 
of $b_{i+1}$.  In Figure \ref{fig:Collapse12433}, we have indicated 
the three such matches in our example by placing stars below the cells in the 
$12433$-matches. In this case, the collapse map just maps $(B,w)$ to the 
word $v =C(B,w,u)$ which is the result of starting with 
$w$ and removing the letters in all such matches that do not 
correspond to the end points of the match.  This process is pictured 
in Figure \ref{fig:Collapse12433} where again we have starred the elements 
in $C(B,w,u)$ that remain from the original $12433$-matches in $w$. In this 
case, the resulting word $C(B,w,u)$ must be weakly increasing.

As in the previous section, we want to 
construct the set of fixed points of $(B,w)$ of $I_u$ such that 
$C(B,w,v)$ is equal to a given word $v=v_1 \ldots v_n$ 
where $v_1 \leq \cdots \leq v_n$.

If $v_s < v_{s+1}$, then we have three possibilities: 
(i) $v_sv_{s+1}$ could lie in the same brick $b_i$ of $B$, 
(ii) $v_s$ could end a brick $b_i$ and $v_{s+1}$ could 
start the brick $b_{i+1}$ in $B$, or (iii) 
$v_sv_{s+1}$ arose from a collapse across two 
bricks $b_i$ and $b_{i+1}$ where there was a decrease 
between bricks $b_i$ and $b_{i+1}$ and 
$v_s$ played the role of 1 in the $u$-match and 
$v_{s+1}$ plays the role of second 3 in the $u$-match 
that must cross the bricks $b_i$ and $b_{i+1}$. 
 For example, 
suppose that the underlying alphabet is $[9]$. If $v_s =8$ and $v_{s+1} =9$, then $v$ 
could not have come from 
the collapse of $12433$-match because we can not add a letter which 
could play the role of 4  in the 12433-match.  Hence, the  
weight associated to a rise 89 is just $1-x$. If 
we consider  the first rise  $13$ in the $C(B,w,u)$ of 
Figure \ref{fig:Collapse12433}, then we see there are many ways that we 
could add the three letters middle letters. That is, the 
original 12433-match could have been any $12c33$ where 
$c \in \{4,5,6,7,8,9\}$. It follows 
that the extra weight from these possibilities that 
is not included in $\overline{z}^{C(B,w,u)}t^{|C(B,w,u)|}$ 
is $-x t^3 z_2z_3\sum_{4 \leq c \leq 9} z_c$. Here the $-x$ comes from 
the fact that we know that the original match straddled two bricks 
and there is a weight of $-x$ associated with the end point of 
the first of those two bricks. Thus the weight 
associated with the rise $13$ is 
$1-x -x t^3 z_2z_3\sum_{4 \leq c \leq 9} z_c$. If 
we consider the second rise  $36$ in the $C(B,w,u)$ of 
Figure \ref{fig:Collapse12433}, then we see there are many ways that we 
could add the three letters middle letters. That is, the 
original 12433-match could have been any $3cd66$ where 
$c \in \{4,5\}$ and $d \in \{7,8,9\}$. It follows 
that the extra weight from these possibilities that 
is not included in $\overline{z}^{C(B,w,u)}t^{|C(B,w,u)|}$ in this case 
is $-x t^3 z_6(z_4+z_5)(z_7+z_8+z_9)$. Thus 
the weight associated with the rise 36 is 
$1-x-x t^3 z_6(z_4+z_5)(z_7+z_8+z_9)$. 
Finally, 
we consider the third rise  $68$ in the $C(B,w,u)$ of 
Figure \ref{fig:Collapse12433}, then there is only one 
way to add the three middle letters. That is, the 
original 12433-match must have been $67988$. It follows 
that the extra weight in this case that 
is not included in $\overline{z}^{C(B,w,u)}t^{|C(B,w,u)|}$  
would be $-x t^3 z_7z_8z_9$. Thus the weight associated 
with the final rise 68 is  
$1-x-x t^3 z_7z_8z_9$.
On the other hand, if $v_s =  v_{s+1}$, 
then we have only two choices. That is, either cell $s$ was the end of a brick 
or cell $s$ was an internal cell of a brick.  This implies that 
each level in $v$ contributes a factor of $(1-x)$ since if 
$s$ is at the end of a brick, there is a weight of $-x$ associated with 
the last cell of a brick. In this way, we can associate a 
weight with each level  or rise of $v$ which 
will allows to compute
$$\sum_{\overset{(B,w) \ \mathrm{is\ a\ fixed\ point\ of} \ I_u}{C(B,w,u) =v}} sgn(B,w) wt(B,w).$$

In our case where $u =12433$ and $k=9$, 
the weights associated with the rises are  
given in the table below.

\begin{center}
\begin{tabular}{|l|l|}
\hline
Rises & $wt_{12433,9}(ij)$\\
\hline
$i9\  (i \leq 7)$ or $i (i+1)$ & $1-x$ \\
\hline
$i8\  (i \leq 6)$ & $1-x -xt^3z_8z_9(\sum_{i < j < 8} z_j)$\\
\hline 
$i7 \ (i \leq 5)$ & $1-x-xt^3 z_7(z_8+z_9) (\sum_{i < j < 7} z_j)$ \\
\hline
$i6 \ (i \leq 4)$ &  $1-x-xt^3 z_6(z_7+z_8+z_9) (\sum_{i < j < 6} z_j)$\\
\hline
$i5 \ (i \leq 3)$ & $1-x-xt^3 z_5(z_6+z_7+z_8+z_9) (\sum_{i < j < 5} z_j)$ \\
\hline 
$i4 \ (i \leq 2)$ & $1-x-xt^3 z_4(\sum_{4 < s \leq 9} z_s)(\sum_{i < j < 4} z_j)$\\
\hline
$13$ & $1-x-xt^3 z_2z_3(\sum_{3 < s \leq 9} z_s)$\\
\hline
\end{tabular}\\
The weights $wt_{12433,9}(ij)$.
\end{center}

However, if $u=12433$ and we want to compute 
$U^{(\PP)}_{u,n}(x,{\bf z}_\infty)$, the 
weights for any rise $ij$ where $i+1 < j$ would be 
$1-x -xt^3z_j (\sum_{i <s < j} z_s)(\sum_{j <d} z_d)$ which 
is an infinite sum. 

Going back to our example where $u =12433$ and $k =9$, 
it follows 
that for any $v \in [9]^+$, 
\begin{equation}\label{eq:v912333}
\sum_{\overset{(B,w) \ \mathrm{is\ a\ fixed\ point\ of} \ I_u}{C(B,w,u) =v}}
sgn(B,w) wt(B,w) =
-x \overline{z}^v t^{|v|} (1-x)^{\lev(v)} 
\prod_{s \in Rise(v)} wt_{12433,9}(v_s v_{s+1}).
\end{equation}
As in the previous section, then 
initial $-x$ comes from the fact that the last cell of $(B,w)$ always 
contributes a $-x$ since the last cell is at the end of a brick. 
But then we know that
\begin{eqnarray}\label{eq:U912433}
U^{(9)}_{12433}(x,{\bf z}_9,t) &=& 1 + \sum_{n \geq 1} U^{(9)}_{12433,n}(x,{\bf z}_9) t^n \nonumber \\
&=& 1 + \sum_{v \in [9]^+,\des(v) =0} -x(1-x)^{\lev(v)} \overline{z}^v t^{|v|}
\prod_{s \in Rise(v)} wt_{12433,9}(v_s v_{s+1}).
\end{eqnarray}
Hence we could compute 
$\displaystyle \mathcal{N}^{(9)}_{12433}(x,{\bf z}_9,t) = 
\frac{1}{U^{(9)}_{12433}(x,{\bf z}_9,t)}$
if we can compute the right-hand side of (\ref{eq:U912433})

As in the previous section, the case of exact matches 
is much simpler. In that case, we want to compute 
$$\sum_{\overset{(B,w) \ \mathrm{is \ a \ fixed \ point \ of} \ J_u}{C(B,w,u) =v}} sgn(B,w) wt(B,w).$$
Going back to our example of $u=12433$ over the alphabet 
$[9]$, we see that the weight associated to a rise $v_s < v_{s+1}$ 
is $1-x$ unless $v_s=1$, $v_{s+1}=3$.  If $v_s=1$, $v_{s+1}=3$, then we  
must have eliminated a 243 from $w$.  Thus if we want to compute 
$EU^{(\PP)}_{12433,n}(x,{\bf z}_\infty)$ or 
$EU^{(k)}_{12433,n}(x,{\bf z}_k)$ for $k \geq 4$, the 
weights associated to rises are given in the following table. 
\begin{center}
\begin{tabular}{|l|l|}
\hline
Rise & weight $ewt_{12433,\PP}(ij)$\\
\hline
$ij\  \mbox{where either} \ i \neq 1 \ \mbox{or} \ j \neq 3$  & $1-x$ \\
\hline 
13 & $1-x -x z_2z_3z_4t^3$ \\
\hline
\end{tabular}\\
The weights $ewt_{12433}(ij)$.
\end{center}
It follows 
that for any $v \in [9]^+$, 
\begin{equation}\label{eq:ev912433}
\sum_{\overset{(B,w) \ \mathrm{is \ a \ fixed \ point\  of} \ J_{12433}}{C(B,w,12433) =v}} sgn(B,w) wt(B,w) =
-x\overline{z}^v t^{|v|}(1-x)^{\lev(v)} \prod_{s \in Rise(v)} ewt_{12433,9}(v_s v_{s+1})\end{equation}
and
\begin{eqnarray}\label{eq:EU912433}
EU^{(9)}_{12433}(x,{\bf z}_9,t) &=& 1 + \sum_{n \geq 1} EU^{(9)}_{12433,n}(x,{\bf z}_9) t^n \nonumber \\
&=& 1 + \sum_{v \in [9]^+,\des(v) =0} -x\overline{z}^v t^{|v|}(1-x)^{\lev(v)} \prod_{s \in Rise(v)} ewt_{12433,9}(v_s v_{s+1}).
\end{eqnarray}

When $u$ does not have the minimal overlapping property, 
we can obtain similar results if $u$ has the $\PP$-weakly increasing 
overlapping property or the $[k]$-weakly increasing 
overlapping property. For example suppose that 
$u =u_1, \ldots, u_j$, $\des(u) =1$, $u_1 < u_j$, and $u$ 
has the $\PP$-weakly increasing overlapping property.  
Now suppose that 
$w = w_1 \ldots w_n$ is a maximal sequence of linked $u$-matches. 
That is, we assume $w$ starts and ends with a $u$-match 
and any two consecutive $u$-matches share at least two letters. 
Then if the $u$-matches in $w$ start at positions 
$1=i_1 < i_2 < \cdots < i_k$, then the $\PP$-weakly increasing 
overlapping property in $w$ ensures that 
$w_1=w_{i_1} \leq \cdots \leq w_{i_k} < w_n$. Thus in a collapse map, 
if  we eliminate $w_2 \ldots w_{n-1}$ we will be left with a rise 
$w_1 w_n$.  This may not happen if $u$ does not have the 
$\PP$-weakly increasing overlapping property. For example, 
suppose $u =2413$, then the words $w^{(1)}=472613$, $w^{(2)}=472614$, 
and $w^{(3)}=472615$ have $u$-matches starting at positions 1 and 3. 
Thus in such a case, we have no control over the relationship 
between first and last letter of a maximal sequence of linked $u$-matches.

\begin{figure}[htbp]
  \begin{center}
    \includegraphics[width=0.6\textwidth]{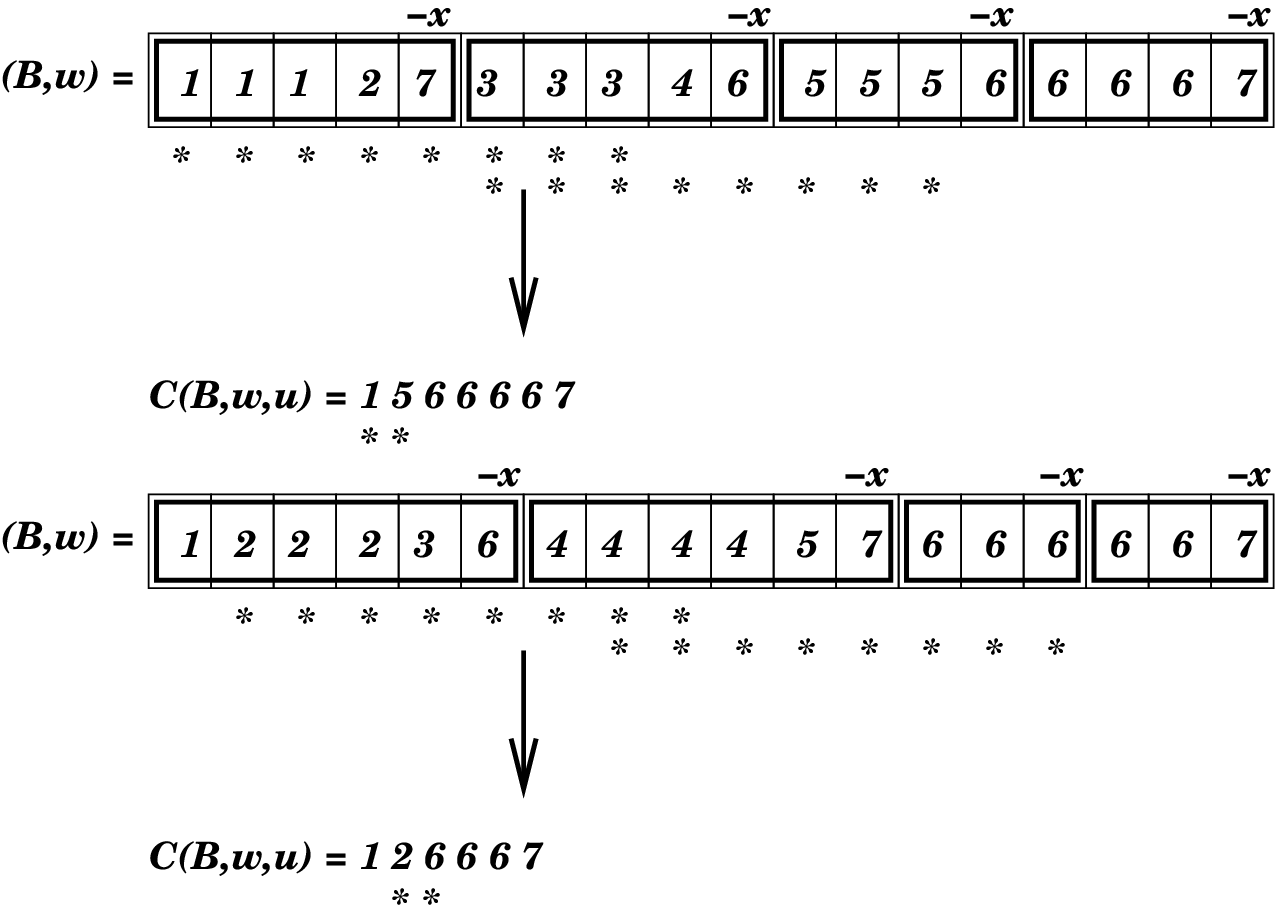}
    \caption{A fixed point of $I_{11124333}$.}
    \label{fig:Collapse11124333}
  \end{center}
\end{figure}

Thus assume that $u =u_1 \ldots u_j$, $\des(u) =1$, $u_1 < u_j$ and 
$u$ has the $\PP$-weakly increasing overlapping property. 
Then we shall see that the collapse map still works but 
the weight function $wt_u(ij)$ is more complicated. 
As we saw in 
the previous section, we must pay attention to 
overlapping $u$-matches that share more than  
one letter.  We will consider the example where 
$u =11124333$ and $k=7$. Clearly $u$ has the weakly increasing overlapping 
property. In this case, $u$-matches can 
overlap in either one, two, or three letters. As in 
the previous section, the collapse map will keep only 
the first and last letters of a consecutive sequence 
of $u$-matches such that each consecutive pair 
share at least two letters. For example, 
at the top of Figure \ref{fig:Collapse11124333}, 
we have given an example where two consecutive 
$u$-matches share 3 letters and at the bottom 
of Figure \ref{fig:Collapse11124333}, 
we have given an example where two consecutive 
$u$-matches share 2 letters.

As before, if we are given a weakly increasing word 
$v = v_1 \ldots v_n \in [7]^+$, 
we want to find the sum of the weights of all fixed points 
$(B,w)$ of $I_u$ such that  
$C(B,w,u) =v$.  Now if $v_s =v_{s+1}$, then either 
$v_sv_{s+1}$ lie in the same brick which contributes a factor of 1 
or $v_sv_{s+1}$ lie in different bricks which contributes a 
factor of $-x$ for the brick that ends at $v_s$.  Thus 
we obtain a factor of $1-x$ for each level of $v$. 
For the rises of $v$, we should observe that 
the start and the end of any two consecutive $u$-matches 
which share more than one letter must differ by at least 4. 
Similarly, the start and the end of any three consecutive $u$-matches 
in which each two consecutive $u$-matches 
share more than one letter must differ by at least 6. 
Hence, for $k=7$, we can not have three consecutive $u$-matches 
in which each two consecutive $u$-matches 
share more that one letter because the smallest starting point 
is 1 the smallest ending point is 7 which leaves no room 
for a letter which is larger than the last three letters in 
such a sequence.  For each pair, $v_s< v_{s+1}$ which occurs 
in $v$, we get a factor of $1-x$ as we did for levels. However 
in this case, we must also consider the possible collapses that 
could give rise to $v_sv_{s+1}$.  These 
are as follows. 
\begin{enumerate}
\item Rises of the form $i(i+1)$ or $i7$ where $1 \leq i \leq 5$ can not 
arise from the collapse map in our case so that 
$wt_{11124333,7}(v_sv_{s+1}) =1-x$ in these cases.

\item $v_sv_{s+1} =13$. In this case, a $u$-match that 
could give  rise to $13$ under the collapse map 
must be of the form $1112a333$ where 
$a \in \{4,5,6,7\}$. Thus 
$$wt_{11124333,7}(v_sv_{s+1})= 1-x-xt^6z_1^2z_2(z_4+z_5+z_6+z_7)z_3^2.$$

\item $v_sv_{s+1} =14$. In this case, a $u$-match that 
could give rise to $14$ under the collapse map 
must be of the form $111ab444$ where 
$a \in \{2,3\}$ and $b \in \{5,6,7\}$. Thus 
$$wt_{11124333,7}(v_sv_{s+1})= 1-x-xt^6z_1^2(z_2+z_3)(z_5+z_6+z_7)z_4^2.$$

\item $v_sv_{s+1} =15$. In this case, a single $u$-match that 
could give rise to $15$ under the collapse map 
must be of the form $111ab555$ where 
$a \in \{2,3,4\}$ and $b \in \{6,7\}$. There are 
also two possibilities for linked $u$-matches that 
could give rise to $15$ under the collapse map, namely, 
(i) $1112a3334b555$ or (ii) $1112a33334b555$ where 
$a \in \{4,5,6,7\}$ and $b \in \{6,7\}$. Thus  
\begin{eqnarray*}
wt_{11124333,7}(v_sv_{s+1})&=& 1-x-xt^6z_1^2(z_2+z_3+z_4)(z_6+z_7)z_5^2-\\
&&xt^{11}z_1^2z_2(z_4+z_5+z_6+z_7)z_3^3z_4(z_6+z_7)z_5^2 - \\
&&xt^{12}z_1^2z_2(z_4+z_5+z_6+z_7)z_3^4z_4(z_6+z_7)z_5^2.
\end{eqnarray*}

\item $v_sv_{s+1} =16$. In this case, a single $u$-match that 
could give rise to $16$ under the collapse map 
must be of the form $111a7666$ where 
$a \in \{2,3,4,5\}$. There are 
also four possibilities for linked $u$-matches that 
could give rise to $16$ under the collapse map, namely, \\
(i) $1112a333b7666$ $a \in \{4,5,6,7\}$ and $b \in \{4,5\}$, 
(ii) $1112a3333b7666$ where $a \in \{4,5,6,7\}$ and $b \in \{4,5\}$, 
(iii) $111ab44457666$ $a \in \{2,3\}$ and $b \in \{5,6,7\}$, or  
(iv) $111ab444457666$ where $a \in \{2,3\}$ and $b \in \{5,6,7\}$. 
Thus  
\begin{eqnarray*}
wt_{11124333,7}(v_sv_{s+1})&=& 1-x-xt^6z_1^2(z_2+z_3+z_4+z_5)z_7z_6^2-\\
&&xt^{11}z_1^2z_2(z_4+z_5+z_6+z_7)z_3^3(z_4+z_5)z_7z_6^2 - \\
&&xt^{12}z_1^2z_2(z_4+z_5+z_6+z_7)z_3^4(z_4+z_5)z_7z_6^2 - \\
&&xt^{11}z_1^2(z_2+z_3)(z_5+z_6+z_7)z_4^3z_5z_7z_6^2 - \\
&&xt^{12}z_1^2(z_2+z_3)(z_5+z_6+z_7)z_4^4z_5z_7z_6^2.
\end{eqnarray*}

\item $v_sv_{s+1} =24$. In this case, a $u$-match that 
could give  rise to $24$ under the collapse map 
must be of the form $2223a444$  
where $a \in \{5,6,7\}$. Thus  
$$wt_{11124333,7}(v_sv_{s+1})= 1-x-xt^6z_2^2z_3(z_5+z_6+z_7)z_4^2.$$

\item $v_sv_{s+1} =25$. In this case, a $u$-match that 
could give rise to $25$ under the collapse map 
must be of the form $222ab555$ where 
$a \in \{3,4\}$ and $b \in \{6,7\}$. Thus  
$$wt_{11124333,7}(v_sv_{s+1})= 1-x-xt^6z_2^2(z_3+z_4)(z_6+z_7)z_5^2.$$

\item $v_sv_{s+1} =26$. In this case, a single $u$-match that 
could give rise to $26$ under the collapse map 
must be of the form $222a7666$ where 
$a \in \{3,4,5\}$. There are 
also two possibilities for linked $u$-matches that 
could give rise to $26$ under the collapse map, namely, 
(i) $2223a44457666$ or (ii) $2223a444457666$ where 
$a \in \{5,6,7\}$.  Thus 
\begin{eqnarray*}
wt_{11124333,7}(v_sv_{s+1})&=& 1-x-xt^6z_2^2(z_3+z_4+z_5)z_7z_6^2-\\
&&xt^{11}z_2^2z_3(z_5+z_6+z_7)z_4^3z_5z_7z_6^2 - \\
&&xt^{12}z_2^2z_3(z_5+z_6+z_7)z_4^4z_5z_7z_6^2.
\end{eqnarray*}

\item $v_sv_{s+1} =35$. In this case, a $u$-match that 
could give rise to $35$ under the collapse map 
must be of the form $3334a555$ where 
$a \in \{6,7\}$. Thus 
$$wt_{11124333,7}(v_sv_{s+1})= 1-x-xt^6z_3^2z_4(z_6+z_7)z_5^2.$$

\item $v_sv_{s+1} =36$. In this case, a $u$-match that 
could give rise to $36$ under the collapse map 
must be of the form $333a7666$ where 
$a \in \{4,5\}$. Thus 
$$wt_{11124333,7}(v_sv_{s+1})= 1-x-xt^6z_3^2(z_4+z_5)z_7z_6^2.$$

\item $v_sv_{s+1} =46$. In this case, a $u$-match that 
could give rise to $46$ under the collapse map 
must be of the form $44457666$. Thus  
$$wt_{11124333,7}(v_sv_{s+1})= 1-x-xt^6z_4^2z_5z_7z_6^2.$$

\end{enumerate}

It follows 
that for any $v \in [7]^+$ such that $v$ is weakly increasing, 
\begin{multline}\label{eq:v711124333}
\sum_{\overset{(B,w) \ \mathrm{is \ a \ fixed \ point \ of} \ I_{11124333}}{C(B,w,11124333) =v}} sgn(B,w) wt_{11124333,7}(B,w) = \\
-x\overline{z}^v t^{|v|}(1-x)^{\lev(v)} \prod_{s \in Rise(v)} wt_{11124333,7}(v_s v_{s+1}).
\end{multline}
and
\begin{multline}\label{eq:U711124333}
U^{(7)}_{11124333}(x,{\bf z}_7,t) = 1 + \sum_{n \geq 1} U^{(7)}_{11124333,n}(x,{\bf z}_7) t^n =  \\
1 + \sum_{v \in [7]^+, \des(v) =0} -x\overline{z}^v t^{|v|}(1-x)^{\lev(v)} \prod_{s \in Rise(v)} 
wt_{11124333,7}(v_s v_{s+1}).
\end{multline}

What we need to be able to compute the right-hand sides of either 
(\ref{eq:U912433}), (\ref{eq:EU912433}), or (\ref{eq:U711124333}), 
is the generating function over all weakly 
increasing words $v \in \PP^*$ where we not only keep track of the rises of $P$ but also 
the type of rises. 

By Theorem \ref{thm:labris}, 
we know that 
\begin{equation}\label{PPRIS1}
\sum_{n \geq 1} t^n \sum_{S \subseteq \PP, |S| =n} 
RXZ(S) = \sum_{w = w_1 \leq \cdots \leq w_n \in \PP^+} 
t^{|w|}\overline{z}^w \prod_{i \in Rise(w)}
x_{w_i w_{i+1}}.
\end{equation}
If we first replace $t$ by $yt$ and $x_{ij}$ by $x_{ij}/y$ in (\ref{PPRIS1}) 
and then divide by $y$,
the right-hand side (\ref{PPRIS1}) becomes 
$$\sum_{w = w_1 \leq \cdots \leq w_n \in \PP^+} 
t^{|w|}\overline{z}^w y^{\lev(w)} \prod_{i \in Rise(w)}
x_{w_i w_{i+1}} $$
and the left-hand side 
becomes 
$$\sum_{n \geq 1} t^n \sum_{S \subseteq \PP, |S| =n} 
RXYZ(S)$$
where 
\begin{equation}
RXYZ(S) = \begin{cases} 
\frac{z_j}{1-z_jyt} & \mbox{if} \ S = \{j\}, \ \mbox{and} \\
\left(\prod_{i=1}^k\frac{z_{j_i}}{1-z_{j_i}yt}\right) \prod_{i=1}^{k-1} x_{j_i j_{i+1}} & 
\mbox{if} \ S =\{j_1 < \cdots < j_k\} \ \mbox{where} \ k \geq 2.
\end{cases}
\end{equation}
Hence 
\begin{multline}\label{PPRIS2}
1-x\sum_{w = w_1 \leq \cdots \leq w_n \in \PP^+} 
t^{|w|}y^{\lev(w)}\overline{z}^w \prod_{i \in Rise(w)} x_{w_i w_{i+1}} = \\
1-x \sum_{n \geq 1} t^n \sum_{S \subseteq \PP, |S| =n} 
RXYZ(S). 
\end{multline}
If we set $z_i=0$ for $i > k$, then we obtain that  
\begin{multline}\label{kRIS3}
1-x\sum_{w = w_1 \leq \cdots \leq w_n \in [k]^+} t^{|w|}y^{\lev(w)}\overline{z}^w \prod_{i \in Rise(w)}x_{w_i w_{i+1}} = \\
1-x \sum_{n = 1}^{k} t^n \sum_{S \subseteq [k], |S| =n} 
RXYZ(S). 
\end{multline}
Note that if we replace $y$ by $(1-x)$ 
and $x_{ij}$ by $wt_u(ij)$, the left-hand side 
of (\ref{PPRIS2}) becomes 
$U^{(\PP)}_u(x,{\bf z}_\infty,t)$ and the left-hand side 
of (\ref{kRIS3}) becomes 
$U^{(k)}_u(x,{\bf z}_k,t)$. Similarly, if we replace $y$ by $(1-x)$ and 
$x_{ij}$ by $ewt_u(ij)$, the left-hand side 
of (\ref{PPRIS2}) becomes 
$EU^{(\PP)}_u(x,{\bf z}_\infty,t)$ and the left-hand side 
of (\ref{kRIS3}) becomes 
$EU^{(k)}_u(x,{\bf z}_k,t)$.
Then using the fact that 
$\mathcal{N}^ {(\PP)}_u(x,{\bf z}_\infty,t) = 
1/U^{(\PP)}_u(x,{\bf z}_\infty,t)$ and that 
$\mathcal{EN}^ {(\PP)}_u(x,{\bf z}_\infty,t) = 
1/EU^{(\PP)}_u(x,{\bf z}_\infty,t)$,  
we have the following 
theorem.

\begin{theorem}\label{thm:u1<ujPP} 
Suppose that $u = u_1 \ldots u_j \in \PP^*$, $\red(u) =u$, 
$\des(u) =1$, $u_1 < u_j$, and $u$ has the $\PP$-weakly increasing 
overlapping property. Then 
\begin{equation} 
\mathcal{N}^ {(\PP)}_u(x,{\bf z}_\infty,t) = 
\frac{1}{1-x \sum_{n \geq 1} t^n \sum_{S \subseteq \PP, |S| =n} 
RXTZ(S)}
\end{equation}
and
\begin{equation}
\mathcal{EN}^ {(\PP)}_u(x,{\bf z}_\infty,t) = 
\frac{1}{1-x \sum_{n \geq 1} t^n \sum_{S \subseteq \PP, |S| =n} 
ERXTZ(S)}
\end{equation}
where 
\begin{equation}
RXTZ_u(S) = \begin{cases} 
\frac{z_j}{1-(1-x)z_jt} & \mbox{if} \ S = \{j\}, \ \mbox{and} \\
\left( \prod_{i=1}^k \frac{z_{j_i}}{1-(1-x)z_{j_i}t}\right)  
\prod_{i=1}^{k-1} (wt_u(j_ij_{i+1})) & 
\mbox{if} \ S =\{j_1 < \cdots < j_k\} \ \mbox{where} \ k \geq 2.
\end{cases}
\end{equation}
and
\begin{equation}
ERXTZ_u(S) = \begin{cases} 
\frac{z_j}{1-(1-x)z_jt} & \mbox{if} \ S = \{j\}, \ \mbox{and} \\
\left(  \prod_{i=1}^k \frac{z_{j_i}}{1-(1-x)z_{j_i}t}\right)  
\prod_{i=1}^{k-1} ewt_u(j_ij_{i+1}) & 
\mbox{if} \ S =\{j_1 < \cdots < j_k\} \ \mbox{where} \ k \geq 2.
\end{cases}
\end{equation}
\end{theorem}

If we specialize the variables so that $z_i =0$ for 
all $i > k$, then we have the following theorem.

\begin{theorem}\label{thm:u1<ujk} 
Suppose that $u = u_1 \ldots u_j \in [k]^*$, $\red(u) =u$, 
$\des(u) =1$, $u_1 < u_j$, and $u$ has the $[k]$-weakly increasing 
overlapping property. Then 
\begin{equation} 
\mathcal{N}^ {(k)}_u(x,{\bf z}_k,t) = 
\frac{1}{1-x \sum_{n =1}^k t^n \sum_{S \subseteq [k], |S| =n} 
RXTZ(S)}
\end{equation}
and
\begin{equation}
\mathcal{EN}^ {(k)}_u(x,{\bf z}_k,t) = 
\frac{1}{1-x \sum_{n=1}^k t^n \sum_{S \subseteq [k], |S| =n} 
ERXTZ(S)}.
\end{equation}
\end{theorem}

It follows from Theorem \ref{thm:u1<ujk} that 
to compute the generating function we need 
to $\mathcal{N}^ {(k)}_u(x,{\bf z}_k,t)$, we need 
only compute sums of the form 
$$P_{n,u}(x,t) = \sum_{S \subseteq [k], |S| =n} RXTZ_u(S)$$
for $1 \leq n \leq k$ and that 
to compute the generating function we need 
to $\mathcal{EN}^ {(k)}_u(x,{\bf z}_k,t)$, we need 
only compute sums of the form 
$$P_{n,u}(x,t) = \sum_{S \subseteq [k], |S| =n} ERXTZ_u(S)$$
for $1 \leq n \leq k$.

For example, suppose that we want to compute 
$\mathcal{EN}^ {(9)}_u(x,{\bf z}_9,t)$ where 
$u =12433$ and we set $z_i =1$ for $i =1, \ldots , 9$. 
For each set singleton $S=\{j\}$, $ERXTZ_u(S) = \frac{1}{(1-(1-x)t)}$. 
For sets $S$ of cardinality greater than 2, there  
are two types of sets $S = \{j_1 < j_2 < \ldots j_n\}$ to consider, 
namely, those where $j_1 =1$ and $j_2 =3$ and those sets 
where it is not the case that $j_1 =1$ and $j_2 =3$. 
If  $S = \{j_1 < j_2 < \ldots j_n\}$ 
where it is not the case that $j_1 =1$ and $j_2 =3$, 
then we know that $ERXTZ_u(S) =\frac{(1-x)^{k-1}}{(1-(1-x)t)^k}$. 
If $S$ is of the form $\{1,3\}\cup T$ where $T \subseteq \{4,5,6,7,8,9\}$, 
then 
$$ERXTZ_u(S) =
(1-x-xt^3)(1-x)^{|T|}\frac{1}{(1-(1-x)t)^{|T|+2}}.$$
If follows that 
\begin{eqnarray*}
\sum_{n=1}^9 t^n \sum_{S \subseteq [9], |S| =n} ERXTZ(S) 
 &=& \sum_{n=1}^p \binom{9}{k} \frac{t^k(1-x)^{k-1}}{(1-(1-x)t)^k} - \\
&&\sum_{j=0}^6 \binom{6}{j} \frac{t^{j+2}(1-x)^{j+1}}{(1-(1-x)t)^{j+2}} + \\
&& \sum_{j=0}^6 \binom{6}{j} 
\frac{t^{j+2}(1-x-xt^3)(1-x)^{j}}{(1-(1-x)t)^{j+2}} \\
&=&  
\sum_{n=1}^p \binom{9}{k} \frac{t^k(1-x)^{k-1}}{(1-(1-x)t)^k}-  \\
&&  \sum_{j=0}^6 \binom{6}{j} 
\frac{t^{j+2}(xt^3)(1-x)^{j}}{(1-(1-x)t)^{j+2}}. 
\end{eqnarray*}
Thus if we let 
$$A_{12433,9}(x,t) = 1-x\left(  
\sum_{k=1}^p \binom{9}{k} \frac{t^k(1-x)^{k-1}}{(1-(1-x)t)^k}- 
\sum_{j=0}^6 \binom{6}{j} \frac{t^{j+2}(xt^3)(1-x)^{j}}{(1-(1-x)t)^{j+2}} 
\right),$$
then 
\begin{equation}\label{2eq:U912433}
\mathcal{EN}^ {(9)}_{12433}(x,{\bf z}_9,t)|_{z_i =1} = 
\frac{1}{A_{12433,9}(x,t)}.
\end{equation}
We have used (\ref{2eq:U912433}) to compute the first few terms 
in the series of $\mathcal{EN}^ {(9)}_{12433}(x,{\bf z}_\infty,t)|_{z_i =1}$.
\begin{eqnarray*}
&&\mathcal{EN}^ {(9)}_{12433}(x,{\bf z}_9,t)|_{z_i =1} = \\
&&1+9 t x+t^2 \left(45 x+36 x^2\right)+t^3 \left(165 x+480 x^2+84 x^3\right)+\\
&& t^4 \left(495 x+3510 x^2+2430 x^3+126 x^4\right)+\\
&&t^5 \left(1287
x+18612 x^2+31212 x^3+7812 x^4+126 x^5\right)+\\
&&t^6 \left(3003 x+79925 x^2+262626 x^3+167826 x^4+17976 x^5+84 x^6\right)+\\
&&t^7 \left(6435 x+294616 x^2+1683386
x^3+2132496 x^4+634446 x^5+31536 x^6+36 x^7\right)+\\
&&t^8 \left(12870 x+965709 x^2+8885187 x^3+19458252 x^4+11854197 x^5+
\right. \\
&&\left. 1826577 x^6+43677 x^7+9 x^8\right)+\\
&& t^9
\left(24310 x+2881330 x^2+40454572 x^3+140542120 x^4+149803150 x^5 + 
\right. \\
&&\left. 49462810 x^6+4200670 x^7+48610 x^8+x^9\right)+ \cdots.  
\end{eqnarray*}

We end this section with a remark about the case where 
$u=u_1 \ldots u_j$, $\des(u) =1$, $u_1 < u_j$, and 
$u$ does not have the weakly increasing overlapping property. There 
are two problems in this case. First, as we saw earlier, it 
is possible that the end points of collapse $u$-match in a fixed $(B,w)$
point of $I_u$ can lead to a rise, a level, or a descent in 
$C(B,w,u)$. This means that the weights $w_{u,\PP}(ij)$ or 
$w_{u,[k]}(ij)$ are much more complicated. The second 
problem is to find $U^{(\PP)}_u(x,\mathbf{z}_\infty,t)$,
we would need to substitute in a generating 
function of the form 
\begin{equation}\label{full}
1 + \sum_{n \geq 1} t^n \sum_{w =w_1 \ldots w_n \in \PP^n} 
\prod_{i=1}^{n-1} x_{w_i w_{i+1}}
\end{equation}
and we do not know of any way to find a compact form for such a generating 
function.

\section{The case $u = u_1 \ldots u_j$, $\des(u)=1$, and $u_1= u_j$}

In this section, we shall consider the problem 
of computing the generating functions \\
$\mathcal{N}^{(\PP)}_u(x,{\bf z}_\infty,t)$, 
$\mathcal{N}^{(k)}_u(x,{\bf z}_k,t)$, $\mathcal{EN}^{(\PP)}_u (x,{\bf z}_\infty,t)$, and 
$\mathcal{EN}^{(k)}_u(x,{\bf z}_k,t)$ for  
$u = u_1 \ldots u_j$ such that  $\des(u)=1$ and $u_1 = u_j$.

As in the previous sections, we need to compute 
$U^{(\PP)}_u(x,{\bf z}_\infty,t)$, 
$U^{(k)}_u(x,{\bf z}_k,t)$, $EU^{(\PP)}_u (x,{\bf z}_\infty,t)$, and 
$EU^{(k)}_u(x,{\bf z}_k,t)$. To compute these generating 
functions, we use Theorem \ref{thm:labwdes} or  \ref{thm:lablev} plus the 
collapse map.

\begin{figure}[htbp]
  \begin{center}
    \includegraphics[width=0.6\textwidth]{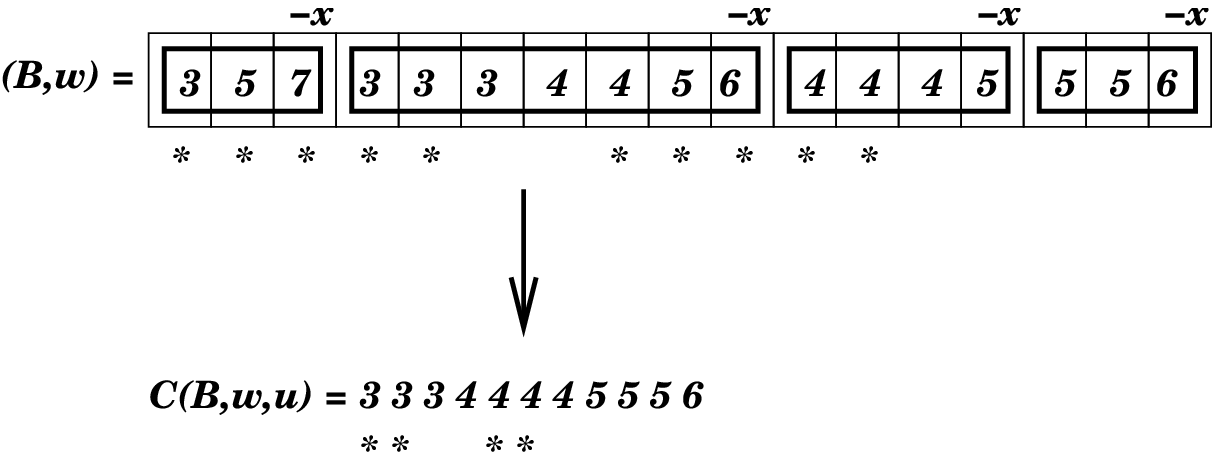}
    \caption{A fixed point of $I_{12311}$.}
    \label{fig:Collapse12311}
  \end{center}
\end{figure}

First assume that  $u = u_1 \ldots u_j$, $\red(u) =u$, $\des(u)=1$,  
$u_1  = u_j$, and $u$ has the $\PP$-
minimal overlapping property.  We can define the 
collapse map to fixed points of $I_u$ or $J_u$ exactly 
as in the previous sections. For example, suppose 
that $u = 12311$ and we want to compute $U^{(7)}_{12311}(x,{\bf z}_7,t)$.
By (\ref{eq:udes1}), we know 
that 
\begin{equation}\label{2eq:kudes1}
U^{(7)}_{12311,n}(x,{\bf z}_7) = \sum_{O \in \mathcal{O}^{(k)}_{12311,n},I_{12311}(O) = O}sgn(O)wt(O).
\end{equation}
As before, we know that if $(B,w)$ is a fixed point 
of $I_{12311}$, then elements in the bricks are weakly increasing 
and if there is a decrease between two brick $b_i$ and $b_{i+1}$, 
there must be a $12311$-match that involves the last 3 cells 
of $b_i$ and the first three cells of $b_{i+1}$.  We have 
pictured such a fixed point in Figure $\ref{fig:Collapse12311}$.

The difference between this case and the previous case 
where $u_1 >u_j$ is that a $12311$-match of 
the form $ijkii$ will just be replaced by $ii$ so that 
only factors of the form $ii$ could have come from 
a $12311$-match in the collapse of a fixed point of $I_{12311}$. The fact 
that $12311$ has the $\PP$-minimal overlapping property ensures 
that any two such $12311$-matches can only intersect 
at the right-hand endpoint of the first match and left-hand endpoint 
of the second match. It follows that $C(B,w,u)$ will always 
be a weakly increasing word. We claim that in this 
case a factor of the form $ii$ must 
have weight $1-x -xt^3 z_i \sum_{i < c <d \leq k} z_c z_d$ 
if we are computing $U^{(k)}_{12311,n}(x,{\bf z}_k)$  and  
 $1-x -xt^3 z_i \sum_{i < c <d } z_c z_d$ 
if we are computing $U^{(\PP)}_{12311,n}(x,{\bf z}_\infty)$. 
That is, the 1 corresponds to the case where $ii$ are in the same 
brick, the $-x$ corresponds to the case where the first $i$ is in 
last cell of some brick $b_j$ and the second $i$ is in the first 
cell of the next brick, and the third term corresponds to the 
cases where we have a decrease between two consecutive bricks 
and we deleted the second, third, and fourth elements 
of the $12311$-match between the two bricks. 
In our example, the weight of the levels 
for computing $U^{(7)}_{12311,n}(x,{\bf z}_7)$ would be as follows.

\begin{center}
\begin{tabular}{|l|l|}
\hline
Levels & $wt_{12311,7}(ii)$\\
\hline
$77$ & $1-x$ \\
\hline
$66$ & $1-x$ \\
\hline 
$55$ & 
$1-x -xt^3z_5z_6z_7$ \\
\hline
44 & $1-x-xt^3z_4 (\sum_{4 < c < d \leq 7} z_cz_d)$\\
\hline
33 & $1-x-xt^3z_3 (\sum_{3 < c < d \leq 7} z_cz_d)$\\
\hline
22 & $1-x-xt^3z_2 (\sum_{2 < c < d \leq 7} z_cz_d)$\\
\hline
11 & $1-x-xt^3z_1 (\sum_{1 <  c < d \leq 7} z_cz_d)$\\
\hline
\hline
\end{tabular}\\
The weights $wt_{12311,7}(ii)$.
\end{center}

In this case, rises in $C(B,w,12311)$ of the form 
$ij$ where $i <j$ correspond to a factor of $1-x$ where 
the 1 comes from the case where $ij$ are in the same brick 
and the $-x$ corresponds to the case where $i$ and $j$ are 
in different bricks. 

It follows 
that for any $v \in [7]^+$ which is weakly increasing, 
\begin{equation}\label{eq:v712311}
\sum_{\overset{(B,w) \ \mathrm{is \ a \ fixed \ point \ of} \ I_{12311}}{C(B,w,12311) =v}} sgn(B,w) wt_{12311}(B,w) =
-x\overline{z}^v t^{|v|}(1-x)^{\ris(v)} \prod_{s \in Lev(v)} wt_{12311,7}(v_s v_{s+1}).
\end{equation}
and
\begin{eqnarray}\label{eq:U712311}
U^{(7)}_{12311}(x,{\bf z}_7,t) &=& 1 + \sum_{n \geq 1} U^{(7)}_{12311,n}
(x,{\bf z}_7) t^n \nonumber \\
&=& 1 + \sum_{v \in [7]^+,\des(v) = 0} -x\overline{z}^v t^{|v|}(1-x)^{\ris(v)} \prod_{s \in Lev(v)} 
wt_{12311,7}(v_s v_{s+1}).
\end{eqnarray}

\begin{figure}[htbp]
  \begin{center}
    \includegraphics[width=0.6\textwidth]{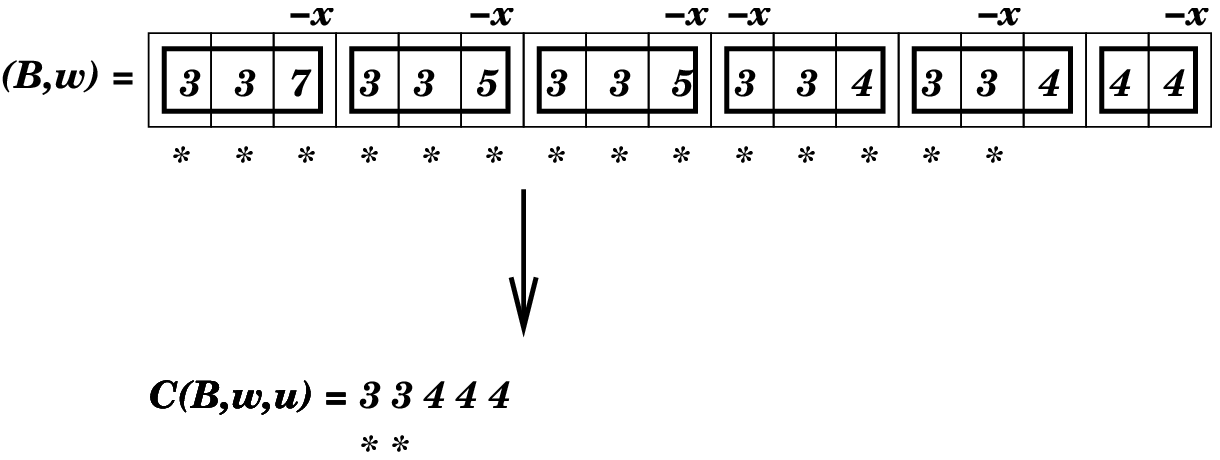}
    \caption{A fixed point of $I_{11211}$.}
    \label{fig:Collapse11211}
  \end{center}
\end{figure}

Next suppose that $u=u_1 \ldots u_j$, $\red(u) =u$, $\des(u)=1$,  $u_1=u_j$, 
and $u$ has the $\PP$-level overlapping property or 
the $[k]$-level overlapping property,  but $u$  
does not have the $\PP$-minimal overlapping property. The fact that 
$u$ has the $\PP$-level overlapping property ($[k]$-level overlapping property) ensures that if $w =w_1 \ldots w_n$ is word which starts and ends with 
a $u$-match and any two consecutive $u$-matches in $w$ share at least two 
letters, then it must be the case that $w_1 =w_n$.  Thus under the collapse 
map, any collapse will end up with a level of the form $ii$. 
The main difference in this case is that it is possible to have 
the weights $wt_{u,k}(ii)$ or $wt_{u,\PP}(ii)$  correspond to infinite 
families of words of different lengths even in the case 
where the alphabet is finite.  For example, suppose that $u =11211$. Then 
it is possible that in a fixed point 
$(B,w)$ of $I_{11211}$, $w$ has a factor where consecutive occurrences 
of the pattern 11211 are linked of the form 
$iiy_1iiy_2iiy_3ii \ldots iiy_nii$ where $y_1, \ldots, y_n >i$ like those that occur 
in the first 14 cells of the fixed point pictured in Figure \ref{fig:Collapse11211}. 
For each given maximal sequence of this type, the collapse map 
would eliminate all the symbols between the first and the last $i$.  In such 
a case, the weight corresponding to the symbols that are eliminated for such a string 
in the collapse map would be $(-x)^nz_i^{2n}z_{y_1}\cdots z_{y_n}t^{3n}$. 
It would follow that if we are working in $\PP^*$, then 
$$wt_{11211,\PP}(ii) = 1 - x + \frac{-xz_i^2\left(\sum_{s > i} z_s\right) t^3}{1+xz_i^2\left(\sum_{s > i} z_s\right) t^3}$$
while if we are working in $[k]^*$, then for $1 \leq i <k$, 
$$wt_{11211,k}(ii) = 1 - x + \frac{-xz_i^2\left(\sum_{s =i+1}^k z_s\right) t^3}{1+xz_i^2\left(\sum_{s =i+1}^k z_s\right) t^3}$$
and 
$$wt_{11211,k}(kk) = 1-x.$$
That is, in each of these expressions the 1 corresponds to the case where both 
$i$s are part of the same brick, the $-x$ corresponds to the case where the two $i$s are the 
last and first elements of two consecutive bricks, and the series 
$\frac{-xz_i^2\left(\sum_{s > i} z_s\right) t^3}{1+xz_i^2\left(\sum_{s > i} z_s\right) t^3}$ 
corresponds the fact that we could have eliminated sequences of the form 
$iy_1iiy_2iiy_3ii \ldots iiy_ni$ for any $n \geq 1$ between the two $i$s.

Nevertheless,  we can still apply the same reasoning as above to prove  
that for any $v \in [7]^+$ which is weakly increasing, 
\begin{equation}\label{eq:v711211}
\sum_{\overset{(B,w) \ \mathrm{a \ fixed \ point \ of} \ I_{11211}}{C(B,w,11211) =v}} sgn(B,w) wt_{11211}(B,w) =
-x\overline{z}^v t^{|v|}(1-x)^{\ris(v)} \prod_{s \in Lev(v)} wt_{11211,7}(v_s v_{s+1}).
\end{equation}
and
\begin{eqnarray}\label{eq:U711211}
U^{(7)}_{11211}(x,{\bf z}_7,t) &=& 1 + \sum_{n \geq 1} U^{(7)}_{11211,n}(x,{\bf z}_7) t^n \nonumber \\
&=& 1 + \sum_{v \in [7]^+,\des(v) =0} -x\overline{z}^v t^{|v|}(1-x)^{\ris(v)} \prod_{s \in Lev(v)} 
wt_{11211,7}(v_s v_{s+1}).
\end{eqnarray}

We should note that as patterns get more complicated, it becomes 
increasingly difficult to compute $wt_{u,\PP}(ii)$ or $wt_{k}(ii)$. 
For example, 
suppose $u =3^545123^5$. Then linked patterns can overlap at either 1,2,3,4, or 5 symbols.

It follows from Theorem \ref{thm:lablev} that 
\begin{equation}\label{inc2}
\sum_{v \in \PP^+,\des(v)=0}
\overline{z}^v t^{|v|} \prod_{i \in Lev(v)} x_{v_iv_i} = -1+ \prod_{i \geq 1} \left(1+\frac{z_it}{1-x_{ii}z_it}\right).
\end{equation}
Replacing $t$ by $yt$ and $x_{ij}$ by $x_{ii}/y$, we see that 
\begin{equation}\label{inc3}
\sum_{\overset{v =v_1 \ldots v_n \in \PP^+}{v_1 \leq v_2 \leq \cdots \leq v_n}}
\overline{z}^v t^{|v|} y^{\ris(v)+1} \prod_{i \in Lev(v)} x_{v_iv_i} = 
-1+ \prod_{i \geq 1} \left(1+\frac{yz_it}{1-x_{ii}z_it}\right).
\end{equation}
Thus 
\begin{equation}\label{eq:PPinc2}
1+ \sum_{\overset{v =v_1 \ldots v_n \in \PP^+}{v_1 \leq v_2 \leq \cdots \leq v_n}}
-x\overline{z}^v t^{|v|} y^{\ris(v)} \prod_{i \in Lev(v)} x_{v_iv_i} = 
1+\frac{-x}{y}\left(-1+\prod_{i \geq 1} \left(1+\frac{yz_it}{1-x_{ii}z_it}\right)\right).
\end{equation}
and 
\begin{equation}\label{eq:kinc2}
1+ \sum_{\overset{v =v_1 \ldots v_n \in [k]^+}{v_1 \leq v_2 \leq \cdots \leq v_n}}
-x\overline{z}^v t^{|v|} y^{\ris(v)} \prod_{i \in Lev(v)} x_{v_iv_i} = 
1+\frac{-x}{y}\left(-1+\prod_{i =1}^k \left(1+\frac{yz_it}{1-x_{ii}z_it}\right)\right).
\end{equation}
But then it follows that if $u =u_1 \ldots u_j$, $\red(u) =u$, $\des(u)=1$, $u_1=u_j$, and $u$ has the $\PP$-level overlapping property, then 
\begin{eqnarray*}
U^{(\PP)}_{u}(x,{\bf z}_\infty,t) &=& 
1+ \sum_{v \in \PP^+,\des(v)=0}
-x\overline{z}^v t^{|v|} (1-x)^{\ris(v)} \prod_{i \in Lev(v)} wt_{u,\PP}(v_i v_i)\\
&=& 
1+\frac{-x}{1-x}\left(-1+\prod_{i \geq 1} \left(1+\frac{((1-x)z_it}{1-wt_{u,\PP}(ii)z_it}\right)\right)\end{eqnarray*}
and, for all $k \geq 1$, if $u =u_1 \ldots u_j$, $\red(u) =u$, $\des(u)=1$, $u_1=u_j$, and $u$ has the $[k]$-level overlapping property, then 
\begin{eqnarray*}
U^{(k)}_{u}(x,{\bf z}_k,t) &=& 
1+ \sum_{v \in [k]^+,\des(v)=0}
-x\overline{z}^v t^{|v|} (1-x)^{\ris(v)} \prod_{i \in Lev(v)} wt_{u,k}(v_i v_i)\\
&=& 
1+\frac{-x}{1-x}\left(-1+\prod_{i=1}^k \left(1+\frac{((1-x)z_it}{1-wt_{u,k}(ii)z_it}\right)\right).
\end{eqnarray*}

Thus we have the following theorem.

\begin{theorem}\label{thm:Ninc} If $u =u_1 \ldots u_j \in \PP^*$ is such 
that $\red(u) =u$, $\des(u) =1$, $u_1 = u_j$, and $u$ has 
the $\PP$-level overlapping property, then  
\begin{equation}\label{eq:NincPP}
\mathcal{N}^{(\PP)}_{u}(x,{\bf z}_\infty,t) = 
\frac{1}{1-\frac{x}{1-x}\left(-1+\prod_{i \geq 1} 
\left(1+\frac{(1-x)z_it}{1-wt_{u,\PP}(ii)z_it}\right)\right)}.
\end{equation}
If $u =u_1 \ldots u_j \in [k]^*$ is such 
that $\red(u) =u$, $\des(u) =1$, $u_1 = u_j$, and $u$ has 
the $[k]$-level overlapping property, then
\begin{equation}\label{eq:Ninck}
\mathcal{N}^{(k)}_{u}(x,{\bf z}_k,t) = \frac{1}{1-\frac{x}{1-x}\left(-1+\prod_{i=1}^k \left(1+\frac{(1-x)z_it}{1-wt_{u,k}(ii)z_it}\right)\right)}
\end{equation}
\end{theorem}

Note that if $u =u_1 \ldots u_j$, $\des(u) = 1$, $u_1 =u_j$,
then  $u$ automatically has the exact 
$\PP$-level overlapping property (exact $[k]$-level overlapping 
property). 

\begin{theorem}\label{thm:elev} If $u =u_1 \ldots u_j \in \PP^*$ is such 
that $\des(u) =1$ and $u_1 = u_j$, then  
\begin{equation}\label{eq:eNincPP}
\mathcal{EN}^{(\PP)}_{u}(x,{\bf z}_\infty,t) = 
\frac{1}{1-\frac{x}{1-x}\left(-1+\prod_{i \geq 1} 
\left(1+\frac{(1-x)z_it}{1-ewt_{u,\PP}(ii)z_it}\right)\right)}.
\end{equation}
and if $u =u_1 \ldots u_j \in [k]^*$ is such that $\des(u) =1$ and $u_1 = u_j$, then
\begin{equation}\label{eq:eNinck}
\mathcal{EN}^{(k)}_{u}(x,{\bf z}_k,t) = \frac{1}{1-\frac{x}{1-x}\left(-1+\prod_{i=1}^k \left(1+\frac{(1-x)z_it}{1-ewt_{u,k}(ii)z_it}\right)\right)}
\end{equation}
\end{theorem}

For example, suppose we want to compute 
$\mathcal{N}^{(7)}_{12311}(x,{\bf z}_7,t)$ where we set $z_i=1$ for all $i$. 
It follows from (\ref{eq:Ninck}) that  
\begin{equation}\label{series12311}
\mathcal{N}^{(7)}_{12311}(x,{\bf z}_7,t) = 
\frac{1}{1-\frac{x}{(1-x)}(-1+\prod_{i=1}^7 Q_i(x,t))}
\end{equation}
where 
\begin{enumerate}
\item $Q_1(x,t) = 1+\frac{(1-x)t}{1-(1-x-15xt^3)t}$, 
\item $Q_2(x,t) = 1+\frac{(1-x)t}{1-(1-x-10xt^3)t}$, 
\item $Q_3(x,t) = 1+\frac{(1-x)t}{1-(1-x-6xt^3)t}$, 
\item $Q_4(x,t) = 1+\frac{(1-x)t}{1-(1-x-3xt^3)t}$, 
\item $Q_5(x,t) = 1+\frac{(1-x)t}{1-(1-x-xt^3)t}$, 
\item $Q_6(x,t) = 1+\frac{(1-x)t}{1-(1-x)t}$, and 
\item $Q_7(x,t) = 1+\frac{(1-x)t}{1-(1-x)t}$. 
\end{enumerate}
We have computed that 
\begin{align*}
&\mathcal{N}^{(7)}_{12311}(x,\bf{z}_7,t) = \\
& 1 +7xt +7(4x+3x^2)t^2+7(12x+32x^2+5x^3)t^3 + \\
&7(30x+190x^2+118x^3+5x^3)t^4 +7(66x+823x^2+1236x^3+268x^3+3x^5)t^5 + \\
&7(132x+2912x^2+8500x^3+4770x^4+422x^5 +x^6)t^6 + \\
&(1716x +62532x^2+312558x^3 +349315x^4+88852x^4+3424x^6+x^7)t^7+\\
&7(429x+24609x^2+194029x^3+374249x^4+197729x^525209x^6+429x^7)t^8 +\cdots.
\end{align*}

Finally we shall consider the case where 
$u =u_1 \ldots u_j$, $\red(u) =u$, $\des(u) = 1$, $u_1 =u_j$ and $u$ does not 
have the $\PP$-level overlapping property ($[k]$-level overlapping 
property).  Given such a $u$, let $s$ be the position such 
that $u_s > u_{s+1}$.  Then we must have that 
$u_{s+1} \leq \cdots \leq u_j =u_1$ and $St^{(\PP)}(u) \subset 
\{s+1, \ldots, j\}$ ($St^{([k])}(u) \subset \{s+1, \ldots, j\}$). 
This means that $u$ automatically has the $\PP$-weakly decreasing 
overlapping property ($[k]$-weakly decreasing overlapping property) and 
$u$ is not $\PP$-minimal overlapping ($[k]$-minimal overlapping). 
Now suppose that 
$w = w_1 \ldots w_n$ is a maximal sequence of linked $u$-matches. 
That is, we assume $w$ starts and ends with a $u$-match 
and any two consecutive $u$-matches share at least two letters. 
Then if the $u$-matches in $w$ start at positions 
$1=i_1 < i_2 < \cdots < i_k$, then the $\PP$-weakly decreasing 
overlapping property ensures that 
$w_1=w_{i_1} \geq \cdots \geq w_{i_k} =w_n$. Thus in a collapse map, 
if  we eliminate $w_2 \ldots w_{n-1}$, then 
 we will be left with a weak descent 
$w_1 w_n$. Thus we must figure out the weights $wt_{u}(ji)$ for $j \geq i$. 

To illustrate the process, we will consider the example 
where $u =2312$ and the alphabet is $[4]$. If $w =w_1w_2w_3w_4 \in [4]^*$ and 
$\red(w) =2312$, then clearly $w$ must start with  either 2  or 3 since 
those are the only letters $a$ which have at least one letter in 
$[4]$ bigger than $a$ and one letter in $[4]$ which is less than $a$. 
It follows that $wt_{2312,4}(44)= wt_{2312,4}(11)=1-x$. Also 
$wt_{2312,4}(4i) =0$ for $i =1,2,3$ and $wt_{2312,4}(j1) = 0$ for 
$j=2,3,4$.

Next consider $wt_{2312,4}(22)$. There are only two 
possible words in $[4]^4$ that reduce to $u$, namely, 
$w=2312$ and $v=2412$.  Since there is no $u$-match that can 
start with 1, there cannot be a pair of linked $u$-matches that 
start with either $w$ or $v$. Thus there can be no maximal sequences 
of linked $u$-matches that start and end with 2. This means 
that when we collapsed to 22,  either we started with 
$2312$ and eliminated 31 or we started with $2412$ and we eliminated 
41. It follows that $wt_{2312,4}(22)=1-x-xz_1(z_3+z_4)t^2$.

Next consider $wt_{2312,4}(33)$. There are only two 
possible words in $[4]^4$ that reduce to $u$, namely, 
$w=3413$ and $v=3423$.  Since there is no $u$-match that can 
start with 1, there cannot  a be pair of linked $u$-matches that 
start with $w$. There is a pair of linked $u$ matches that 
start with $v$, namely, $342312$. However this pair can not  be 
extended. Thus there can be no maximal sequences 
of linked $u$-matches that start and end with 3. 
This means 
that when we collapsed to 33 either we started with 
$3413$ and eliminated 41 or we started with $3423$ and eliminated 
42. It follows that $wt_{2312,4}(33)=1-x-xz_4(z_1+z_2)t^2$.

Finally we consider $wt_{2312,4}(32)$. In this case, the only possible 
way to have a maximal sequence $w$ 
of linked $u$-matches starting with a $u$-match 
whose first letter is 3 and ending with a $u$-match whose last letter is 2 
is  $w = 342312$.  
Since in the fixed points of $I_{2312}$, the sequences in the bricks are 
weakly increasing, the only way that 32 occurs in the collapse of fixed 
point $(B,w)$ of $I_{2312}$ is if we started with $342312$, which 
means that a brick ended after 4 and a brick ended after the second 
3, and 
eliminated $4231$. Hence $wt_{2312,4}(32) = 
x^2z_1z_2z_3z_4t^4$. 

Thus we have the following table for $wt_{2312,4}(32)$

\begin{center}
\begin{tabular}{|l|l|}
\hline
Weak Descents & $wt_{2312,4}(ji)$\\
\hline
$44$ & $1-x$ \\
\hline
$4i\  (i < 4)$ & $0$\\
\hline 
$33$ & $1-x-xz_4(z_1+z_2)t^2$ \\
\hline
$32$ & $x^2 z_1z_2z_3z_4t^4$ \\
\hline
$31$ & 0 \\
\hline
$22$ &  $1-x-xz_1(z_3+z_4)t^2$\\
\hline
$21$ & 0 \\
\hline
$11$ & $1-x$ \\
\hline
\end{tabular}\\
The weights $wt_{2312,4}(ji)$.
\end{center}

It follows 
that for any $v \in [4]^+$, 
\begin{equation}\label{eq:v42312}
\sum_{\overset{(B,w) \ \mathrm{is\ a\ fixed\ point\ of} \ I_{2312}}
{C(B,w,2312) =v}}
sgn(B,w) wt(B,w) =
-x \overline{z}^v t^{|v|} (1-x)^{\ris(v)} \prod_{s \in WDes(v)} 
wt_{2312,4}(v_s v_{s+1}).
\end{equation}
Here the initial $-x$ comes from the fact that the last cell of $(B,w)$ always 
contributes a $-x$ since the last cell is at the end of a brick. 
It follows that 
\begin{eqnarray}\label{eq:U42312}
U^{(4)}_{2312}(x,{\bf z}_4,t) &=& 1 + \sum_{n \geq 1} U^{(4)}_{2312,n}(x,{\bf z}_4) t^n \nonumber \\
&=& 1 + \sum_{v \in [4]^+} -x(1-x)^{\ris(v)} \overline{z}^v t^{|v|} \prod_{s \in WDes(v)} wt_{2312,4}(v_s v_{s+1}).
\end{eqnarray}
Hence we could compute 
$\displaystyle \mathcal{N}^{(4)}_{2312}(x,{\bf z}_4,t) = \frac{1}{U^{(4)}_{2312}(x,{\bf z}_4,t)}$
if we can compute the right-hand side of (\ref{eq:U42312})

What we need to be able to compute the right-hand side of  
(\ref{eq:U42312}) is the generating function over all words $v \in \PP^*$ where we not only keep track of the weak descents of $P$ but also 
of type of weak descents of $P$. 

By Theorem \ref{thm:labwdes}, 
we know that 
\begin{equation}\label{PPWDES1}
\frac{1}{1 - \sum_{n \geq 1} t^n \sum_{v \in  WD\PP^*, |v| =n} 
WDXZ(v)} = 1+\sum_{w = w_1 \ldots w_n \in \PP^+} t^{|w|}\overline{z}^w \prod_{i \in WDes(w)} x_{w_i w_{i+1}}.\end{equation}
Hence 
\begin{align}\label{PPWDES2}
\sum_{w = w_1 \ldots w_n \in \PP^+} t^{|w|}\overline{z}^w \prod_{i \in WDes(w)}
x_{w_i w_{i+1}} &= 
\left(\frac{1}{1 - \sum_{n \geq 1} t^n \sum_{v \in WD\PP^*, |v| =n} 
WDXZ(v)}\right) -1 \nonumber \\
&=\frac{\sum_{n \geq 1} t^n \sum_{v \in WD\PP^*, |v| =n} 
WDXZ(v)}{1 - \sum_{n \geq 1} t^n \sum_{v \in WD\PP^*, |v| =n} 
WDXZ(v)}.
\end{align}
Next suppose that  we replace $t$ by $yt$ and $x_{ij}$ by 
$\frac{x_{ij}}{y}$.  Under this substitution 
the left-hand side in (\ref{PPWDES2}) becomes 
 
$$\sum_{w = w_1 \ldots w_n \in \PP^+} t^{|w|}y^{\ris(w)+1}\overline{z}^w \prod_{i \in WDes(w)}x_{w_i w_{i+1}}.$$
Note that for $v =j_1 \geq  \cdots \geq  j_k$ where $k \geq 2$, our 
substitution replaces 
$t^k WDXZ(v)$ by 
$$y^kt^kz_{j_1} \cdots z_{j_k} \prod_{i=1}^{k-1} 
\left(\frac{x_{j_{i+1}j_i}}{y} -1\right) = 
yt^k z_{j_1} \cdots z_{j_k} \prod_{i=1}^{k-1} 
(x_{j_{i+1}j_i} -y).$$
Thus if we let 
\begin{equation}
WDXYZ(v) = \begin{cases} 
z_j & \mbox{if} \ v = j, \ \mbox{and} \\
z_{j_1} \cdots z_{j_k} \prod_{i=1}^{k-1} (x_{j_{i+1}j_i}-y) & 
\mbox{if} \ v = j_1 \geq  \cdots \geq  j_k \ \mbox{where} \ k \geq 2,
\end{cases}
\end{equation}
then we see that the right-hand side of (\ref{PPWDES2}) becomes
$$
\frac{y \sum_{n \geq 1} t^n \sum_{v \in WD\PP^*, |v| =n} 
WDXYZ(v)}{1 - y \sum_{n \geq 1} t^n \sum_{v \in  WD\PP^*, |v| =n} 
WDXYZ(v)}.
$$
It follows that
\begin{equation*}
-x\sum_{w = w_1 \ldots w_n \in \PP^+} t^{|w|}y^{\ris(w)}\overline{z}^w \prod_{i \in WDes(w)}x_{w_i w_{i+1}} =
\frac{-x\sum_{n \geq 1} t^n \sum_{v \in WD\PP^*, |v| =n} 
WDXYZ(v)}{1 - y \sum_{n \geq 1} t^n \sum_{v \in WD\PP^*, |v| =n} 
WDXYZ(v)}.
\end{equation*}
Thus 
\begin{multline}\label{PPWDES3}
1-x\sum_{w = w_1 \ldots w_n \in \PP^+} t^{|w|}y^{rise(w)}\overline{z}^w \prod_{i \in WDes(w)}x_{w_i w_{i+1}} =\\
\frac{1-(x+y)\sum_{n \geq 1} t^n \sum_{v \in  WD\PP^*, |v| =n} 
WDXYZ(v)}{1 - y \sum_{n \geq 1} t^n \sum_{v \in WD \PP, |v| =n} 
WDXYZ(v)}.
\end{multline}
By setting $z_i=0$ for $i > k$, we also 
obtain that 
\begin{multline}\label{kWDES3}
1-x\sum_{w = w_1 \ldots w_n \in [k]^+} t^{|w|}y^{\ris(w)}\overline{z}^w \prod_{i \in WDes(w)}x_{w_i w_{i+1}} = \\
\frac{1-(x+y)\sum_{n = 1}^{k} t^n \sum_{v \in WD[k]^*, |v| =n} 
WDXYZ(v)}{1 - y \sum_{n = 1}^{k} t^n \sum_{v \in WD[k]^*, |v| =n} 
WDXYZ(v)}.
\end{multline}

Note that if we replace $y$ by $(1-x)$ 
and $x_{ji}$ by $wt_u(ji)$, the left-hand side 
of (\ref{PPWDES3}) becomes 
$U^{(\PP)}_u(x,{\bf z}_\infty,t)$ and the left-hand side 
of (\ref{kWDES3}) becomes 
$U^{(k)}_u(x,{\bf z}_k,t)$. 
Then using the fact that 
$\mathcal{N}^ {(\PP)}_u(x,{\bf z}_\infty,t) = 
1/U^{(\PP)}_u(x,{\bf z}_\infty,t)$,  
we have the following 
theorem.

\begin{theorem}\label{thm:u1=ujWDPP} 
Suppose that $u = u_1 \ldots u_j \in \PP^*$, $\red(u) =u$, 
$\des(u) =1$, $u_1 =u_j$, and $u$ does not have the $\PP$-level overlapping property (so it automatically has the $\PP$-weakly decreasing property). Then 
\begin{equation} 
\mathcal{N}^ {(\PP)}_u(x,{\bf z}_\infty,t) = 
\frac{1 - (1-x) \sum_{n \geq 1} t^n \sum_{v \in  WD\PP^*, |v| =n} 
WDXTZ_u(v)}{1-\sum_{n \geq 1} t^n \sum_{v \in WD\PP^*, |v| =n} 
WDXTZ_u(v)}
\end{equation}
where 
\begin{equation}
WDXTZ_u(v) = \begin{cases} 
z_j & \mbox{if} \ v = j, \ \mbox{and} \\
z_{j_1} \cdots z_{j_k} \prod_{i=1}^{k-1} (wt_u(j_{i+1}j_i)+x-1) & 
\mbox{if} \ v= j_1 \geq \cdots \geq j_k \ \mbox{where} \ k \geq 2.
\end{cases}
\end{equation}
\end{theorem}

If set $z_i =0$ for 
all $i > k$, then we obtain  the following theorem.

\begin{theorem}\label{thm:u1=ujWDk} 
Now suppose that $u = u_1 \ldots u_j \in [k]^*$, $\red(u) =u$, 
$\des(u) =1$, $u_1= u_j$, and $u$ does not have the $[k]$-level overlapping property (so it automatically has the $[k]$-weakly decreasing property). Then 
\begin{equation} \label{u1equjnl}
\mathcal{N}^ {(k)}_u(x,{\bf z}_k,t) = 
\frac{1 - (1-x) \sum_{n = 1}^{k} t^n \sum_{v \in WD[k]^*, |v| =n} 
WDXTZ_u(v)}{1-\sum_{n = 1}^{k} t^n \sum_{v \in WD[k]^*, |v| =n} 
WDXTZ_u(v)}.
\end{equation}
\end{theorem}

The key to be able to compute 
$\mathcal{N}^ {(k)}_u(x,{\bf z}_k,t)$ in the case of 
Theorem \ref{thm:u1=ujWDk} is to be able to compute 
$\sum_{n \geq 1} t^n \sum_{v \in WD[k]^*, |v| =n} WDXTZ_u(v)$.
This is often complicated because of the large number 
of weakly decreasing words in $WD[k]^*$, but for certain patterns 
we can compute it. For example, consider the 
case where $u =2312$, $k=4$, and we set $z_i=1$ for $i =1, \ldots ,4$. 
With this substitution, the table of $WDXTX_{2312}(ij)$ becomes 

\begin{center}
\begin{tabular}{|l|l|}
\hline
Weak Descents & $WDXTX_{2312}(ij)$\\
\hline
$44$ & $0$ \\
\hline
$4i\  (i < 4)$ & $x-1$\\
\hline 
$33$ & $-2xt^2$ \\
\hline
$32$ & $x-1 +x^2t^4$ \\
\hline
$31$ & $x-1$ \\
\hline
$22$ &  $-2xt^2$\\
\hline
$21$ & $x-1$ \\
\hline
$11$ & $0$ \\
\hline
\end{tabular}\\
The weights $WDXTZ_{2312}(ij)$ in the case $z_i=0$ for $i =1, \ldots, 4$.
\end{center}

Because $WDXTZ_{2312}(44)=WDXTZ_{2312}(11)=0$, it follows 
that the only words that we have to consider 
are $1$, $4$, $41$ and words in 
$ (\epsilon +4) (\{2\}^+ +\{3\}^+ +\{3\}^*32\{2\}^*)(\epsilon+1)$. 
It is easy to see that 
$$\sum_{n \geq 1} t^n \sum_{v \in \{3\}^+, |v| =n} WDXTZ_u(v) = 
\frac{t}{1+2xt^3}.$$ That is, the first 3 gives a factor of $t$ and each 
additional 3 gives a factor of $-2xt^3$.  
Similarly, 
$$\sum_{n \geq 1} t^n \sum_{v \in \{2\}^+, |v| =n} WDXTZ_u(v) = 
\frac{t}{1+2xt^3}.$$
When considering words in $\{3\}^* 32 \{2\}^*$, the 
$32$ gives a factor of $(x-1)t^2+x^2t^6$ and each 
additional 3 to the left gives a factor of $-2xt^3$ and each additional 
2 to the right gives a factor of $-2xt^3$. Thus 
$$\sum_{n \geq 1} t^n \sum_{v \in \{3\}^*32\{2\}^*, |v| =n} WDXTZ_u(v) = 
\frac{(x-1)t^2+x^2t^6}{(1+2xt^3)^2}.$$
Thus 
$$\sum_{n \geq 1} t^n \sum_{v \in \{3\}^+ + \{2\}^+ + \{3\}^*32\{2\}^*, |v| =n} WDXTZ_u(v) = \frac{2t+4xt^4+(x-1)t^2 +x^2t^6}{(1+2xt^3)^2}.$$
Hence if $E = (\epsilon +4) (\{3\}^+ + \{2\}^+ + \{3\}^*32\{2\}^*)(\epsilon+1)$, it follows 
that 
$$\sum_{n \geq 1} t^n \sum_{v \in E, |v| =n} WDXTZ_{2312}(v) = \frac{(1+(x-1)t)^2(2t+4xt^4+(x-1)t^2 +x^2t^6)}{(1+2xt^3)^2}$$
since adding a 4 to the left of a word $w \in \{3\}^+ + \{2\}^+ + \{3\}^*32\{2\}^*$ gives rise to a factor of $(x-1)t$ and 
adding a 1 to the right of a word $w \in \{3\}^+ + \{2\}^+ + \{3\}^*32\{2\}^*$ gives rise to a factor of $(x-1)t$.
It follows that  
\begin{eqnarray*}
\sum_{n \geq 1} t^n \sum_{v \in WD[4]^*, |v| =n} WDXTZ_{2312}(v)  &=& 
2t +(x-1)^t + \frac{(1+(x-1)t)^2(2t+4xt^4+(x-1)t^2 +x^2t^6)}{(1+2xt^3)^2} \\
&=& \frac{P(x,t)}{(1+2xt^3)^2}.
\end{eqnarray*}
where 

\begin{eqnarray*} P(x,t) &=& 4 t+(-6+6 x) t^2+\left(4-8 x+4 x^2\right) t^3+\left(-1+15 x-3 x^2+x^3\right) t^4+\\
&&\left(-12 x+12 x^2\right) t^5+\left(4 x-7 x^2+4 x^3\right)
t^6+\\
&&\left(6 x^2+2 x^3\right) t^7+\left(-3 x^2+2 x^3+x^4\right) t^8.
\end{eqnarray*}
Thus
\begin{equation}
\mathcal{N}^ {(4)}_{2312}(x,1,1,1,1,t) = 
\frac{1-(x-1) \frac{P(x,t)}{(1+2xt^3)^2}}{1-\frac{P(x,t)}{(1+2xt^3)^2}}.
\end{equation}
We have used Mathematica to compute the first few terms in 
this series:
\begin{eqnarray*}
&&1+4 x t+2 \left(5 x+3 x^2\right) t^2+4 \left(5 x+10 x^2+x^3\right) t^3+\left(35 x+151 x^2+65 x^3+x^4\right) t^4+\\
&&4 \left(14 x+109 x^2+111
x^3+14 x^4\right) t^5+\left(84 x+1068 x^2+2009 x^3+716 x^4+28 x^5\right) t^6+\\
&&2 \left(60 x+1166 x^2+3561 x^3+2535 x^4+362 x^5+4 x^6\right) t^7+\\
&& \left(165
x+4670 x^2+21400 x^3+25650 x^4+8172 x^5+486 x^6+x^7\right) t^8 + \cdots .
\end{eqnarray*}

\section{The proofs of Theorems \ref{thm:labdes}, \ref{thm:labwdes}, \ref{thm:labris}, and \ref{thm:lablev}}

In this section, we shall prove Theorems \ref{thm:labdes}, \ref{thm:labwdes}, \ref{thm:labris}, 
and \ref{thm:lablev}. 

We start with the proof of Theorem \ref{thm:labdes}. \\
\ \\
{\bf Proof of Theorem \ref{thm:labdes}.}\\

Recall that given a set  
$S \subseteq \PP$, we let 
\begin{equation}
DXZ(S) = \begin{cases} 
z_j & \mbox{if} \ S = \{j\}, \ \mbox{and} \\
z_{j_1} \cdots z_{j_k} \prod_{i=1}^{k-1} (x_{j_{i+1}j_i}-1) & 
\mbox{if} \ S =\{j_1 < \cdots < j_k\} \ \mbox{where} \ k \geq 2.
\end{cases}
\end{equation}

Define a ring homomorphism $\Gamma:\Lambda \rightarrow Q[\bf{x},\bf{z}]$ by 
defining $\Gamma(e_0) =1$ and, for $n \geq 1$,  
\begin{equation}\label{Gammadef}
\Gamma(e_n) = (-1)^{n-1} \sum_{S \subseteq \PP, |S|= n} DXZ(S).
\end{equation}

Then we claim that 
\begin{equation}\label{des1}
\Gamma(h_n) = \sum_{w \in \PP^n} \overline{z}^w \prod_{i <j} 
x_{ji}^{\mathbf{ji}(w)}.
\end{equation}
That is,
\begin{eqnarray}\label{des2}
\Gamma(h_{n}) &=& 
\sum_{\mu\vdash n} (-1)^{n-\ell(\mu)}B_{\mu, (n)}\Gamma(e_{\mu}) \nonumber \\
&=& \sum_{\mu \vdash n} (-1)^{n-\ell(\mu)} \sum_{(b_1,
\ldots, b_{\ell(\mu)}) \in \mathcal{B}_{\mu,(n)}}
\prod_{j=1}^{\ell(\mu)} (-1)^{b_j-1} 
\sum_{S_j \subseteq  \PP, |S_j| =b_j} DXZ(S_j) \nonumber \\
&=& \ \sum_{\mu \vdash n} \sum_{(b_1, \ldots, b_{\ell(\mu)}) \in \mathcal{B}_{\mu,(n)}} \prod_{j=1}^{\ell(\mu)}  \sum_{S_j \subseteq  \PP, |S_j| =b_j} DXZ(S_j)
\end{eqnarray}

Next we want to give a combinatorial interpretation to
(\ref{des2}).  First we pick a brick tabloid $B = (b_1, \ldots, b_k)$ 
of length $n$. 
Then we interpret 
$\prod_{j=1}^{\ell(\mu)}  \sum_{S_j \subseteq  \PP, |S_j| =b_j} DXZ(S_j)$ 
as picking a sequence of subsets of $\PP$, $(S_1, \ldots, S_{\ell(\mu)})$,  
such that $S_j$ has size $b_j$ and placing the elements of $S_j$ in the cells of $b_j$ in decreasing order for $j=1, \ldots, \ell(\mu)$. If 
$S_j = \{a_1 >  \cdots > a_{b_j}\}$, then we interpret the 
factor $DXZ(S_j) = z_{a_1} \cdots z_{a_{b_j}} \prod_{i=1}^{b_j-1} (x_{a_i a_{i+1}}-1)$ as the ways of labeling the cells of $b_j$ that contain  
$a_i$ where $i <b_j$ with either 
$z_{a_i}x_{a_i,a_{i+1}}$ or with $-z_{a_i}$ and labeling 
the last of cell $b_j$ with $z_{a_{b_j}}$.  
We shall call all such objects created in this way \emph{filled labeled
brick tabloids} and let $\mathcal{H}_{n}$ denote the set of all
filled labeled brick tabloids that arise in this way.  
Thus $\mathcal{H}_{n}$ consists of all triples $(B,w,L)$ such that 
\begin{enumerate}
\item $B =(b_1, \ldots, b_k)$ is a brick tabloid of length $n$, 
\item $w=w_1 \ldots w_n$ is a word in $\PP^n$ such that 
$w$ is strictly decreasing in each brick, and 
\item $L$ is a labeling of the cells of $B$ such that 
$L(i)$ is equal to $z_a$ if $i$ is the last cell of some brick 
$b_j$ which contains $a$ and $L(i) =-z_a$ or $L(i) = x_{ab}z_a$ if 
$i$ is not the last cell of a brick, cell $i$ contains $a$ and cell 
$i+1$ contains $b$. 
\end{enumerate} 
We then define the weight of $(B,w,L)$, $wt(B,w,L)$, to be the product of all
the $x_{ab}$ and $z_a$ labels in $L$ and the sign of 
$(B,w,L)$, $sgn(B,w,L)$, to be
the product of all the $-1$ factors in the labels in $L$. This 
process is illustrated in Figure \ref{fig:fdes1} to construct 
an element $(B,w,L)$ of $H_{12}$ such that 
$$wt(B,w,L) = z_1^3 z_2^2 z_3^2 z_4^2 z_5 z_6^2 x_{54} x_{41} x_{64} 
x_{32} x_{63}$$
and $sgn(B,w,L) = -1$. 

Thus
\begin{equation}\label{des4}
\Gamma(h_{n}) = \sum_{(B,w,L) \in \mathcal{H}_{n}}
sgn(B,w,L) wt(B,w,L).
\end{equation}

\fig{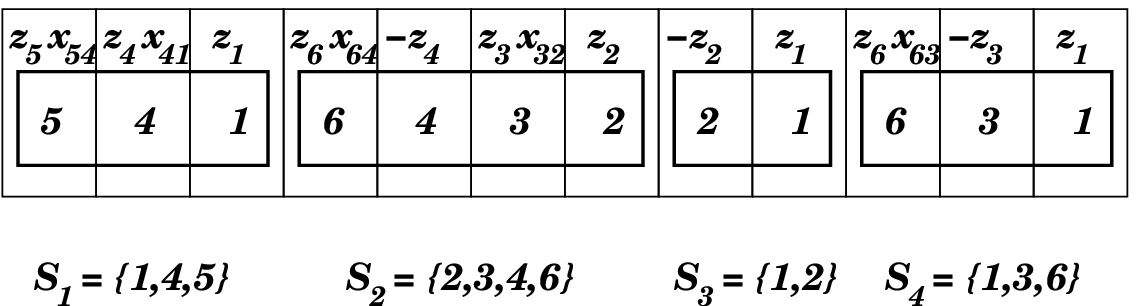}{A element $(B,w,L) \in \mathcal{H}_{12}$.}

Next we define a weight-preserving sign-reversing involution
$I:\mathcal{H}_{n} \rightarrow \mathcal{H}_{n}$.  
To define
$I(C)$, we scan the cells of $C =(B,w,L)$ from left  to right
looking for the leftmost cell, $t$, such that either (i) $t$ is
labeled with $-z_{w_t}$ or (ii) $t$ is at the end of a brick, $b_j$, 
there is a brick $b_{j+1}$ immediately following $b_j$, and 
$w_t > w_{t+1}$.  In case (i), $I(C)
=(B,w',L')$ where $B'$ is the result of  replacing the brick
$b$ in $B$ containing $t$ by two bricks $b^*$ and $b^{**}$, where $b^*$ 
contains all the cells of $b$ weakly to the left of cell $t$ and 
$b^{**}$ contains all the cells of $b$ strictly to the right of
$t$, $w' =w$, and $L'$ is the labeling that results
from $L$ by changing the label of cell  $t$ from $-z_{w_t}$ to $z_{w_t}$. In
case (ii), $I(C) =(B',w',L')$ where $B$ is the result of
replacing the bricks $b_j$ and $b_{j+1}$ in $B$ by a single brick
$b,$ $w' = w$, and $L'$ is the labeling that results
from $L$ by changing the label of cell $t$ from $z_{w_t}$ to $-z_{w_t}$. If
neither case (i) or case (ii) applies, then we let $I(C) =C$. For
example, if $C$ is the element of $\mathcal{H}_{12}$ pictured in
Figure \ref{fig:fdes1}, then $I(C)$ is pictured in Figure
\ref{fig:fdes2}.

\fig{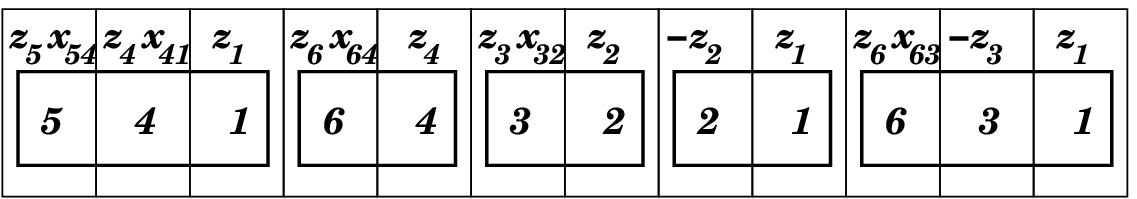}{$I(C)$ for $C$ in Figure \ref{fig:fdes1}.}

It is easy to see that $I^2(C) =C$ for all 
$C \in \mathcal{H}_{n}$ and that if 
$I(C) \neq C$, then $sgn(C)w(C) = -sgn(I(C))w(I(C))$. Hence  
$I$ is a weight-preserving and sign-reversing involution that shows 
\begin{equation}\label{des5}
\Gamma(h_n) = \sum_{C \in \mathcal{H}_n,I(C) = C}
sgn(C) w(C).
\end{equation}

Thus, we must examine the fixed points, $C = (B,w,L)$, of $I$.
First, there can be no $-z_a$ labels in $L$ so that $sgn(C) =1$.
Moreover,  if $b_j$ and $b_{j+1}$ are two consecutive bricks in $B$
and $t$ is the last cell of $b_j$, then it cannot be the case that
$w_t > w_{t+1}$ since otherwise we
could combine $b_j$ and $b_{j+1}$. Thus for each cell $t$ such 
that $w_t > w_{t+1}$, it must be the case that 
cells $t$ and $t+1$ lie in the same brick and, hence, 
cell $t$ is labeled with $z_{w_t}x_{w_t w_{t+1}}$.

It follows that $sgn(C)w(C) = \overline{z}^w \prod_{1 \leq i < j} 
x_{ji}^{\mathbf{ji}(w)}$.
For example, Figure \ref{fig:fdes3} shows a 
fixed point of $I$ in $H_{12}$. 

Vice versa, if
$w \in \PP^n$, then we can create a fixed point, $C
=(B,w,L)$, by having the bricks of $B$ end at cells $t$ such 
that either $w_t \leq w_{t+1}$ or $t=n$, labeling 
each cell $t$ such that $w_t > w_{t+1}$ with $z_{w_t}x_{w_t w_{t+1}}$ 
and labeling the remaining cells $t$ with $z_{w_t}$.
Thus we have shown that
\begin{equation*}
\Gamma(h_n) = \sum_{w \in \PP^n} \overline{z}^w \prod_{i <j} 
x_{ji}^{\mathbf{ji}(w)}.
\end{equation*}
as desired.

\fig{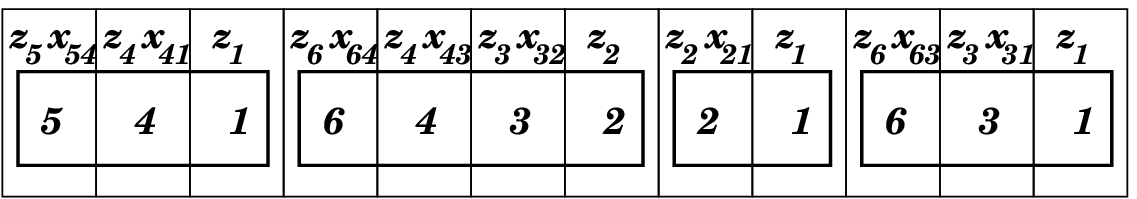}{A fixed point of $I$.}

Applying $\Gamma$ to the identity $H(t) = \frac{1}{E(-t)}$, we get
\begin{eqnarray*}
\sum_{n \geq 0} \Gamma(h_n) t^n &=& 1+ \sum_{n\geq 1} t^n  
\sum_{w \in \PP^n} \overline{z}^w \prod_{i <j} 
x_{ji}^{\mathbf{ji}(w)} \\
&=& \frac{1}{1+\sum_{n\geq 1} (-t)^n\Gamma(e_n)} \\
&=& \frac{1}{1+\sum_{n\geq 1}(-1)^{n} t^{n}
(-1)^{n-1} \sum_{S \subseteq \PP, |S|=n} DXZ(S)}\\
&=& \frac{1}{1-\sum_{n\geq 1} t^n
 \sum_{S \subseteq \PP, |S|=n} DXZ(S)} 
\end{eqnarray*}
which proves (\ref{eq:labdes}). \\
\ \\
{\bf Proof of Theorem \ref{thm:labwdes}.}\\

One can easily modify the proof of Theorem \ref{thm:labdes} 
to prove Theorem \ref{thm:labwdes}. 

 Recall that given a weakly decreasing word $w$ from 
$\PP^*$, we let 
\begin{equation}
WDXZ(w) = \begin{cases} 
z_j & \mbox{if} \ w = j, \ \mbox{and} \\
z_{j_1} \cdots z_{j_k} \prod_{i=1}^{k-1} (x_{j_{i}j_{i+1}}-1) & 
\mbox{if} \ w =j_1 \geq  \cdots \geq j_k \ \mbox{where} \ k \geq 2.
\end{cases}
\end{equation}

Define a ring homomorphism $\Gamma_w:\Lambda \rightarrow Q[\bf{x},\bf{z}]$ by 
defining $\Gamma(e_0) =1$ and, for $n \geq 1$,  
\begin{equation}\label{Gammawdef}
\Gamma_w(e_n) = (-1)^{n-1} \sum_{w \in WD\PP^*, |w| =n} WDXZ(w).
\end{equation}

Then we claim that 
\begin{equation}\label{wdes1}
\Gamma_w(h_n) = \sum_{w \in \PP^n} \overline{z}^w \prod_{i \leq j} 
x_{ji}^{\mathbf{ji}(w)}.
\end{equation}
That is,
\begin{eqnarray}\label{wdes2}
\Gamma_w(h_{n}) &=& 
\sum_{\mu\vdash n} (-1)^{n-\ell(\mu)}B_{\mu, (n)}\Gamma_w(e_{\mu}) \nonumber \\
&=& \sum_{\mu \vdash n} (-1)^{n-\ell(\mu)} \sum_{(b_1,
\ldots, b_{\ell(\mu)}) \in \mathcal{B}_{\mu,(n)}}
\prod_{j=1}^{\ell(\mu)} (-1)^{b_j-1} 
\sum_{w_j \subseteq WD\PP^*, |w_j| =b_j} WDXZ(w_j) \nonumber \\
&=& \ \sum_{\mu \vdash n} \sum_{(b_1, \ldots, b_{\ell(\mu)}) \in \mathcal{B}_{\mu,(n)}} \prod_{j=1}^{\ell(\mu)}  \sum_{w_j  \in WD\PP^*, |w_j| =b_j} WDXZ(w_j)
\end{eqnarray}

Next we want to give a combinatorial interpretation to
(\ref{wdes2}).  First we pick a brick tabloid $B = (b_1, \ldots, b_k)$ 
of length $n$. 
Then we interpret 
$\prod_{j=1}^{\ell(\mu)}  \sum_{w_j \subseteq  WD\PP^*, |w_j| =b_j} WDXZ(w_j)$ 
as picking a sequence of words in $WD\PP^*$, $(w_1, \ldots, w_{\ell(\mu)})$,  
such that $w_j$ has length $b_j$ and placing the elements of $w_j$ in the cells of $b_j$ in $j=1, \ldots, \ell(\mu)$. If 
$w_j = a_1 \geq   \cdots \geq a_{b_j}$, then we interpret the 
factor $WDXZ(w_j) = z_{a_1} \cdots z_{a_{b_j}} \prod_{i=1}^{b_j-1} (x_{a_i a_{i+1}}-1)$ as the ways of labeling the cells of $b_j$ that contain  
$a_i$ where $i <b_j$ with either 
$z_{a_i}x_{a_ia_{i+1}}$ or with $-z_{a_i}$ and labeling 
the last of cell $b_j$ with $z_{a_{b_j}}$.  
We shall call all such objects created in this way \emph{filled labeled
brick tabloids} and let $\mathcal{WDH}_{n}$ denote the set of all
filled labeled brick tabloids that arise in this way.  
Thus $\mathcal{WDH}_{n}$ consists of all triples $(B,w,L)$ such that 
\begin{enumerate}
\item $B =(b_1, \ldots, b_k)$ is a brick tabloid of length $n$, 
\item $w=w_1 \ldots w_n$ is a word in $\PP^n$ such that 
$w$ is weakly decreasing in each brick, and 
\item $L$ is a labeling of the cells of $B$ such that 
$L(i)$ is equal to $z_a$ if $i$ is the last cell of some brick 
$b_j$ which contains $a$ and $L(i) =-z_a$ or $L(i) = x_{ab}z_a$ if 
$i$ is not the last cell of a brick, cell $i$ contains $a$ and cell 
$i+1$ contains $b$. 
\end{enumerate} 
We then define the weight of $(B,w,L)$, $wt(B,w,L)$, to be the product of all
the $x_{ab}$ and $z_a$ labels in $L$ and the sign of 
$(B,w,L)$, $sgn(B,w,L)$, to be
the product of all the $-1$ factors in the labels in $L$. This 
process is illustrated in Figure \ref{fig:wfdes1} to construct 
an element $(B,w,L)$ of $H_{12}$ such that 
$$wt(B,w,L) = z_1 z_2^3 z_3 z_4^3 z_5 z_6^3 x_{54} x_{44} x_{64} x_{32} x_{66} $$
and $sgn(B,w,L) = -1$.

\fig{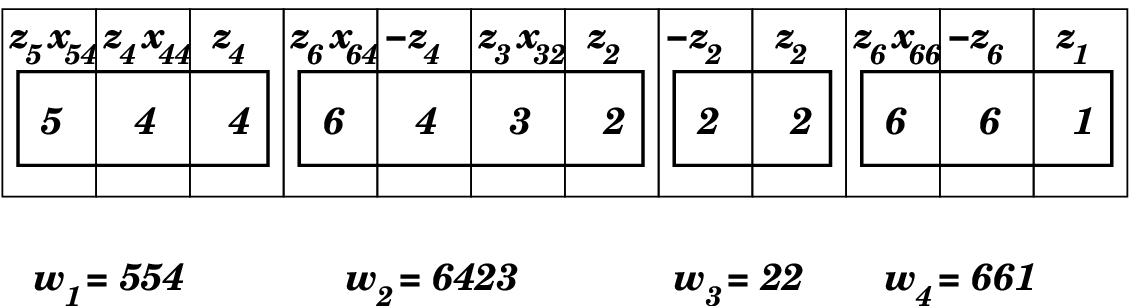}{A element $(B,w,L) \in \mathcal{WDH}_{12}$.}

At this point, the only difference in the proof is that we are dealing 
with filled brick tabloids which have weakly decreasing sequences in the 
 bricks  
rather than strictly decreasing sequences in the bricks.  This means that we can modify the involution $I$ of Theorem \ref{thm:labdes} by splitting 
bricks at cells labeled with $-z_i$ or combining two bricks 
such that the elements in the two bricks form a weakly decreasing sequence. 
Then essentially the same proof will show that (\ref{wdes1}) holds. \\
\ \\
{\bf Proof of Theorem \ref{thm:labris}.}\\

Given any weakly increasing word $w=w_1 \ldots w_n$, 
we let $S(w)$ denote the set of letters that appear 
in $W$.  For example, if 
$w =11123555$, then $S(w) = \{1,2,3,5\}$. We claim 
that for any non-empty set $S = \{j_1< \cdots < j_k\}$ 
contained in $\PP$, 
$$\sum_{w \in \PP^+,S(w) =S} t^{|w|} \overline{z}^w 
\prod_{i < j} x_{ij}^{\mathbf{ij}(w)} = t^{|S|} RXZ(S).$$
That is, if $S=\{j\}$, then $w$ must be of the 
for $j^k$ for some $k \geq 0$ so that in this 
case 
$$\sum_{w \in \PP^+,S(w) =S} t^{|w|} \overline{z}^w 
\prod_{i < j} x_{ij}^{\mathbf{ij}(w)}
= \frac{z_jt}{1-z_jt} = t^{|S|}\frac{z_j}{1-z_jt} =
t^{|S|}RXZ(S).$$
If $S(w) =\{j_1 < \cdots < j_k\}$ where $k \geq 2$, 
then $w$ must be of the form 
$w = j_1^{a_1} j_2^{a_2} \ldots j_k^{a_k}$ where $a_i \geq 1$ for 
$i=1, \ldots, k$. 
For any such word, it is easy to see 
that  
$$\prod_{i < j} x_{ij}^{\mathbf{ij}(w)} = 
\prod_{i=1}^{k-1} x_{j_i j_{i+1}}.$$
Hence, 
\begin{eqnarray*}
\sum_{w \in \PP^+,S(w) =S} t^{|w|} \overline{z}^w 
\prod_{i < j} x_{ij}^{\mathbf{ij}(w)}
&=& \left(\prod_{i=1}^k \frac{z_{j_i}t}{1-z_{j_i}t}\right)
\prod_{i=1}^{k-1} x_{j_i j_{i+1}}\\
&=& t^{|S|} \left(\prod_{i=1}^k\frac{z_{j_i}}{1-z_{j_i}t}\right)
\prod_{i=1}^{k-1} x_{j_i j_{i+1}}\\
&=& t^{|S|} RXZ(S).
\end{eqnarray*}
Thus 
$$\mathcal{R}(\mathbf{x}_\infty,\mathbf{z}_\infty,t) 
= 1+\sum_{n \geq 1} t^n \sum_{S \subseteq \PP,|S|=n} 
RXZ(S).$$
\ \\
{\bf Proof of Theorem \ref{thm:lablev}.}\\

Consider a factor $\left(1+\frac{z_it}{1-x_{ii}z_it}\right)$. One 
can think of the choice of 1 in that factor as not choosing $i$ to 
occur in the word where as the factor $\frac{z_it}{1-x_{ii}z_it}$ corresponds 
to choosing one of $i, ii, iii, iiii, \ldots $ in word.  Equation 
(\ref{eq:lablev}) easily follows.

\section{Possible extensions}

The methods that we have used in this paper can be modified 
to find generating functions of the form 
\begin{eqnarray*}
&& \sum_{w \in \PP^*,\umch(w)=0} t^{|w|} x^{\lev(w)+1}\overline{z}^w, \ \ \ \ \ \ \ \
\sum_{w \in [k]^*,\umch(w)=0} t^{|w|} x^{\lev(w)+1}\overline{z}^w,\\
&& \sum_{w \in \PP^*,\eumch(w)=0} t^{|w|} x^{\lev(w)+1}\overline{z}^w \ \mbox{and} \ \ 
\sum_{w \in [k]^*,\eumch(w)=0} t^{|w|} x^{\lev(w)+1}\overline{z}^w.
\end{eqnarray*}
in the case where $\lev(u)=1$ and generating functions of 
the form 
\begin{eqnarray*}
&& \sum_{w \in \PP^*,\umch(w)=0} t^{|w|} x^{\wdes(w)+1}\overline{z}^w, \ \ \ \ \ \ \ \ 
\sum_{w \in [k]^*,\umch(w)=0} t^{|w|} x^{\wdes(w)+1}\overline{z}^w, \\
&& \sum_{w \in \PP^*,\eumch(w)=0} t^{|w|} x^{\wdes(w)+1}\overline{z}^w \ \mbox{and} \ \ 
\sum_{w \in [k]^*,\eumch(w)=0} t^{|w|} x^{\wdes(w)+1}\overline{z}^w.
\end{eqnarray*}
in the case where $\wdes(u) =1$.
The idea is that one can modify the reciprocal method presented in Section 3 
to replace the statistic $\des(w)+1$ by $\lev(w)+1$ or $\wdes(w)+1$.  Then 
one can modify the collapse map appropriately. Finally, one 
needs appropriate modifications of Theorems \ref{thm:labdes}, \ref{thm:labwdes},
\ref{thm:labris}, and \ref{thm:lablev} to produce generating 
functions which keep track of labeled rises, levels, or descents that 
can be specialized to compute the generating functions of interest. 
For example, in the case where we study the distribution 
of $\lev(w)+1$ over words with no $u$-matches, the collapse map 
only produces words which have no consecutive repeated letters so 
that we need generating functions which keep track of labeled descents 
and rises over words which have no 11-match. 
These modifications will appear in the thesis of the second author. 

By the isomorphism which sends a word $w =w_1 \ldots w_n$ to its 
reverse, $w^r = w_n \ldots w_1$, one can automatically produce 
similar generating functions where the statistics $\des(w)+1$ and $\wdes(w)+1$ 
are replaced by  $\wris(w)+1$ and $\ris(w)+1$, respectively.

One can also easily modify the methods to keep track of restricted 
sets of descents.  For example, given a word $w=w_1 \ldots w_n \in \PP^*$, 
let $\mathrm{edes}(w) = |\{i:w_i > w_{i+1} \ \mbox{and} \ w_i \ \mbox{is even}\}|$. 
Then the techniques of this  paper can be easily modified to find 
closed expressions for 
\begin{eqnarray*}
&& \sum_{w \in \PP^*,\umch(w)=0} t^{|w|} x^{\mathrm{edes}(w)+1}\overline{z}^w, \ \ \ \ \ \ \ \
\sum_{w \in [k]^*,\umch(w)=0} t^{|w|} x^{\mathrm{edes}(w)+1}\overline{z}^w,\\
&& \sum_{w \in \PP^*,\eumch(w)=0} t^{|w|} x^{\mathrm{edes}(w)+1}\overline{z}^w \ \mbox{and} \ \ 
\sum_{w \in [k]^*,\eumch(w)=0} t^{|w|} x^{\mathrm{edes}(w)+1}\overline{z}^w
\end{eqnarray*}
in the case were $\mathrm{edes}(u) =1$.

Finally, one can extend the reciprocal methods in this paper 
to give a combinatorial interpretation of $U^{(\PP)}_{u,n}(x,{\bf z}_\infty,t)$, 
$U^{(k)}_{u,n}(x,{\bf z}_k,t)$, 
$EU^{(\PP)}_{u,n}(x,{\bf z}_\infty,t)$, and $EU^{(k)}_{u,n}(x,{\bf z}_k,t)$ 
in the case where $\des(u) > 1$.  Basically one has to modify 
the involution $I_u$ presented in Section 3 appropriately.  This has 
been done in the case of permutations by Quang Bach and the first 
author in the case of permutations \cite{Bach1,Bach2}.  However, 
in the case where $\des(u) >1$, the corresponding set of fixed points 
are much more complicated.  For example, it will no longer be the 
case that in fixed points of the modified version of $I_u$ that the underlying 
word will be weakly increasing in bricks. These more complicated fixed 
points then require a more complicated version of the collapse map. Nevertheless, 
one can still come up with closed formulas for the generating functions 
$U^{(\PP)}_{u}(x,{\bf z}_\infty,t)$, $U^{(k)}_{u}(x,{\bf z}_k,t)$, 
$EU^{(\PP)}_{u}(x,{\bf z}_\infty,t)$, and $EU^{(k)}_{u}(x,{\bf z}_k,t)$.
This work will appear in a subsequent paper.


\begin{thebibliography}{60}

\bibitem{AAM} R.E.L. Aldred, M. Atkinson, and D.J. McCaughan, 
Avoiding consecutive patterns in permutations, 
{\em Advances in Applied Mathematics}, {\bf 45} (2010), 449-461.


\bibitem{Bach1} Q. Bach and J.B. Remmel, Generating functions for permutations 
which avoid consecutive patterns with multiple descents, 
Australian Journal 
of Combinatorics. {\bf 64} (2016), 194-231. 


\bibitem{Bach2} Q. Bach and J.B. Remmel, Decent c-Wilf equivalence, 
arXiv:1520.07190. 

\bibitem{Bec4}
D. Beck and J. Remmel, Permutation enumeration of the symmetric
  group and the combinatorics of symmetric functions, {\em J. Combin. Theory Ser. A}, 
  \textbf{72} (1995), no.~1, 1--49.

\bibitem{Bre1}
F. Brenti,Unimodal polynomials arising from symmetric functions,
 {\em  Proc. Amer. Math. Soc.}, \textbf{108} (1990), no.~4, 1133--1141.

\bibitem{Bre}
F. Brenti, Permutation enumeration symmetric functions, and
  unimodality, {\em Pacific J. Math.}, \textbf{157} (1993), no.~1, 1--28.


\bibitem{DuaneR}
A. Duane and J. Remmel, Minimal overlapping patterns in colored 
permutations, Electronic J. Combinatorics, {\bf 18(2)} (2011), P25, 34 pgs.


\bibitem{Eliz}
S. Elizalde and M. Noy, Consecutive patterns in permutations, {\em Adv.
  in Appl. Math.}, \textbf{30} (2003), no.~1-2, 110--125, 
Formal power series and
  algebraic combinatorics (Scottsdale, AZ, 2001).

\bibitem{Eliz2} S. Elizalde, The most and least avoided consecutive pattern, 
Proc. Lond. Math. Soc., \textbf{106} (2013), 957-979.

\bibitem{Eg1}
O. E{\u{g}}ecio{\u{g}}lu and J. B. Remmel, Brick tabloids and the
  connection matrices between bases of symmetric functions, {\em Discrete Appl.
  Math.}, \textbf{34} (1991), no.~1-3, 107--120, Combinatorics and theoretical
  computer science (Washington, DC, 1989).





\bibitem{GJ} I.P. Goulden and D.M. Jackson, {\em Combinatorial 
Enumeration}, John Wiley \& Sons 1983.


\bibitem{HM} S. Heubach and T. Mansour, {\em Combinatorics 
of Compositions and Words}, Discrete Mathematics and Its Applications, 
Chapman \& Hall/CRC, Taylor 
\& Francis Group, Boca Raton, London, New York, (2009).


\bibitem{JR} M. Jones and J.B. Remmel, Pattern Matching in the Cycle Structures of Permutations, {\em Pure Math. and Applications}, {\bf 22} (2011), 
173-208.


\bibitem{JR11} M. Jones and J.B. Remmel, 
A reciprocity approach to computing 
generating functions for permutations with no pattern matches, 
{\em Discrete Mathematics and Theoretical Computer Science}, 
DMTCS Proceedings, 23 International Conference on Formal 
Power Series and Algebraic Combinatorics (FPSAC 2011), {\bf 119} (2011), 
551-562.

\bibitem{JR1} M. Jones and J.B. Remmel, 
A reciprocity method for computing generating functions over the set of permutations with no consecutive occurrence of $\tau$, Discrete Mathematics, {\bf 313} Issue 23 (2013), 2712-2729. 


\bibitem{JR2} M. Jones and J. Remmel, Generating functions 
for the number of permutations with no 
consecutive occurrences of  $1p23 \ldots (p-1)$ or $13\ldots (p-1)2p$,  
to appear in Pure Mathematics and Applications.




\bibitem{Kit} S. Kitaev, {\em Patterns in permutations and words}, Springer-Verlag, 2011.


\bibitem{KS} A. Khoroshkin and B. Shapiro, Using homological duality 
in consecutive pattern avoidance, {\em Electronic. J. Combinatorics}, {\bf 18} (2011), \# P9



\bibitem{Mac}
I. G. Macdonald, \emph{Symmetric Functions and Hall Polynomials.}
2nd ed. Oxford University Press, 1995.



\bibitem{LieseRem} 
J. Liese and J.B. Remmel, Generating functions for permutations 
avoiding a consecutive pattern, {\em Annals of Combinatorics}, 
{\bf 14} (2010), 103-121. 


\bibitem{Lan2}
T. Langley and J.B. Remmel, Enumeration of \begin{math}m\end{math}-tuples of
  permutations and a new class of power bases for the space of symmetric
  functions, {\em Adv. Appl. Math.}, {\bf 36} (2006), 30-66.


\bibitem{MenRem1} A. Mendes and J.B. Remmel, 
Generating functions for statistics on 
\begin{math}C_k \wr S_n\end{math}, {\em S\'eminaire Lotharingien de Combinatoire}, B54At, (2006), 40 pp.

\bibitem{MenRem2}  A. Mendes and J.B. Remmel, 
Permutations and words counted by 
consecutive patterns, {\em Adv. Appl. Math.}, {\bf 37} 4, (2006) 443-480.

\bibitem{MenRem3} A. Mendes and J.B. Remmel, Descents, major indices, and inversions in permutation groups, {\em Discrete Mathematics}, Vol. 308, Issue 12, (2008), 2509-2524. 





\bibitem{MRR}  A. Mendes, J.B. Remmel, and A. Riehl, 
Permutations with \begin{math}k\end{math}-regular  descent patterns, 
\emph{ Permutation Patterns} (S. Linton, N. Ruskuc, 
and V. Vatter, eds.), London Math. Soc. Lecture Notes 376, 259-286, (2010).

\bibitem{MenRembook}
 A. Mendes and J. Remmel, Counting with symmetric functions, 
Developments in Mathematics, vol. 43, Springer (Cham, Heidelberg, New York, Dordrecht, London), ISBN 978-3-309-23617-9 (ISBN 978-3-23619-6 Electronic).  


\bibitem{RT} D. Rawlings and M. Tiefenbruck, Consecutive Patterns: From Permutations to Column-Convex Polyominoes and Back, \textit{Electronic Journal of Combinatorics}, {\bf 17 (1)} (2010), \#R62


\bibitem{Stan} R.P. Stanley, Enumerative Combinatorics, vol. 2, Cambridge University Press, (1999).



\bibitem{St} E. Steingr\'{\i}msson: Generalized permutation patterns -- a short
survey,
\emph{ Permutation Patterns, St Andrews 2007}, S.A. Linton, N.
Ruskuc, V. Vatter (eds.), LMS Lecture Note Series 376, Cambridge
University Press, (2010), 137-152.

\end{thebibliography}
\end{document}